\documentclass[a4paper,final]{amsart} 
\usepackage{amsmath}
\usepackage{epsfig}
\usepackage{amssymb} 
\theoremstyle{remark}

\newtheorem{remark}{Remark}

\theoremstyle{definition}

\theoremstyle{theorem}

\newtheorem*{acknow}{Acknowledgements}
\newcommand{\R}{\mathbb{R}}

\newcommand{\T}{\mathbb{T}}

\begin{document}

% Definition of title page:
\title[Numerical study of KP equations]{Numerical study of blow up and stability of solutions 
of generalized Kadomtsev-Petviashvili equations} 

\author{C.~Klein}
\address{Institut de Math\'ematiques de Bourgogne,
                Universit\'e de Bourgogne, 9 avenue Alain Savary, 21078 Dijon
                Cedex, France}
    \email{Christian.Klein@u-bourgogne.fr}

\author{J.-C.~Saut}
    \address{Laboratoire de Math\' ematiques, UMR 8628,
    Universit\' e Paris-Sud et CNRS,
    91405 Orsay, France}
    \email{Jean-Claude.Saut@math.u-psud.fr}

\date{\today}    % optional

\begin{abstract}
We first review the known mathematical results concerning the KP type equations. Then we perform numerical simulations to analyze various qualitative properties of the equations : blow-up versus long time behavior, stability and instability of solitary waves.
\end{abstract}


%\thanks{We thank M.~Haragus for useful discussions and hints, L. Molinet and N. Tzvetkov for a lengthy and fruitful joint research with the second Author on KP equations. 
%This work has been supported in part by the project FroM-PDE funded by the European
%Research Council through the Advanced Investigator Grant Scheme, the Conseil R\'egional de Bourgogne
%via a FABER grant, the ANR via the program ANR-09-BLAN-0117-01 and the Wolfgang Pauli Institute in Vienna. }

\maketitle
\section{Introduction}
The main goal of this paper is an  analysis of qualitative properties of solutions to generalized Kadomtsev-Petviashvili (KP) equations via a numerical approach.

The  (classical) Kadomtsev-Petviashvili (KP) equations

\begin{equation}\label{KP}
(u_t+u_{xxx} +u u_x )_x \pm u_{yy} =0
\end{equation}

were introduced in \cite{KaPe} to study the transverse stability of the solitary wave solution of the 
Korteweg- de Vries equation which reads in the context of water-waves

\begin{equation}\label{KdV}
u_t+u_x+uu_x+(\frac{1}{3}-T)u_{xxx}=0, \;x\in \R,\;t\geq0.
\end{equation}

Here $T\geq 0$ is the Bond number\footnote {Note that here $T=0$ corresponds to the absence of surface tension. In the Fluid Mechanics community, the Bond number is often defined as the inverse of our $T$.}  which measures  surface tension effects in the context of surface hydrodynamical waves.

Actually the (formal) analysis in \cite{KaPe} consists in looking for a {\it  weakly transverse} perturbation 
of the one-dimensional transport equation

\begin{equation}\label{transp}
u_t+u_x=0.
\end{equation}

This perturbation, which is obtained by a Taylor expansion of the dispersion relation  $\omega(k_1,k_2)=\sqrt{k_1^2+k_2^2}$ of the two-dimensional linear wave equation assuming $|k_1|\ll 1 $ and $\frac{|k_2|}{|k_1|}\ll 1,$ amounts to adding a nonlocal term, leading to 

\begin{equation}\label {perttransp}
u_t+u_x+\frac{1}{2}\partial_x^{-1}u_{yy}=0.
\end{equation}

Here the operator $\partial_x^{-1}$ is defined via Fourier transform,

$$\widehat{\partial_x^{-1}f}(\xi)=\frac{i}{\xi_1}\widehat{f}(\xi),\,\text{where}\;\xi=(\xi_1,\xi_2).$$

\begin{remark}\label{LDPE}
Equation \eqref{perttransp} is reminiscent of the {\it linear diffractive pulse equation}

$$2u_{tx}=\Delta_y u,$$

where $\Delta_y$ is the Laplace operator in the transverse variable $y,$
studied in \cite{AR}.
\end{remark}

The same formal procedure is applied to the KdV equation \eqref{KdV}, 
assuming that the transverse dispersive effects are of the same order as the x-dispersive and nonlinear terms,   yielding the KP equation in the form

\begin{equation}\label{KPbrut}
u_t+u_x+uu_x+(\frac{1}{3}-T)u_{xxx}+\frac{1}{2}\partial_x^{-1}u_{yy}=0.
\end{equation}

By change of frame and scaling, \eqref{KPbrut} reduces to \eqref{KP} with the $+$ sign (KP II) when $T<\frac{1}{3}$ and the $-$ sign (KP I) when $T>\frac{1}{3}$.

Note however that $T>\frac{1}{3}$ corresponds to a layer of fluid of depth smaller than $0.46$ cm, and in this situation viscous effects due to the boundary layer at the bottom cannot be ignored. One could then say that ``the KP I equation does not exist in the context of water waves", but it appears naturally in other contexts (see for instance Remark 2 below).

\vspace{0.3cm}
Of course the same formal procedure could also be applied to {\it any} one-dimensional weakly nonlinear dispersive  equation of the form

\begin{equation}\label {gen}
u_t+u_x+f(u)_x-Lu_x=0, \;x\in \R,\;t\geq0,
\end{equation}

where $f(u)$ is a smooth real-valued function (most of the time polynomial) and $L$ a linear operator taking into account the dispersion and defined in Fourier variable by

\begin{equation}\label {L}
\widehat{Lu)}(\xi)=p(\xi)\mathfrak F u(\xi),
\end{equation}

where the symbol $p(\xi)$ is real-valued. The KdV equation corresponds for instance to $f(u)=\frac{1}{2}u^2$ and $p(\xi)= -\xi^2.$  Examples with a fifth order dispersion in $x$ are considered in \cite{AbSt}, \cite{KaBe1}, \cite{KaBe2}.

This leads to a class of generalized KP equations

\begin{equation}\label {genKP}
u_t+u_x+f(u)_x-Lu_x+\frac{1}{2}\partial_x^{-1}u_{yy}=0, \;x\in \R,\;t\geq0.
\end{equation}

Let us notice, at this point, that models alternative to KdV-type equations
\eqref{gen} are the equations of Benjamin--Bona--Mahony (BBM) type \cite{BBM}
\begin{equation}\label{1.2}
u_{t}+u_{x}+uu_{x}+Lu_{t}=0
\end{equation}
with corresponding two-dimensional ``KP--BBM-type models'' (in the case $p(\xi)\geq 0$)
\begin{equation}\label{1.5}
u_{t}+u_{x}+uu_{x}+Lu_{t}+\partial_{x}^{-1}\partial_{y}^{2}u=0
\end{equation}
or, in the {\it derivative form}
\begin{equation}\label{1.10}
(u_{t}+u_x+uu_{x}+Lu_{t})_{x}+ \partial_{y}^{2}u=0
\end{equation}
and free evolution group
\[
S(t)=e^{-t(I+L)^{-1}[\partial_{x}+ \partial_{x}^{-1}\partial_{y}^{2}]}\,\, .
\]

It was only after the seminal paper \cite{KaPe} that Kadomtsev-Petviashvili type equations have been  derived as asymptotic  models to describe the propagation of long, quasi-unidirectional waves of small amplitude (under an appropriate scaling) in various physical situations. This was done for instance formally in \cite{AS}  in the context of water waves (see  \cite{L}, \cite {LS} for a rigorous approach in the same context).

A rigorous approach (with error estimates)  was also made by Ben-Youssef and Lannes \cite{BY-L} for a class of hyperbolic quasilinear systems. We refer also to \cite {GS}, \cite{Pau} for a rigorous derivation from  a Boussinesq system or the Benney-Luke equation. In all these works, the convergence rate is proven to be low (we will be a bit more precise in Subsection 2.1), a phenomena which is deeply analysed in \cite{L} (see also \cite{LS}) and which is mainly due to the singularity at $\xi_1=0$ of the  KP dispersion relation. 

On the other hand, Dryuma \cite{Dr} found a Lax pair to the KP- I/II equations, proving the ``integrability" of the KP equations (see \cite{NMPZ} for a precise description of the Inverse Scattering aspects of the KP equations).
\begin{remark}
In some physical contexts (not in water waves!) one could consider higher dimensional transverse perturbations, which amounts to replacing $\partial_x^{-1}u_{yy}$ in \eqref{genKP} by $\partial_x^{-1}\Delta^{\perp}u$, where $\Delta^{\perp}$ is the Laplace operator in the transverse variables.
%\end{remark}

For instance, the KP I equation (in both two and three dimensions) also describes after a suitable scaling the long wave {\it transonic} limit of the Gross-Pitaevskii equation (see \cite {BGS} for rigorous results  for the solitary waves (ground states)  in $2D$   and \cite{CR} for the Cauchy problem in dimensions two and more).

More precisely, let the Gross-Pitaevskii equation be

\begin{equation}\label{GP}
i\psi_t+ \Delta \psi+(1-|\psi|^2)\psi=0, \;\text{in}\;\R^d\times \R,
\end{equation}

for functions $\psi$ with finite {\it Ginzburg-Landau energy}

\begin{equation}\label{EGP}
E_{GL}(t)=\frac{1}{2}\int_{\R^d}\lbrack|\nabla \psi|^2+\frac{1}{2}(1-|\psi|^2)^2 \rbrack dx.
\end{equation}

The Gross-Pitaevskii equation can be written in hydrodynamic form via the Madelung transform provided $\psi$ does not vanish. This is the case in the {\it transonic regime} where in particular  $|\psi|$ is close to one. In this regime, the KP I equation describes, after a suitable rescaling, the behavior of $\eta =1-|\psi|^2.$ 

\end{remark}

Note again that in the classical KP equations, the distinction  between KP I and KP II arises from the {\it sign} of the dispersive term in $x$, which characterizes what physicists call ``positive" or ``negative" dispersion media.

\vspace{0.3cm}
A noteworthy fact implied by the derivation of KP equations is that, as far as the {\it transverse stability} of the KdV {\it  solitary wave  ``1-soliton"} 
\begin{equation}
    \phi_c(x,y)=\frac{3c}{2}\,{\rm cosh}^{-2}\Big(\frac{\sqrt{c}\, x}{2}\Big)
    \label{kdvsol}
\end{equation}
is concerned, the natural initial condition 
associated to \eqref {KP} should be of the type $u_0=\phi_c+v_0$ where $v_0$ is either ``localized" in $x$ and $y$, or localized in $x$ and $y$-periodic.

The same observation is of course valid for \eqref{genKP} in connection with the transverse stability of solitary waves of \eqref{gen}.

This paper is organized as follows. In the next Section we review  the main rigorous mathematical results concerning KP type equations. They are mainly obtained by PDE techniques with the exception of results obtained via the Inverse Scattering machinery for the {\it classical} KP I and KP II equations.  This Section has an interest in itself but it  leads also to {\it conjectures} on the long time dynamics (or blow-up in finite time when it is expected) for KP type equations. To keep this survey short we will only discuss the main results and ideas and refer to the original papers for details and proofs. Also we will mainly comment on results obtained by PDE or Nonlinear Analysis techniques. We refer for instance to the excellent survey article \cite{Fo} which describes in particular  results on the Cauchy problem obtained by Inverse Scattering methods for the classical KP I and KP II equations. 

The aim of the next Section 3 which is the heart of this work  is to give numerical evidence of those conjectures and to suggest further theoretical investigations. We will also provide precise numerical decay rates of various norms of the solutions.

\vspace{0.5cm}
\centerline{Notations}

For any real number $s_0$, the notation $s_0^+$ means ``for any $s>s_0$". The norm in the classical Lebesgue spaces $L^p(\R^d)$ will be denoted by $|.|_p.$

We denote, for any real number $s,$  by $H^s(\R^d)$ the Sobolev space of distributions  $u$ in the Schwartz space ot tempered distributions $\mathcal S'(\R^d)$ equipped with the norm $||u||_s=|(1+|\xi|^2)^{\frac{s}{2}}\hat{u}|_2$ where $\hat{u}$ denotes the Fourier transform of $u$.

We will use sometimes the anisotropic Sobolev spaces $H^{s_1,s_2}(\R^2)$, equipped with the norm

$||u||_{s_1,s_2}=|(1+|\xi_1|^2)^{\frac{s_1}{2}}(1+|\xi_2|^2)^{\frac{s_2}{2}} \hat{u}|_2$. 

For the {\it homogeneous} version $\dot{H}^s$ of those spaces, we just replace the weight $(1+|\xi|^2)^{\frac{s}{2}}$ by $|\xi|^s$.

\section{A survey of theoretical results on KP equations}

We survey here various mathematical results concerning KP type equations. They will serve as a guide for our numerical simulations, together with the (many) open problems.

\subsection{The zero-mass constraint in x}
\label{constraint}

In \eqref{genKP}, it is implicitly assumed that the operator $\partial_{x}^{-1}\partial_{y}^{2}$
is well defined, which a~priori imposes a constraint on the solution $u$, which, in some sense, has to be an $x$-derivative.
This is achieved, for instance, if  $u\in {\mathcal S}'(\R^2)$ is such that
\begin{equation}\label{1.6}
\xi_{1}^{-1}\,\xi_{2}^{2}\,\widehat{u}(t,\xi_1,\xi_2)\in {\mathcal S}'(\R^2)\, ,
\end{equation}
thus in particular if $\xi_{1}^{-1}\,\widehat{u}(t,\xi_1,\xi_2)\in {\mathcal S}'(\R^2)$.
Another possibility to fulfill the constraint is to write $u$ as
\begin{equation}\label{1.7}
u(t,x,y)=\frac{\partial}{\partial x}v(t,x,y),
\end{equation}
where $v$ is a continuous function having a classical derivative with respect
to $x$, which, for any fixed $y$ and $t\neq 0$, vanishes when $x\rightarrow \pm \infty$.
Thus one has
\begin{equation}\label{1.8}
\int_{-\infty}^{\infty}u(t,x,y)dx=0,\qquad y\in\R,\,\,\, t\neq 0,
\end{equation}
in the sense of generalized Riemann integrals.
Of course the differentiated version of \eqref{genKP}, namely
\begin{equation}\label{1.9}
(u_{t}+u_{x}+uu_{x}-Lu_{x})_{x}+ \partial_{y}^{2}u=0,
\end{equation}
can make sense without any constraint of type (\ref{1.6}) or  (\ref{1.8}) on $u$, 
and so does the Duhamel integral representation of \eqref{genKP},
\begin{equation}\label{1.11}
u(t)=S(t)u_{0}-\int_{0}^{t}S(t-s)(u(s)u_{x}(s))ds,
\end{equation}
where $S(t)$ denotes the (unitary in all Sobolev spaces $H^{s}(\R^2)$) group associated with \eqref{genKP},
\begin{equation}\label{1.12}
S(t)=e^{-t(\partial_{x}-L\partial_{x}+\partial_{x}^{-1}\partial_{y}^{2})}\,\, .
\end{equation}

In particular, the results established on the Cauchy problem for KP type equations  which use the Duhamel (integral ) formulation (see for instance \cite{Bou} \cite{ST4})  are valid {\it without} any constraint on the initial data.

 One has to be careful however {\it in which sense} the {\it differentiated equation} is taken. For instance let us consider the integral equation
 
 \begin{equation}\label {dudu}
 u(x,y,t)=S(t) u_0(x,y)-\int_0^t S(t-t')\lbrack u(x,y,t')u_x(x,y,t)\rbrack dt',
 \end{equation}
 
where $S(t)$ is here the KP II group, for initial data $u_0$ in $H^s(\R^2),$ $s>2$, (thus $u_0$ does not satisfy any zero mass constraint), yielding a local solution $u\in C(\lbrack 0,T\rbrack; H^s(\R^2))$.

By differentiating \eqref{dudu} first with respect to $x$ and then with respect to $t,$ one obtains the equation

$$\partial _t\partial_xu+\partial_x(uu_x)+\partial_x^4u+\partial_y^2u =0 \quad\text{in} \quad C(\lbrack 0,T\rbrack; H^{s-4}(\R^2)).$$

However, the identity $ \partial_t\partial_xu =\partial_x\partial_t u$ holds only in a very weak sense, for example in $\mathfrak D'((0,T)\times \R^2).$

On the other hand, a constraint has to be imposed when using the Hamiltonian
formulation of the equation. In fact, the Hamiltonian for \eqref{1.9} is
\begin{equation}\label{1.13}
\frac{1}{2}\int\left[-u\,Lu+(\partial_{x}^{-1}u_y)^2+u^2+\frac{u^3}{3}\right]
\end{equation} 
and the  Hamiltonian associated with (\ref{1.10}) is
\begin{equation}\label{1.14}
\frac{1}{2}\int\left[(\partial_{x}^{-1}u_y)^2+u^2+\frac{u^3}{3}\right]\, .
\end{equation}

It has been established in \cite{MST1} that, for a rather  general class of KP or KP--BBM
equations,  the solution of the Cauchy problem obtained for 
\eqref{1.9}, \eqref{1.10} (in an appropriate functional setting) satisfies the zero-mass constraint in $x$  for any $t \neq 0$ (in a sense to be precised below), even if the initial data does not. This is a manifestation of the infinite speed of propagation inherent to KP equations. Moreover, KP type equations display a striking {\it smoothing effect} (different from the one reviewed in Subsection 2.3 though) : if the initial data belongs to the space $L^1(\R^2)\cap H^{2,0}(\R^2)$ and if $u\in C(\lbrack 0,T\rbrack; H^{2,0}(\R^2))$  \footnote {We will see in Subsection 2.4  that KP type equations (in particular the classical KP I and KP II equations) do possess solutions in this class.} is a solution in the sense of distributions, then, for any $t>0,$  $u(.,t)$ becomes a {\it continuous} function of $x$ and $y$ (with zero mean in $x$).  Note that the space  $L^1(\R^2)\cap H^{2,0}(\R^2)$ is not included in the space of continuous functions.

The key point when proving those results is a careful analysis of  the fundamental solution 
of KP-type equations \footnote{ In the case of KP II, one can use the explicit form of the fundamental solution found in \cite{Re}.} which turns out to be  a $x$ derivative of a continuous function of $x$ and $y$, $C^1$ with respect to  $x,$ which, for fixed $t\neq 0$ and $y$, tends to zero as
$x\rightarrow \pm \infty$. Thus its generalized Riemann integral in $x$
vanishes for all values of the transverse variable $y$ and of $t\neq 0$. A
similar property can be established for the solution of the nonlinear
problem \cite{MST1}. Those results have been checked in the numerical simulations of \cite{klsp}, \cite{klspma}. 

 We refer  \cite{FoSu1}, \cite{Su}  for a   rigorous approach of the Cauchy problem with initial data which do not satisfy the zero-mass condition via the Inverse Spectral Method  in the integrable case.

Nevertheless, the singularity at $\xi_1=0$ of the dispersion relation of KP type  equations make them rather {\it bad} asymptotic models. First the singularity at $\xi_1=0$ yields a very bad approximation of the dispersion relation of the original system (for instance the water wave system) by that of the KP equation. 

Another drawback is the poor error estimate between the KP solution and the solution of the original problem. This has been established clearly in the context of water waves as we briefly recall now.

Let us first recall the  KP scaling for water waves. We denote by 
$a$  a  typical amplitude, $h$ the mean depth, $\lambda$ (resp. 
$\mu$) the  typical wavelength in $x$, (resp. $y$). Then we  set

$$\frac{a}{h}= \epsilon,\; \frac{\lambda^2}{h^2}=\frac{S_1}{\epsilon}, \; \frac{\mu^2}{h^2}=\frac{S_2}{\epsilon^2},$$ 

where $\epsilon \ll1,$ $S_1\sim 1,$ $S_2\sim 1.$\footnote { $\epsilon$ is thus the small parameter which measures the weak nonlinear and long wave effects.}

It has been established rigorously in \cite{L}, \cite{LS} that the KP 
II equation yields a poor error estimate when used as an asymptotic 
model of the water wave system.  Roughly speaking,  the error 
estimates with the relevant solution of the (Euler) water waves 
system reads: 

$$||U_{Euler}-U_{KP}||=o(1).$$

So the error is $o(1),$ ($O(\sqrt {\epsilon})$ with some additional 
constraint) instead of  $O(\epsilon ^2 t), $ which should be the optimal rate in this regime (as it is the case for the KdV, Boussinesq equations or systems).  Nevertheless the KP II equation  reproduces (qualitatively) observed features of the water wave theory. For instance  the well-known picture below displays the interaction of two oblique ``line solitary waves" on the Oregon coast which shows a striking resemblance with the so-called  KP II  2-soliton.
 
\begin{figure}[htbp]
   \begin{center}
\includegraphics[width=90mm]{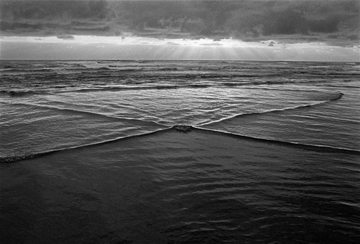}
   \end{center}
   \caption{\footnotesize Interaction of line solitons. Oregon coast}
\end{figure}

We refer to \cite{LS} for the derivation of a new class of weakly transverse approximations of the water waves system in the Boussinesq regime which are purely local (without the anti-derivative in $x$), have a dispersion relation which fitting well with that of the original Euler system,  and which provide {\it optimal} error estimates. No mass constraint has to be imposed there.

\subsection{Localized solitary waves}

Solitary waves are solutions of KP equations of the form

$$u(x,y^{\perp}, t)= \psi_c(x-ct, y^{\perp}),$$

where $y^{\perp}$ ($=y$ or $(y,z))$ is the transverse variable and $c>0$ is  the solitary wave velocity.

The solitary wave is said to be {\it localized} if $\psi_c$ tends to zero at infinity in all directions. For such solitary waves, the following {\it energy space} is natural

$$Y=Y(\R^2)=\lbrace u\in L^2(\R^2),\; u_x \in L^2(\R^2),\;\partial _x^{-1}u_y \in L^2(\R^2)\rbrace,$$

(with an obvious modification in the three dimensional case),

and throughout this section we will deal only with {\it finite energy solitary waves}.

Due to its integrability properties, the KP I equation possesses a localized, finite energy, explicit solitary wave, the {\it lump} \cite{manakov} :

\begin{equation}\label{Lump}
\phi_c(x-ct,y)=\frac{8c(1-\frac{c}{3}(x-ct)^2+\frac{c^2}{3}y^2)}{(1+\frac{c}{3}(x-ct)^2+\frac{c^2}{3}y^2)^2}.
\end{equation}

Another interesting explicit solitary wave of the KP I equation which is {\it localized in $x$ and periodic in $y$} has been found by Zaitsev \cite{Z}. It reads

\begin{equation}\label{Za}
Z_c(x, y)=12\alpha^2\frac{1-\beta \cosh (\alpha x)\cos (\delta y)}{\lbrack \cosh (\alpha x)-\beta\cos (\delta y)\rbrack^2},
\end{equation}

where 

$$(\alpha, \beta)\in ( 0, \infty)\times (-1, +1),$$

and the propagation speed is given by

$$c=\alpha ^2\frac{4-\beta ^2}{1-\beta^2}.$$

Let us observe that the transform $\alpha\rightarrow i\alpha$,
$\delta\rightarrow i\delta$, $c\rightarrow ic$ produces solutions of
the KP I equation which are periodic in $x$ and localized in $y$.

\vspace{0.3cm}
We address first the question of {\it non existence} of solitary waves for the following classes of KP type equations

\begin{equation}\label{GKP2D}
  u_t+u^pu_x+u_{xxx}+\epsilon \partial_x^{-1}u_{yy}=0,
 \end{equation}
 
 where $\epsilon =\pm1$, in the two-dimensional case, and
 
 \begin{equation}\label{GKP3D}
  u_t+u^pu_x+u_{xxx}+a \partial_x^{-1}u_{yy} +b\partial_x^{-1}u_{zz} =0,
 \end{equation}
 
 where $a, b =\pm1$, in the three-dimensional case.

 By establishing suitable Pohojaev type identities, it has been proved in \cite{deBS} that (under mild regularity conditions), no localized solitary waves exist when $c\leq 0$ in both dimensions and for $c>0$ when 
 
 $$ \;\epsilon =-1,\; p\geq 4\quad\text{or}\quad \epsilon =1\;\text{and}\; p\;\text{arbitrary},$$
 
 $$ab=-1\;(\text{resp.}\;a=b=1)\;\text{and}\; p\;\text{arbitrary},$$
 
 \begin{center}
  or
  \end{center}
  
  $$a=b=-1,\;p\geq \frac{4}{3}.$$
  
  In particular, the 2D  KP II-type equations ($\epsilon =1,$ $p\geq1$) {\it do not} possess any localized solitary waves. This was precised by de Bouard and Martel \cite{BoMa} who proved that the classical KP II equation does not possess $L^2$- compact solutions (those finite energy solutions are more general than solitary waves).
  
  We refer to \cite{deBS} for nonexistence results of localized solitary waves of KP type equations with a cubic-quintic dispersion in $x$.

 \vspace{0.3cm}
 We turn now to the {\it existence} of localized, finite energy solitary waves. We will consider as a basic example a class of generalized KP I equations in two dimensions of the form

 \begin{equation}\label{lump}
 u_t+u^pu_x+u_{xxx}-\partial_x^{-1}u_{yy}=0,
 \end{equation}
 
 where $p=\frac{m}{n},$\; $n$ odd, $m$ and $n$ relatively prime.

 Solitary waves are looked for in the energy space $Y(\R^2)$ which can also be defined (see \cite{deBS})  as the closure of the space of $x$ derivatives of smooth and compactly supported functions in $\R^2$ for the norm
 
$$\| \partial_x f \|_{Y(\R^2)} \equiv \Big( \| \nabla f \|_{L^2(\R^2)}^2 + \| \partial_x^2 f \|_{L^2(\R^2)}^2 \Big)^\frac{1}{2}.$$

The equation of a solitary wave $\psi$ of speed $c$ is given by

\begin{equation}
\label{SW}
c\partial_x \psi - \psi \partial_x \psi - \partial_x^3 \psi + \partial_x^{- 1} (\partial_y^2 \psi) = 0,
\end{equation}

which implies 
\begin{equation}
\label{SW2}
c\partial_{xx} \psi - (\psi \partial_x \psi)_x - \partial_x^4 \psi +  \partial_y^2 \psi = 0,
\end{equation}

When $\psi \in Y(\R^2)$, the function $\partial_x^{- 1} \partial_y ^2\psi$ is well-defined (see \cite{deBS}), so that \eqref{SW} makes sense.

Given any $c > 0$, a solitary wave $\psi_c$ of speed $c$ is deduced from a solitary wave  $\psi_1$ of velocity $1$ by the scaling
\begin{equation}
\label{scalingsw}
\psi_c(x, y) = c \psi_1(\sqrt{c} x, c y).
\end{equation}

\vspace{0.3cm}
We now introduce the important notion of {\it ground state} solitary waves.
% Solitary waves may be obtained in dimension two minimizing the Hamiltonian keeping the $L^2$-norm fixed (see \cite{deBoSau3,deBoSau1}). Like \eqref{GP}, equation \eqref{KP} is 

We set

$$E_{KP}(\psi) = \frac{1}{2} \int_{\R^2} (\partial_x \psi)^2 + \frac{1}{2} \int_{\R^2} (\partial_x^{-1}\partial_y \psi)^2 - \frac{1}{2(p+2)} \int_{\R^2} \psi^{p+2},$$
%and the $L^2$-norm of $\psi$ is conserved as well. Setting
and we define the action

$$S(N) = E_{KP}(N) + \frac{c}{2} \int_{\R^2} N^2.$$

We term {\it ground state}, a solitary wave $N$ which minimizes the action $S$ among all finite energy non-constant solitary waves of speed $c$ (see \cite{deBS} for more details). 
%In dimension two, a solitary wave is a ground state if and only if %it minimizes the Hamiltonian $E_{KP}$ keeping the $L^2$-norm fixed (see \cite{deBoSau3}). The constant $\boS_{KP}$, which appears in Theorem \ref{dim2}, denotes the action $S(N)$ of the ground states $N$ of speed $\sigma = 1$.

It was proven in \cite{deBS} that ground states exist if and only if $c>0$ and $1\leq p<4$. Moreover, when $1\leq p<\frac{4}{3},$ the ground states are minimizers of the Hamiltonian  $E_{KP}$ with prescribed mass ($L^2$ norm).

\begin{remark}
When $p=1$ (the classical KP I equation), it is unknown (but conjectured) whether the lump solution is a ground state.
\end{remark}

In the three dimensional case, ground states exist when $a=b=-1$ and $1\leq p<\frac{4}{3}$ \cite{deBS}. They are {\it never} minimizers of the corresponding energy with fixed $L^2$ norm (actually this minimization problem has no solutions in this case).

It turns out that qualitative properties of solitary waves can be established for a large class of KP type equations. Ground state solutions are shown in \cite{deBS2} to be {\it cylindrically symmetric}, that is radial with respect to the transverse variable up to a translation of the origin.

On the other hand, {\it any} finite energy solitary wave is infinitely smooth (at least when the exponent $p$ is an integer) and decay with an {\it algebraic} rate $r^{-2}$ in two dimension and $r^{-3+\epsilon}$ in dimension three \cite{deBS2}. Actually the decay rate is sharp in the sense that a solitary wave {\it cannot} decay faster that $r^{-d},$ $d=2, 3.$

Moreover a precise {\it asymptotic expansion} of the solitary waves has been obtained by Gravejat (\cite{Gr}).

The fact whether the ground states are or are not minimizers of the Hamiltonian is strongly linked to the {\it orbital stability} of the set $\mathcal G_c$ of ground states of velocity $c$. Note that the {\it uniqueness}, up to the obvious symmetries, of the ground state of velocity $c$ is a difficult open problem, even for the classical KP I equation.

Saying that $\mathcal G_c$ is orbitally stable in $Y$ means that for $\phi_c\in \mathcal G_c,$ then for all $\epsilon >0, $ there exists $\delta>0$ such that if $u_0\in Y$\footnote {Some extra regularity on $u_0$ is actually needed.} is such that $||u_0-\phi_c||_Y\leq \delta,$ then the solution $u(t)$ of the Cauchy problem initiating from $u_0$ satisfies

$$\sup_{t\geq 0} \inf_{\psi \in \mathcal G}||u(t)-\psi||_Y\leq \epsilon.$$

Actually it is proved in \cite{deBS3} that the ground state solitary waves are orbitally stable in dimension two if and only if $1\leq p<\frac{4}{3}$ while there are {\it always unstable} in dimension three.

We will see below the link between stability and  blow-up in finite time of the solutions to the Cauchy problem.

\begin{remark}\label{stab}
No rigorous stability (or instability) result seems to be known for the Zaitsev solitary wave solution of the KP I equation.
\end{remark}

\subsection{The linear group}

As previously mentioned,  the linear part associated to {\it any} KP type equations \eqref{1.9} generates a group $S(t)= \exp(t(\partial_x L+\partial ^{-1}_x\partial^2_{yy}))$ which is unitary in all Sobolev spaces, in particular in $L^2$. For instance in the case of KP I/II equations, we have

$$S_{KP I/II}(t)=\exp( t(\partial_{xxx}\mp \partial_x^{-1}\partial^2_{yy}))$$

When the problem is set in $\R^2$, $S_{KP I/II}(t)$ has nice {\it dispersive} properties. In fact an analysis of the fundamental solution leads to the dispersive estimate \cite{Sa}

\begin{equation}\label{disp}
|| S_{KP I/II}(t)\phi||_{L^{\infty}(\R^2)}\leq \frac{C}{|t|}||\phi||_{L^1(\R^2)}, \quad \forall t\neq 0.
\end{equation}

Such an estimate leads to $L^p_t L^q_{xy}$ {\it Strichartz} estimates \cite{Sa} which are useful to solve the nonlinear Cauchy problem.  Typically, one gets estimates of the type (here $S(t)$ stands for the KP I{ \it or} the KP II group) :

\begin{equation}\label{Strich}
||S(t)\phi||_{L^q(\R;L^r(\R^2)}\leq C||\phi||_{L^2(\R)},
\end{equation}

for any $r\in \lbrack 2, +\infty),$ and where $q=q(r)$ is defined by $\frac{1}{q}=\frac{1}{2}-\frac{1}{r}.$

Similar estimates hold for the ``Duhamel integral"

$$\Lambda f(t)=\int_0^tS(t-s)f(s)ds.$$

Of course an estimate like \eqref{disp} is false when the spatial domain is the torus $\T^2$ or  the ``mixed " case $\R\times \T$ or $\T\times \R.$

On the other hand, dispersive estimates of the type \eqref{disp} hold true when the $u_{xxx}$ dispersive term in the KP equation is replaced by $\partial _x P(\partial_x)u,$ where $P$ is an even polynomial (see \cite {BAS}), in particular for the ``cubic-quintic" KP equations. In this case  the corresponding Strichartz estimates involve a {\it global} smoothing effect  \cite{BAS}.

Another interesting {\it dispersive} property of the KP group is the {\it local smoothing} property. 

It has been known since the seminal paper of Kato \cite{Ka} that 
the Airy group, that is $S(t)= e^{it \partial_x^3},$ satisfies the following local smoothing estimate

$$||\partial_xS(t)\phi)||_{L^2(0,T;L^2_{loc}(\R)}\leq C||\phi||_{L^2(\R))}.$$

In other words there is (locally in space) a gain of one derivative.

A similar property was established in \cite{Sa} for the KP II group. Actually, for any initial data $\phi\in L^2(\R^2),\;\partial_x^{-1}\phi \in L^{2}(\R^2)$, the corresponding solution $u$ of the {\it  linear} KP II equation satisfies the estimate, for any $T>0$ and $R>0$,

$$\int _{-T}^{T} \int _{-R}^{R}\int^{\infty}_{-\infty}\lbrack u_x^2(x,y,t)+|\partial _x^{-1}u_y(x,y,t)|^2\rbrack dxdydt\leq C(R,T)||\phi||^2_{L^2(\R^2)},$$

which also displays a local gain of one derivative.

A stronger smoothing effect holds for the {\it nonlinear} KP II equation when the  initial data decay fast to the right \cite{Le}. More precisely, if the initial data satisfies

$$\int_{\R^2} \lbrack \phi^2+\phi^2_{xxx}+(\partial_x^{-1}\phi_{yy})^2+x_{+}^L\phi^2+x_{+}^L\phi^2_x\rbrack dx dy<+\infty$$

for all positive integers $L$,  then the corresponding solution $u$ (which exists globally, see the next Section) satisfies
$u(.,t)\in C^{\infty}(\R^2)$  for all positive $t$'s.

A result of this type (with different hypothesis on the initial data) was established in \cite{LSV} for the nonlinear KP I equation.

\subsection{The nonlinear Cauchy problem}

All the KP type equations can be viewed as a linear skew-adjoint perturbation of the Burgers equation. Using this structure, it is not difficult (for instance by a compactness method)  to prove that the Cauchy problem is locally well-posed for data in the Sobolev spaces $H^s,$ $s>2$ in both dimensions (see \cite{U}, \cite{Sa}, \cite{IN} for results in this direction).

Unfortunately, this kind of result does not suffice to obtain the {\it global} well-posedness of the Cauchy problem. This would need to use the {\it conservation laws} of the equations. For {\it general} KP type equations, there are only two of them, the conservation of the $L^2$ norm and the conservation of the energy (Hamiltonian). For the general equation \eqref{genKP} where $f(u) =\frac{1}{p+1}u^{p+1},$ and without the transport term $u_x$ (which can be eliminated by a change of variable), the Hamiltonian reads

\begin{equation}\label{Ham}
E(u)=\frac{1}{2}\int\left[-u\,Lu+(\partial_{x}^{-1}u_y)^2+\frac{u^{p+2}}{p+2}\right],
\end{equation} 

and for the classical KP I/II equations

\begin{equation}\label{HamKP}
E(u)=\frac{1}{2}\int\left[u_x^2\pm(\partial_{x}^{-1}u_y)^2-\frac{u^{3}}{3}\right],
\end{equation} 

where the $+$ sign corresponds to KP I and the $-$ sign to KP II.

Note that the ``integrable" KP I and KP II equations possess more conservation laws, but not {\it infinitely many} as it is often claimed (see below).

In any case, it is clear that for KP II type equation, the Hamiltonian is useless to control any Sobolev norm, and to obtain the {\it global}  well-posedness of the Cauchy problem one should consider $L^2$ solutions, a very difficult task. On the other hand, for KP I type equations, one may hope (for a not too strong nonlinearity) to have a global control in the {\it energy} space, that is

\begin{equation}\label{energysp}
E= \lbrace u\in L^2,\; uLu\in L^2,\;\partial_x^{-1}u_y\in L^2 \rbrace.
\end{equation}

For the usual KP I equation, $E=Y$, which was defined in the previous subsection.

The problem is thus reduced to proving the local well-posedness of the Cauchy problem in spaces of very low regularity,  a difficult matter  which has attracted a lot of efforts in the recent years.

\begin{remark}
By a standard compactness method, one obtains easily the existence of global weak finite energy solutions (without uniqueness) to the KP I equation (see {\it eg} \cite{Tom}).
\end{remark}

At this point, there is a striking difference between the KP I and the KP II equation, linked to their different dispersion relations.
The  KP II equation can be solved by a Picard scheme implemented on the integral (Duhamel) formulation, implying that the flow map $u_0 \mapsto u(.,t) $ is {\it smooth} in natural Sobolev spaces. In this sense the KP II equation can be viewed as a {\it semi-linear} equation.

On the other hand, it has been proved in \cite{ST} for the periodic case and in \cite{MST} for the Cauchy problem in $\R^2$ that the KP I equation {\it cannot} be solved by such a procedure, in {\it any} natural Sobolev class. This implies that the flow map cannot be smooth, and actually it cannot be even uniformly continuous on bounded sets of the natural energy space $Y$ (see \cite{KTz}). This is a typical property of quasi-linear hyperbolic systems ({\it eg} the Burgers equation). It is thus natural to consider the KP I equation as a {\it quasi-linear} equation.

Those {\it ill/well posedness by iteration} properties  have a huge importance when solving the Cauchy problem since the methods are quite different in both cases. We will first focus on the classical KP I and KP II equations. 

For the KP II equation, a breakthrough was made by Bourgain \cite{Bou} who proved that the Cauchy problem for the KP II equation is locally (thus globally in virtue of the conservation of the $L^2$ norm) for data in $L^2(\R^2),$ and even in $L^2(\T^2).$ This result was later improved (\cite{Ta}, \cite{Tz}, \cite{Tz1} \cite{TaTz}, \cite{IM}, \cite{Ha}, \cite{HHK}) to allow the case of initial data in negative order Sobolev spaces.\begin{remark}
Bourgain's proof uses in a crucial way the fact  (both in the periodic and full space case) that the dispersion relation of the KP II equation induces the triviality of a certain {\it resonant set}.  A similar property was used by Zakharov \cite{Za2} to construct a Birkhoff formal form for the {\it periodic} KP II equation with small initial data. On the other hand, the fact that for the KP I equation the corresponding resonant set is non trivial is crucial in the construction of the counter-examples of \cite{ST}, \cite{MST} and is apparently an obstruction to the Zakharov construction for the periodic KP I equation.
\end{remark}

The KP II equation in three space dimensions has been studied in \cite{Tz}, \cite {ILM1}. The Cauchy problem is locally well-posed for initial data in anisotropic Sobolev spaces having $1^+$ derivative in $x$ and $0^+$ derivatives in $y,z$ in $L^2(\R^3)$. It is an open problem whether or not those solutions are globally defined.

We refer to \cite{Sa}, \cite{IN}, \cite{Ha}, \cite{KZ}, \cite{KZ2}, \cite{Gru}, \cite{HNS} for similar results on {\it generalized} KP II type equations. Note that for nonlinearities higher than quadratic, that is $f(u)= u^{p+1},\;p>1$ in \eqref{genKP}, and  $Lu=-u_{xx},$ one gets only the {\it local well-posedness}. Neither global existence nor blow-up in finite time are established.

On the other hand, global results are known to be true for higher order KP II equations, {\it eg} for the fifth order equation corresponding to  $Lu=-u_{xxxx}$ \cite{ST2}.

Most of those results are based on  iterative methods on the Duhamel formulation, using the {\it Fourier restriction} $X^{s,b}$ spaces of Bourgain, or variant of them, possibly with the injection of various linear dispersive estimates (see \cite{Sa}, \cite{BAS}).

We also would like to mention the paper by Kenig and Martel \cite{KM} where the Miura transform is used to prove the global well-posedness of a modified KP II equation.

In view of \cite{MST} the situation is quite different for the KP I equation. Note however that the KP I equation with a {\it fifth order} dispersion in $x$  is globally well-posed in the corresponding energy space \cite{ST3} (see also \cite{GPW} for local well-posedness results without the ``zero mass constraint" on the initial data).

For the classical KP I equation, the first global well-posedness result for arbitrary large initial data in a suitable Sobolev type space was obtained by Molinet, Saut and Tzvetkov \cite{MST2}. The solution is uniformly bounded in time and space.

The proof is based on  a rather sophisticated compactness method and 
uses the first invariants of the KP I equation to get global in time 
bounds. It is worth noticing that, while the recursion formula in 
\cite{ZS} gives formally a infinite number of  invariants, except for the first ones, those invariants do not make sense for functions belonging to $L^2(\R^2)$ based Sobolev spaces.

 For instance, the invariant which should control
$||u_{xxx}(.,t)||_{L^2}$ contains the $L^2$ norm  of  $\partial_x^{-1}\partial_y(u^2)$ which does not make sense  for a non zero function $u$ in the Sobolev space $H^3(\R^2)$. Actually (see \cite{MST1}, \cite{MST2}), one checks easily that if  $\partial_x^{-1}\partial_y(u^2)\in L^2(\R^2),$ then $\int_{\R^2} \partial _y(u^2) dx = \partial_y \int_{\R^2}u^2 dx \equiv 0,\; \forall y\in \R,$  which, with $u\in L^2(\R^2),$ implies that $u\equiv 0.$ Similar obstructions occur for the higher order ``invariants".

One is thus led to introduce a {\it quasi-invariant} (by skipping the non defined terms) which eventually will provide the desired bound. There are also serious technical difficulties to {\it justify} rigorously the conservation of the meaningful invariants along the flow and to control the remainder terms

The result of \cite{MST2} was extended by Kenig \cite{K} (who considered initial data in a larger space), and by Ionescu, Kenig and Tataru \cite{IKT} who proved that the KP I equation is globally well-posed in the energy space $Y$.

\begin{remark}
For both the KP I and KP II equations, the Inverse Scattering method provides only global well-posedness for {\it small} initial data. We refer to \cite{W} for the KP II equation and to \cite{Zh}, \cite{Su} for the KP I equation.
\end{remark}

Concerning the generalized KP I equations

\begin{equation}\label{genpKPI}
u_t+u_{xxx}- \partial_x^{-1}u_{yy}+u^pu_x=0,
\end{equation}

an important remark is that the anisotropic Sobolev embedding (see \cite{BIN})

$$\int_{\R^2} |u|^{p+2} dx dy\leq C ||u||_{L^2}^{\frac{4-p}{2}}||u_x||_{L^2}^p||\partial_x^{-1}u_y||_{L^2}^{\frac{p}{2}},$$

which is valid for $0\leq p\leq 4,$ implies that the energy norm $||u||_Y$ is controlled in term of the $L^2$ norm and of the Hamiltonian 

$$\int_{\R^2} 
\left(\frac{1}{2}u_x^2+\frac{1}{2}(\partial_x^{-1}u_y)^2-\frac{|u|^{p+2}}{(p+1)(p+2)}\right) dx dy,$$

if and only if $p<\frac{4}{3}.$

This suggest global well-posedness when $p<\frac{4}{3}$ and blow-up when $p\geq \frac{4}{3}$, a fact which will be proved in the next Section.

\subsection{Blow-up issues}

We consider here a class of generalized KP equations

\begin{equation}\label{genpKP}
u_t+u_{xxx}+\epsilon \partial_x^{-1}u_{yy}+u^pu_x=0,
\end{equation}

 where $\epsilon =1$ for the KP II equations and $-1$ for the KP I equations.

The KP I equations are {\it focusing} while the KP II are {\it 
defocusing} and to this aspect they share some properties with the focusing (resp.~defocusing) nonlinear Schr\"{o}dinger equation (NLS). In particular one could expect the blow-up in finite time of the solutions of generalized KP I equations with a sufficiently high nonlinearity.

For the focusing NLS, one can prove the blow-up in finite time of the 
$L^2$ norm of the spatial gradient by using a virial identity (see 
{\it eg} \cite{SuSu}). For the generalized KP I equation, one can use a
 ``transverse" virial identity to prove the blow up in finite time of the $L^2$ norm of the {\it transverse} gradient (that is to say in $2D$ of $||u_y(.,t)||_{L^2(\R^2)}).$ This identity has been formally derived in \cite {TuFa} and rigorously justified in \cite{Sa} (one needs to prove a local well-posedness theory in a weighted $L^2$ space). It reads for  solutions of the generalized KP equations in a suitable functional class:
 
 \begin{equation}\label{virial}
 \frac{d^2}{dt^2} I(t)=4\epsilon\lbrack -pE(t)+\frac{p}{2}\int_{\R^2}u_x^2dxdy+\frac{(4-p)}{2}\epsilon \int_{\R^2}(\partial_x^{-1}u_y)^2dxdy\rbrack,
 \end{equation}
 
% where $\epsilon =1$ for the KP II equations and $-1$ for the KP I equations,
\begin{equation}
    E(t)=\int_{\R^2} \lbrack \frac{1}{2}(u_x^2-\epsilon (\partial_x^{-1}u_y)^2 -\frac{u^{p+2}}{(p+1)(p+2)}\rbrack dx dy,
    \label{genKPener}
\end{equation}

 and 
 
 $$I(t)=\int_{\R^2}y^2u^2 dxdy.$$
 
 A first direct consequence of \eqref{virial} for KP I equations ($\epsilon =-1$) is that, when the initial data $u_0$ (is ``regular" and) belongs to the space 
 
 $$\Sigma =\lbrace \phi \in E, \; y\phi\in L^2(\R^2),\; \phi \neq 0,\; E(\phi)\leq 0 \rbrace,$$

 the corresponding solution of the Cauchy problem cannot remain smooth for all times if $p\geq4.$ More precisely, there should exist $T^{*}>0$ such that
 
 \begin{equation}\label{BU1}
 \lim_{t\rightarrow T^*}||u_y(.,t)||_{L^2}=+\infty.
 \end{equation}
  
  A similar result is obtained in the three-dimensional case {\it ie} when $\partial_x^{-1}u_{yy}$ in \eqref{genpKP} is replaced by 
  $\partial_x^{-1}\Delta^{\perp}u,$ where $\Delta^{\perp}$ is the 
  transverse Laplacian in the $y$ and $z$ variables. 
  
  The critical value is then $p=2$ and \eqref{BU1} should be replaced by 
  
   \begin{equation}\label{BU2}
 \lim_{t\rightarrow T^*}||\nabla^{\perp}u(.,t)||_{L^2}=+\infty.
 \end{equation}
 
 \begin{remark}
 It is readily seen that the set $\Sigma$ is not empty. It contains in particular functions of the type $\pm \lambda\partial_x^2 \exp (-x^2-y^2)$ when $ \lambda $ is a large positive number, the sign depending on the ``parity" of $p$.
 \end{remark}
 
 Those results are not fully satisfactory since we already observed that the natural energy norm is controlled by he $L^2$ norm and the Hamiltonian if and only if $0<p<\frac{4}{3}$ in the two-dimensional case. In the three dimensional case, the critical exponent is  $ p= \frac{4}{5}$ which suggests that the ``standard"  ($p=1)$ 3D- KP I equation is not supposed to have global solutions for arbitrary initial data.
 
 This result was improved by Liu \cite{L1} who used invariant sets of the generalized KP  I flow together with the virial argument above to prove the existence of initial data leading to blow-up in finite time of $||u_y(.,t)||_{L^2(\R^2)}$ when $p\geq \frac{4}{3}$. He also proved a strong instability result of the solitary waves when $2<p<4.$
 
 Similar results are established in \cite{L2} for the three-dimensional generalized KP I equation, in the range $1\leq p<\frac{4}{3}.$
 
 \begin{remark}
 A detailed analysis of the blowing-up solution (profile, blow-up rate,..) is lacking.
 \end{remark}
 
 \begin{remark}
 No rigorous finite time blow-up has been established for the {\it periodic} generalized KP I equation.
 \end{remark}
 
 \begin{remark}
No result seems to be known concerning the long time behavior (finite time blow-up versus global existence)
of the solutions to the generalized KP II equations \eqref{genpKP} when $\epsilon =1$ and $p\geq 2$).
 \end{remark}
 
 \begin{remark}
 The finite time blow- up does not {\it a priori} invalidate the 
 applicability of KP equations as asymptotic models.  As all asymptotic models, they are supposed to describe the dynamics of a more general system, for a suitable scaling regime  (here weakly nonlinear, long wave with a wavelength anisotropy) and on relevant (finite!) time scales.
  For instance, in the water wave context, the KP equations cease to be relevant models for times larger than $O(\frac{1}{\epsilon^2})$.
 \end{remark}

\subsection{The Cauchy problem in the background of a non localized solution.}

As was previously noticed, it is quite natural in view of  transverse stability issues  to consider the Cauchy problem for KP equations on the background of a solitary wave of the underlying KdV equation.

The Inverse Scattering method has been used formally in \cite{AV} and 
rigorously in \cite {AV2} to study the Cauchy problem for the KP II 
equation with non decaying data along a line, that is $u(0,x,y)= 
u_{\infty}(x-vy)+ \phi(x,y)$ with $\phi(x,y)\to 0$ as $x^2+y^2\to \infty$ and $u_{\infty}(x)\to 0 $ as $|x|\to \infty$.  Typically, $u_{\infty}$ is the profile of a traveling wave solution $U({\bf k}.{\bf x}-\omega t)$ with its peak localized on the moving line ${\bf k}.{\bf x}=\omega t.$ It is a particular case of the $N$- soliton of the KP II equation discovered by Satsuma \cite{Sat} (see the derivation and the explicit form when $N=1,2$ in the Appendix of \cite{NMPZ}). As for all results obtained for KP equations by using the Inverse Scattering method, the initial perturbation of the non-decaying solution is supposed to be small enough in a weighted $L^1$ space (see \cite{AV2} Theorem 13). A similar result has been established by Fokas and Pogrebkov \cite{FP} for the KP I equation. 

On the other hand, PDE techniques allow to consider {\it arbitrary large perturbations} of a  class of non-decaying  solutions of the KP I/II equations.

\vspace{0.5cm}
We will therefore  study  the initial value problem for the KP I and KP II  equations
\begin{equation}
(u_t+u_{xxx} +u u_x )_x \pm u_{yy} =0,
\label{KP2}
\end{equation}
where $u=u(t,x,y)$, $(x,y)\in\R^2$, $t\in\R$, with initial data
\begin{equation}\label{1.2a}
u(0,x,y)=\phi(x,y)+\psi_c(x,y),
\end{equation}
where $\psi_c$ is the profile\footnote{This means that
$\psi(x-ct,y)$ solves (\ref{KP2}).} of a non-localized (i.e., not
decaying in all spatial directions) traveling wave of the KP I/II
equations moving with speed $c\neq 0$.

We recall that, contrary to the KP I equation, the KP II equation does not possess any {\it localized in both directions} traveling wave solution.

The background solution $\psi_c$ could be for instance the line soliton (1-soliton) of the Korteweg-
de Vries (KdV) equation \eqref{kdvsol},
% \begin{equation}\label{KdV}
% s_c(x,y)=\frac{3c}{2}\,{\rm cosh}^{-2}\Big(\frac{\sqrt{c}\, x}{2}\Big),
% \end{equation}
or the N-soliton solution of the KdV equation, $N\geq 2,$ or in the case of the KP I equation, the Zaitsev soliton.

The KdV N-soliton is of course considered as a two dimensional (constant in $y$) object. 

This problem can be attacked in two different ways. Either one considers a localized (in $x$ and $y$) perturbation, or the perturbation is localized in $x$ and periodic in $y$.

Both frameworks have been applied to the KP I and KP II equations, and they lead to {\it global well-posedness} of the Cauchy problem  for arbitrary large initial perturbations. More precisely, for the KP I equation, the Cauchy problem  with initial data which are localized perturbations of the KdV N-soliton or of the Zaitsev soliton is globallly well-posed \cite{MST3}. The corresponding result for the KP II equation is proved in \cite{MST4}. 

The global well-posedness of the KP II equation with arbitrary initial data in $L^2(\R\times\T)$ is proved in \cite{MST4}. The corresponding result for the KP I equation with initial data in a suitable Sobolev space  is established in \cite{IK}.

\subsection{Transverse stability issues}

The KP I and KP II equations behave quite differently with respect to the {\it transverse} stability of the KdV $1$-soliton.

Zakharov \cite{Za} has proven, by exhibiting an explicit perturbation 
using the integrability, that the  KdV $1$-soliton is {\it nonlinearly} unstable for the KP I flow. Rousset and Tzvetkov \cite{RTz} have given an alternative proof of this result, which does not use the integrability, and which can  thus be implemented on other problems ({\it eg} for nonlinear Schr\"{o}dinger  equations).

The {\it nature} of this instability is not known (rigorously) and it will be investigated in our numerical simulations.

On the other hand, Mizomachi and Tzvetkov \cite{MiTz} have recently proved the $L^2(\R\times \T)$ orbital stability of the KdV $1$-soliton for the KP II flow. The perturbation is thus localized in $x$ and periodic in $y$.  The proof involves  in particular the Miura transform used in \cite{KM} to established the global well-posedness for a modified KP II equation. 

Such a result is not known (but expected) for a perturbation which is localized in $x$ and $y$.

\subsection{Comparison between the KP and KP/BBM type equations}

As previously noticed, the KP/BBM equations \eqref{1.5} are relevant counterparts to KP type equations, for instance to the KP II equation \footnote{Note however that when the {\it BBM trick} is applied to the KdV or to the KP equation with strong surface tension (Bond number greater than $\frac{1}{3}$) one gets a {\it ill-posed} equation, so that what is called {\it KP-I/BBM} equation is merely a mathematical object without, as far as we know, modelling relevance...}  

The Cauchy problems for those {\it regularized} KP equations has been studied in \cite{ST4}, for initial data in the natural energy space. We refer also to \cite{BLT} for more regular initial data and for a stability analysis of the solitary waves in the {\it focusing} case.

The comparison of the dynamics of the KdV and BBM equation has been 
investigated both theoretically and numerically in \cite{BPS}. A 
similar study was made in \cite{Ma} between the KP and KP/BBM equations. Recall that KdV and BBM (resp. KP and KP/BBM) equations are supposed to describe the {\it same} dynamics, on relevant time scales.

\subsection{Long time dynamics}

KP type equations can be classified into two categories, the {\it defocusing} ones (KP II type), and the {\it focusing } ones (KP I type). For the former, the long time dynamics is expected to be governed by dispersion and scattering. In particular, at least for small initial data the $\sup$ norm should decay with the linear rate, that is $1/t$ (see \cite{Sa}).

On the other hand, the dynamics of the latter case is expected to be 
governed by the solitary waves and blow-up in finite time is expected 
in the {\it supercritical} case. To this extent the situation is reminiscent of that of the nonlinear Schr\"{o}dinger equations \cite{SuSu}.

Relatively few mathematical results are know concerning those issues, 
and one goal of the present paper is to present numerical simulations 
in support of (or in motivation of) various conjectures on the long time dynamics.

For the generalized KP equation

$$u_t+u_{xxx}\pm\partial_x^{-1}u_{yy}+u^pu_x=0,$$

it was proven in \cite{HNS} that if $p\geq 2$ and if the initial data is sufficiently small, then 

$$||u(t)||_{\infty}\leq C(1+|t|)^{-1}(\log (2+|t|)^{\kappa},\quad ||u_x(t)||_{\infty}\leq C(1+|t|)^{-1},$$

for all $t\in \R$, where $\kappa=1$, if $p=2$ and $\kappa =0$ if $ p\geq 3.$

Moreover, a large time asymptotic for $u_x$ can be obtained.

\begin{remark}
No such result seems to be known for the KP  II equation, $p=1$ (see 
however Theorem 9.3 in \cite{Za2}  for a scattering result with small initial data by Inverse Scattering methods). We also mention \cite{HHK} where  a scattering result for small data in the scale invariant non isotropic homogeneous Sobolev space $\dot{H}^{-1/2,0}(\R^2)$ is obtained. 

The large time behavior of solutions to the KP II equation with arbitrary large data is unknown. An interesting open problem is that of possible uniform  in time bounds on the higher Sobolev norms. It has been proven in \cite{IMTz} that the eventual growth in time of Sobolev norms is at most polynomial (see also  \cite{Tz2} for long time bounds for the periodic KP II equation).
\end{remark}

The long time dynamics of KP  I type equations for arbitrarily large data is also an important open problem. For the standard KP I equation, it is expected that the ground state solutions (whether or not the lump is one of them) should play a role in this dynamics as the KdV solitary waves for the corresponding equation.

\subsection{Initial data in the Schwartz class}
\label{schwarzini}

One can prove that the Schwartz space $\mathcal S(\R^d)$ of $C^{\infty}$ functions with rapid decay is not preserved by the (linear and nonlinear) flows of KP type equations, in fact the {\it rapid decay} is not preserved as it is easily seen on the  solution of the linear KP I or KP II equations which reads in Fourier variables 

$$\hat{u}(\xi_1,\xi_2,t)=\hat{u_0}(\xi_1,\xi_2)\exp (it(\xi_1^3\pm \frac{\xi_2^2}{\xi_1}))$$

Taking for instance   $u_0$ equal to a gaussian function, one sees that $u(.,t)$ is not even in $L^1(\R^2)$ for $t\neq 0$ since this would imply by Riemann-Lebesgue theorem that  $\hat{u}(.,t)$ is a continuous function of $\xi$ which obviously is not the case.  This proves in particular that the solution at time $t$ cannot decay faster than $\frac{1}{r^2}$ at infinity, where $r=\sqrt{x^2+y^2}.$
The same conclusion can be drawn in the nonlinear case. In fact the solution of the nonlinear equation with the same gaussian initial data has the Duhamel representation in Fourier variables

$$\hat{u}(\xi_1,\xi_2,,t)=\hat{u_0}(\xi_1,\xi_2)\exp (it(\xi_1^3\pm \frac{\xi_2^2}{\xi_1})) +\frac{i}{2}\int_0^t  \xi_1\exp (i(t-s)(\xi_1^3\pm \frac{\xi_2^2}{\xi_1}))\widehat{u^2}(\xi_1, \xi_2,s) ds$$

and one checks that the integral part defines, for any fixed $t>0$, a continuous function of $\xi$ while the free part is not continuous as shown previously.

A stronger decay of solutions is obtained when extra (rather unphysical) conditions are imposed on the initial data $u_0$, namely the vanishing of $\widehat{u_0}(\xi_1,\xi_2)$ at higher and higher order as $\xi_1 \to 0.$

The proofs above extend to {\it any} KP type equation like \eqref{genKP}, the crux of the matter being the singularity at $\xi_1=0$ of the kernel 
$\frac{\xi_2^2}{\xi_1}.$

The fact that the solutions of KP type equations do not decay fast at infinity yields   difficulties discussed in Section 3.1 for the realistic numerical simulations of the solutions in the whole space.

\subsection {Varia} We mention briefly here a few other mathematical questions concerning KP type equations. 

Many papers have been recently devoted to the study of initial boundary-value problems for various dispersive equations such as the KdV equation, motivated in particular by applications to control theory.

A IBV problem for the KP II equation on a strip or a half space has been studied in \cite{MP}. We do not know of any similar result for the KP I equation. 

Many works have been recently devoted to the control or stabilization of linear or nonlinear dispersive equation. 

We are not aware of any work on the control of (linear or nonlinear) KP type equations. A unique continuation property (which could be useful in control theory) for the KP  II equation is established in \cite{Pant}. A corresponding result for the KP I equation does not seem to be known.

As was previously alluded to  in the context of water waves, one sometimes need for physical reasons to add a {\it viscous} term to the KP equations leading to interesting mathematical problems. For a viscous term of the Burgers type, that is $-u_{xx}$, Molinet and Ribaud \cite {MoRi} have established the global well-posedness of the Cauchy problem for data in $L^2(\R^2)$ for the KP-II-Burgers equation and in the energy space $Y$ for the KP-I-Burgers equation. This result has been extended by Kojok \cite{Ko} in the case of the KP-II-Burgers equation for initial data in the anisotropic Sobolev space $H^{s,0}(\R^2)$, for any $s>-\frac{1}{2}.$

A sharp decay asymptotic analysis of the solutions of KP-Burgers equations as $t\to +\infty$ is carried out in  \cite{Mol}.

\section{Numerical study of generalized KP equations}
In this section we will numerically study the stability of exact 
solutions to the KP equations for various perturbations, and the 
blow-up phenomenon in generalized KP equations. 

\subsection{General setting}
We will now study numerically the solution to initial value 
problems for the generalized KP equations (\ref{genpKP}).
We will only consider functions $u$ that are periodic in $x$ and $y$ with 
period $L_{x}2\pi$ and $L_{y}2\pi$ respectively or in 
the Schwartzian space of rapidly decreasing functions, which are 
essentially periodic for numerical purposes if restricted to an 
interval of length $L$ on which the functions decrease below values that can be 
numerically distinguished from zero ($10^{-16}$ in double precision). 
As noted in subsection~\ref{constraint},  all solutions of 
(\ref{genpKP}) satisfy in this case for $t>0$ the 
constraint     $\int_{-\pi L_{x}}^{\pi L_{x}}u_{yy}dx=0$.
To obtain solutions which are smooth at time $t=0$ (for convenience 
we only consider times $t\geq 0$), the initial data also have to 
satisfy this constraint. It is possible to  numerically solve 
initial value problems which do not satisfy the constraint, see for 
instance \cite{klspma}, but the non-regularity in time is 
numerically difficult to resolve. Therefore we will only consider 
initial data that satisfy this condition, typically by imposing that 
$u_{0}(x,y)$ is the $x$-derivative of a periodic or Schwartzian function. This even holds for 
perturbations of exact solutions we will consider in the following.

The imposed periodicity on the solutions allows the use of Fourier 
spectral methods. Equation \eqref{genpKP} is equivalent to 
\begin{equation}
\partial_t \, \widehat{u} +\frac{i\xi_{1}}{p+1}\widehat{u^{p+1}}-i \xi_{1}^{3}\, \widehat{u}+ 
\lambda \frac{i\xi_{2}^{2}}{\xi_{1}+i0} \, \widehat u=0.
\label{KPf}
\end{equation}
The singular multiplier $-i /\xi_{1}$ is regularized in standard way 
as $-i/(\xi_{1}+i0)$. 
Numerically this is achieved by computing $-i/(\xi_{1}+i\epsilon)$ 
where $\epsilon$ is some small number of the order of the rounding 
error. Since we only consider functions $u$ satisfying 
the constraint $\xi_{2}^{2}\widehat{u}(t,\xi_{1}=0,\xi_{2})=0$, 
this does not lead to any numerical problems.

For a numerical solution equations (\ref{KPf}) are approximated by a 
truncated Fourier series, i.e., a discrete Fourier transform for 
which fast algorithms exist. With this approximation one obtains a 
system of ordinary differential equations for which an initial value 
problem has to be solved. This system is stiff due to the high order 
spatial derivatives, especially since we want to resolve strong 
gradients. Standard explicit schemes would require prohibitively 
small time steps for stability reasons in this case, whereas 
unconditionally stable implicit schemes are typically not efficient, 
especially if one aims at high precision. In \cite{etna} we compared 
several fourth order schemes for the KdV equation, and in \cite{KR} 
a similar analysis was performed for KP. It was found that 
\textit{exponential time differencing} schemes are most efficient in this 
case, see for instance the reviews \cite{minwri},\cite{ho}. There are several 
such schemes known in the literature which show similar performance. We 
use here the fourth order scheme proposed by Cox and Matthews 
\cite{coxma}. 

As noted in subsection~\ref{schwarzini}, solutions to generalized KP 
equations for  Schwartzian initial $u_{0}$ in general will not
stay in the Schwartz class, unless $u_{0}$ satisfies in addition an infinite number of constraints. 
For generic initial data, the solutions $u(t)$ will develop `tails' 
decaying only algebraically in certain spatial directions. This 
implies for the used periodic setting that Gibbs phenomena will appear on the borders of the 
computational domain for larger times. These jump discontinuities affect the numerical 
accuracy. Thus a much higher resolution than for KdV 
(where solutions to Schwartzian initial data stay in the Schwartz 
space) would be needed to obtain high numerical precision. In the 
below computations this effect is addressed by choosing a larger 
computational domain than would be necessary for similar KdV 
computations. The algebraic tails will lead due to the periodicity 
condition to `echoes' which are, however, in most cases of comparatively small 
amplitude.

The numerical accuracy of the computations is controlled by 
propagating exact solutions and by computing conserved quantities of 
KP which are numerically functions of time due to unavoidable 
numerical errors, see \cite{KR}. As discussed before generalized KP equations 
conserve the $L^{2}$-norm (or mass) and the energy (\ref{genKPener}).
Since the latter quantity contains the term $\partial_x^{-1} 
\partial_y u(t,x,y)$ which implies in Fourier space a division by 
$\xi_{1}$, an operation that affects accuracy in the determination of 
the quantity, we will trace mass conservation. In 
\cite{KR} it was shown that relative mass conservation 
$\Delta:=|1-M(t)/M(0)|$ ($M$ being the numerically $L^{2}$-norm)
overestimates typically the accuracy of the numerical solution in an 
$L^{\infty}$ sense by 2 orders of magnitude. Thus we always aim at a 
$\Delta$ smaller than $10^{-4}$ to ensure accuracy of plotting 
precision. Obviously the quantity $\Delta$ only makes sense if the 
spatial numerical resolution is of the same order, i.e., if the 
Fourier coefficients $\widehat{u}(t,\xi_{1},\xi_{2})$ decrease for large 
$\xi_{1},\xi_{2}$ to the same order of magnitude. The precision of the computations 
will therefore always be controlled by mass conservation and 
Fourier coefficients. In addition we check whether the found 
numerical results remain the same within the considered accuracy if 
computed with higher or lower resolution. 

Exact solutions to the KdV equation are automatically solutions to 
both KP equations without explicit $y$-dependence. Thus KdV solitons are 
\emph{line solitons} of KP, i.e., solutions that are exponentially 
localized in one spatial direction and infinitely extended in the 
other. As mentioned above the KdV soliton (\ref{kdvsol}) is  unstable for KP I and 
lineary stable for KP II. 
To test the code we first propagate initial data from the KdV 
soliton $u_{0}=u_{sol}(x-x_{0})$ for $c=4$ and 
$x_{0}=-2L_{x}$ ($L_{x}=8$ and $L_{y}=8$, 
$x\in[-\pi,\pi]L_{x}$, $y\in[-\pi,\pi]L_{y}$)
with $N_{x}=2^{8}$ and $N_{y}=2^{10}$ Fourier modes with $N_{t}=3200$ 
time steps from $t=0$ until $t=6$. For both KP I and KP II we find at $t=6$ that the relative 
mass conservation $\Delta=3.89*10^{-8}$ and that the 
$L^{\infty}$-norm of the difference between numerical and 
exact solution is equal to $9.96*10^{-7}$. This shows that the code 
is well able to propagate the exact solution even for the unstable 
case. It also illustrates that the quantity $\Delta$ can be used as 
expected as an indicator for the numerical accuracy. 
In the following subsections we will test how periodic
perturbations of order unity of the KdV soliton are propagated by the 
KP equations.

\subsection{Perturbations of the KdV soliton for KP II}
The similarity between KP I and KP II for the exact KdV soliton initial 
data disappears if we perturb these data for instance by 
\begin{equation}
    u_{p}=6(x-x_{1})\exp(-(x-x_{1})^{2})\left(\exp(-(y+L_{y}\pi/2)^{2})
    +\exp(-(y-L_{y}\pi/2)^{2})\right).
    \label{pert}
\end{equation}
The perturbations are 
Schwartzian in both variables and satisfy the constraint 
(\ref{constraint}). They are of the same order of magnitude as the KdV 
soliton, i.e., of order $0(1)$, and thus
test the nonlinear stability of the KdV soliton. 

Notice that the computations are carried out on a doubly periodic 
setting, i.e., $\mathbb{T}^{2}$ and not on $\mathbb{R}^{2}$. Thus the 
perturbations cannot be radiated to infinity, but will stay in the 
computational domain. Since they are of a similar size as the initial 
soliton, they will lead to a perturbed soliton though the KdV soliton is 
known to be stable for KP II. One possibility is that the initial 
perturbation is essentially homogeneously distributed over the whole 
computational domain. Another possibility would be that they form 
eventually a KdV soliton which is, however, not observed in the below 
examples.

\begin{figure}
[!htbp]
\begin{center}
\includegraphics[width=\textwidth]{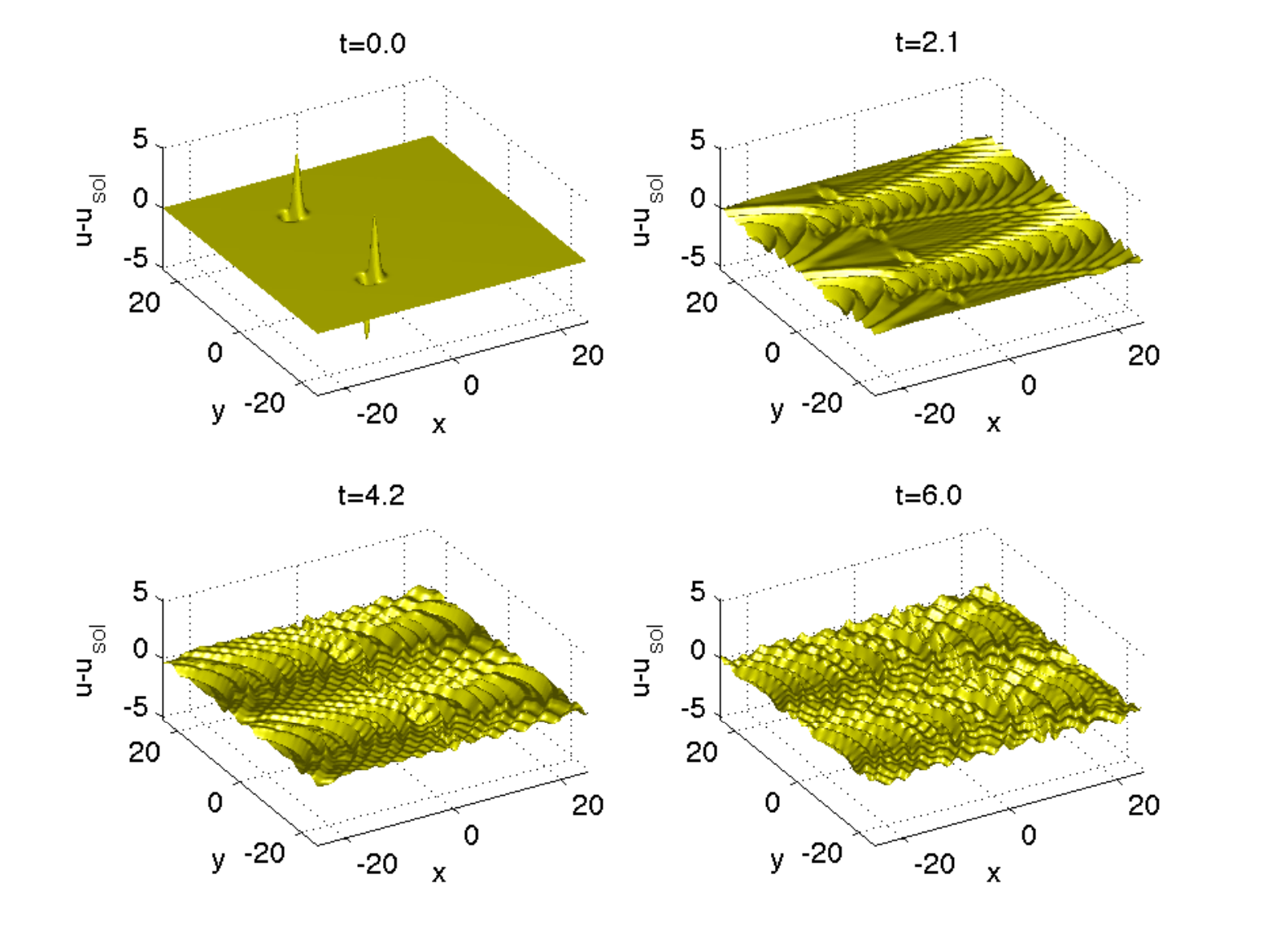} 
\caption{Difference of the solution to the  KP II equation  for initial 
data given by the KdV soliton plus perturbation, 
$u(0,x,y)=u_{sol}(x+2L_{x},0)$ and perturbation $u_{p}=6(x-x_{1})\exp(-(x-x_{1})^{2})\left(\exp(-(y+L_{y}\pi/2)^{2})
    +\exp(-(y-L_{y}\pi/2)^{2})\right)$, $x_{1}=-L_{x}$ and the 
KdV soliton for various values of $t$.}
\label{kpIIsolperd}
\end{center}
\end{figure}

For KP II we use the same parameters for the computation as for the 
tests for the KdV soliton, but now for initial data 
$u_{0}(x,y)=12\mbox{sech}^{2}(x+2L_{x})-u_{p}$ ($x_{1}=x_{0}/2$), 
i.e., a 
superposition of the KdV soliton and the not aligned perturbation. At time $t=6$ 
we obtain for the numerical conservation of the mass 
$\Delta=4.29*10^{-8}$,  almost the same precision as for the 
propagation of the soliton. The result can be seen in Fig.~\ref{kpIIsolperd}. The 
perturbation is dispersed in the form of tails to infinity which 
reenter the computational domain because of the imposed periodicity. 
The soliton appears to be unaffected by the perturbation which 
eventually seems to be smeared out in the background of the soliton.

\begin{figure}
[!htbp]
\begin{center}
\includegraphics[width=\textwidth]{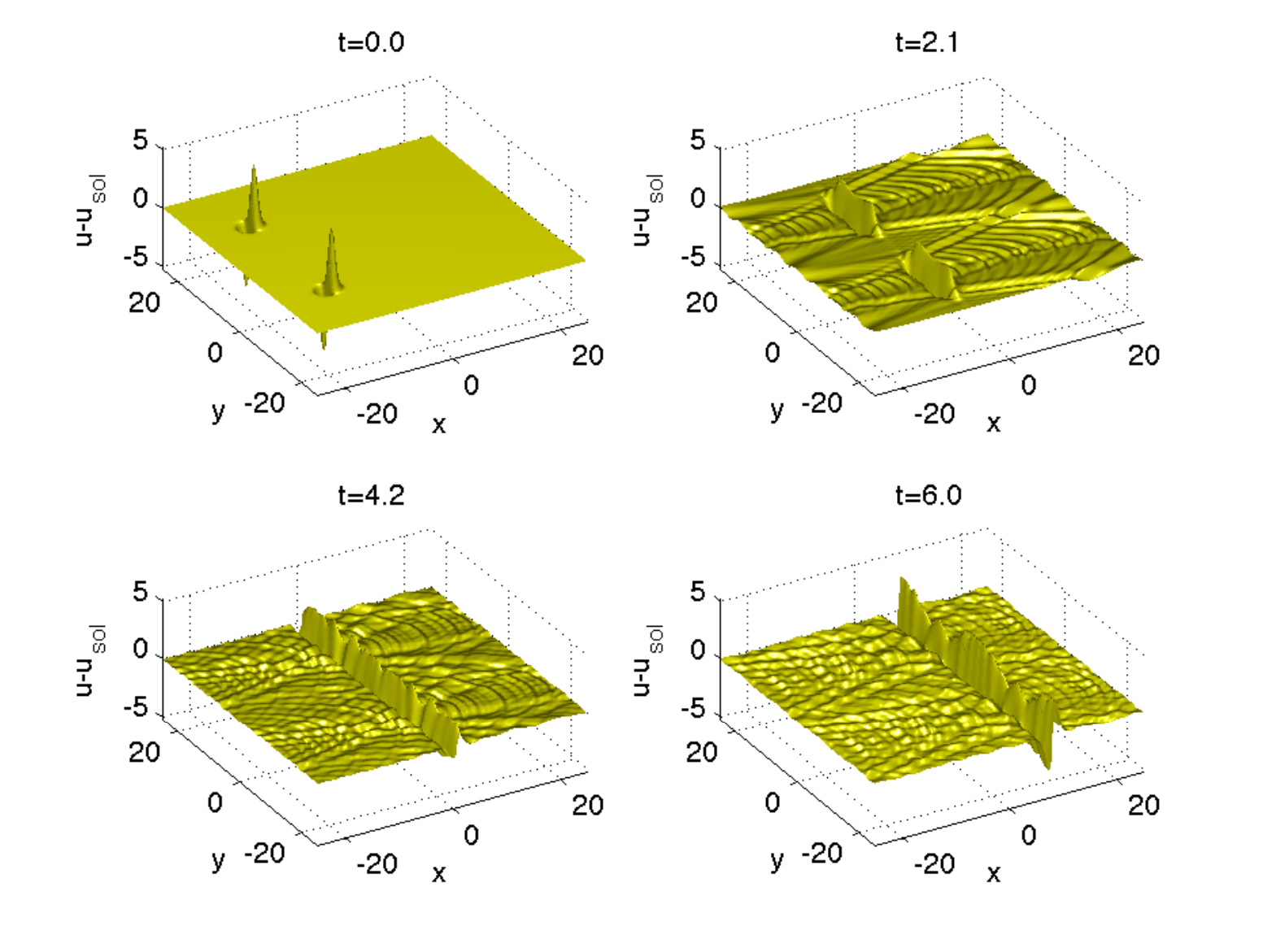} 
\caption{Difference of the solution to the  KP II equation  for initial 
data given by the KdV soliton $u_{sol}(x+2L_{x},0)$ 
plus perturbation $u_{p}=6(x-x_{1})\exp(-(x-x_{1})^{2})\left(\exp(-(y+L_{y}\pi/2)^{2})
    +\exp(-(y-L_{y}\pi/2)^{2})\right)$, $x_{1}=-2L_{x}$, and the 
KdV soliton for various values of $t$.}
\label{kpIIsolper}
\end{center}
\end{figure}

\begin{figure}
[!htbp]
\begin{center}
\includegraphics[width=0.8\textwidth]{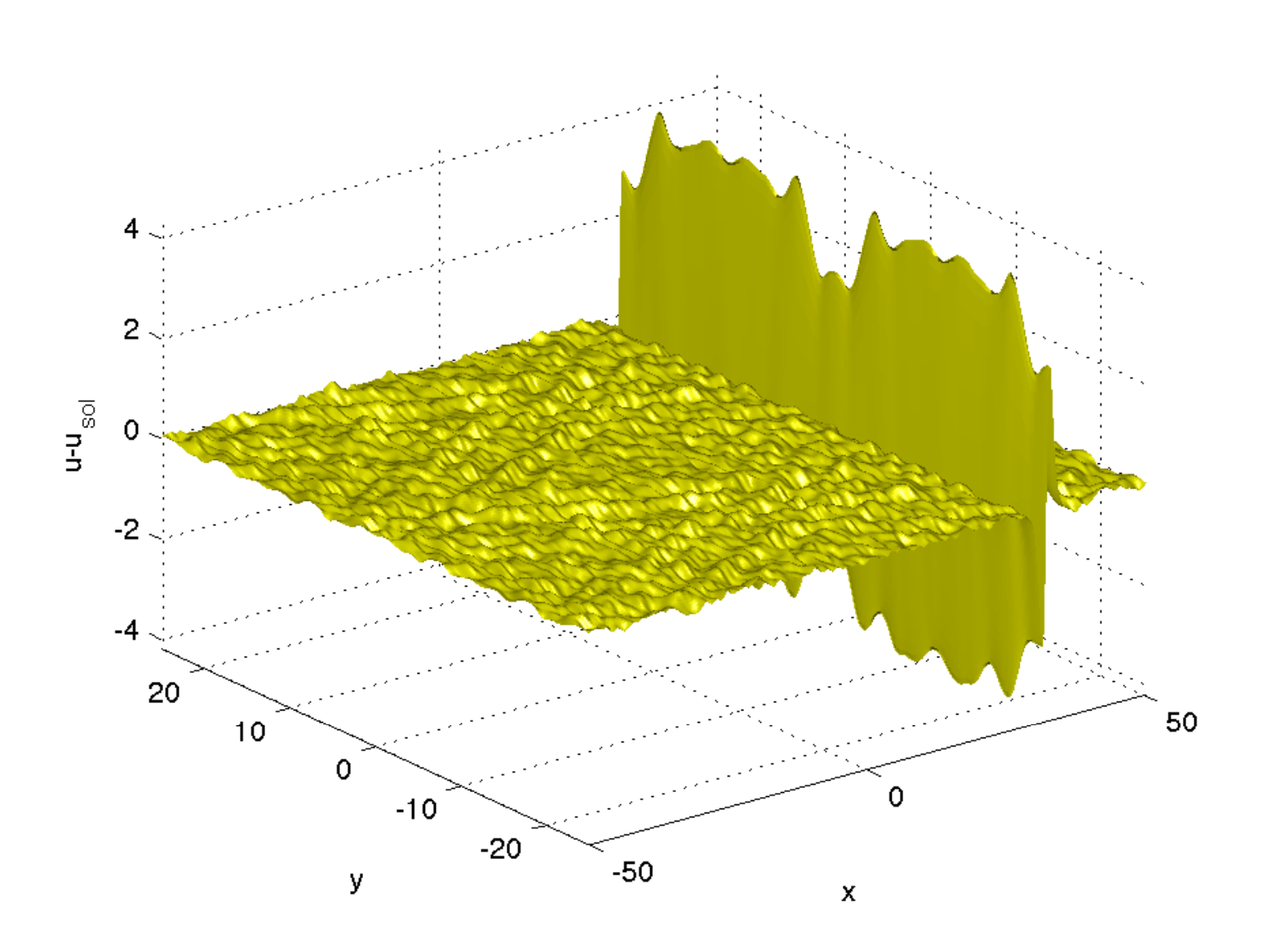} 
\caption{Difference of the solution to the  KP II equation  for the 
situation of Fig.~\ref{kpIIsolper} and the initial
KdV soliton for $t=16$.}
\label{kpIIsolperfit}
\end{center}
\end{figure}

The situation changes slightly if the perturbation and the initial 
soliton are centered around the same $x$-value initially, i.e., the 
same situation as above with $x_{1}=x_{0}=-2L_{x}$. 
In Fig.~\ref{kpIIsolper} we show the 
difference between the numerical solution and the KdV soliton 
$u_{sol}$ for several times for this case.
It can be seen that the initially localized perturbations spread in 
$y$-direction, i.e., orthogonally to the direction of propagation and 
take finally themselves the shape of a line soliton. The 
perturbations appear to modulate the soliton here. It is difficult 
numerically to determine the long time behavior. Extending the 
computation till time $t=16$, one finds a similar behavior as for 
$t=6$ as can be seen in  Fig.~\ref{kpIIsolperfit}.
Thus it appears that the perturbations of the KdV soliton travel here 
partly with the soliton and lead to oscillations around it. 
Since the speed of the soliton is not affected, it is still the 
original, not a soliton with slightly increased $c$ which would change 
both height and speed of the traveling wave. But these perturbations 
traveling with the soliton seem to be present for very long times.

\begin{figure}
[!htbp]
\begin{center}
\includegraphics[width=\textwidth]{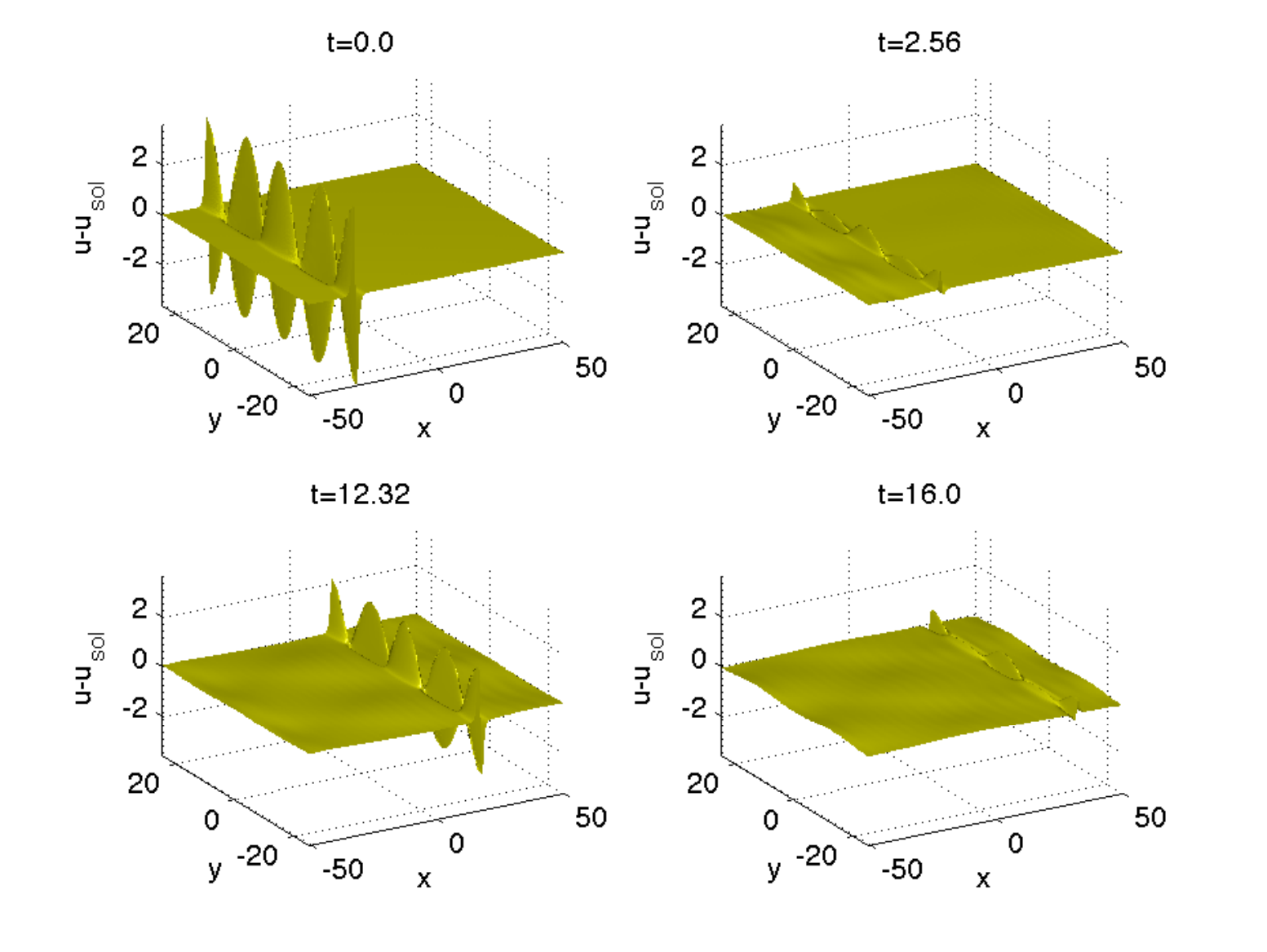} 
\caption{Difference of the solution to the  KP II equation  for the initial 
data $u_{0}(x) = 12\mbox{sech}^{2}(x+0.4\cos(2y/L_{y}))$ and 
the 
KdV soliton for several values of $t$.}
\label{kpIIsolper2}
\end{center}
\end{figure}

\begin{figure}
[!htbp]
\begin{center}
\includegraphics[width=0.8\textwidth]{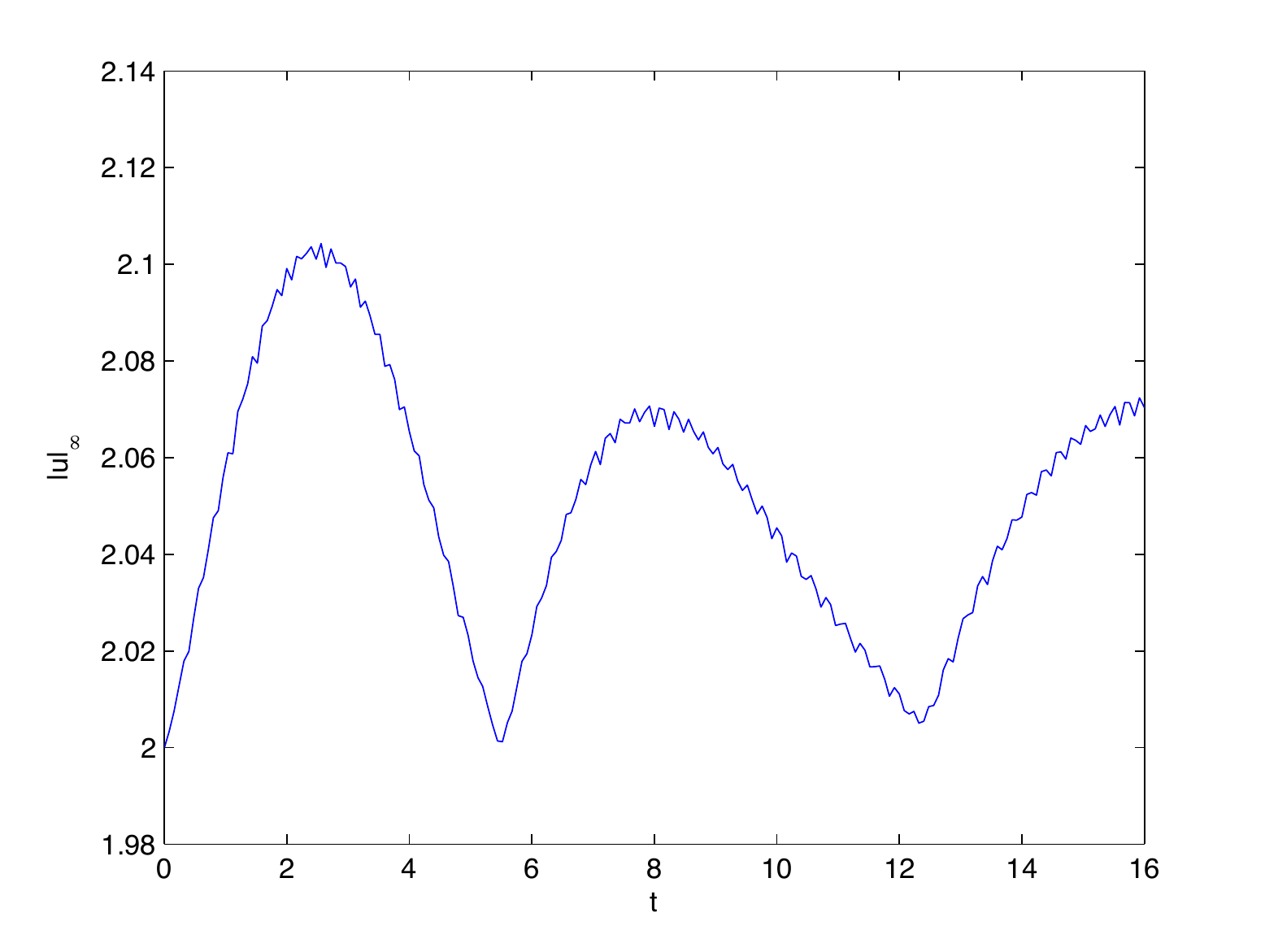} 
\caption{$L^{\infty}$-norm of the solution to the  KP II equation  for the 
situation of Fig.~\ref{kpIIsolper2} in dependence of $t$.}
\label{kpIIsolpercinf}
\end{center}
\end{figure}

A similar behavior is observed if the soliton shape itself is 
deformed, for instance in the form of the initial data 
$u_{0}(x) = 12\mbox{sech}^{2}(x+0.4\cos(2y/L_{y}))$. In 
Fig.~\ref{kpIIsolper2} the difference of the solution for KP II for this 
initial data and the line soliton can be seen. The computation is 
carried out with $N_{x}=2^{10}$, $N_{y}=2^{7}$, $L_{x}=16$, $L_{y}=8$ 
and $N_{t}=6400$ with a final $\Delta\sim 10^{-7}$.  It can be seen 
that the solution approaches the KdV soliton with perturbations 
oscillating around it.
This is even more visible in Fig.~\ref{kpIIsolpercinf} where the 
$L^{\infty}$-norm of the solution is given. It can be seen that the 
latter oscillates around the value for the KdV soliton. Numerically 
it is difficult to decide whether these oscillations will finally 
disappear.

\subsection{Perturbations of the KdV soliton for KP I}

Perturbations of the KdV soliton are propagated in a completely 
different way by the KP I equation. It is known that the line soliton 
is unstable for this equation. One of the reasons for this is 
the existence of lump solitons (\ref{Lump}).
Note that the lump soliton only decays algebraically in both spatial 
directions, as $|x|,|y|\to \infty$. Because of this and the ensuing 
Gibbs phenomena in the periodic setting used here for the  numerical 
simulations, the appearance of lumps 
will always lead to a considerable drop off in accuracy. In contrast to the line 
solitons it is localized in both spatial directions. Explicit 
solutions for multi-lumps are known, see \cite{manakov}. It was
shown by Ablowitz and Fokas \cite{AF}  that solutions to the KP I 
equation with small norm will asymptotically decay into lumps and radiation.

\begin{figure}
[!htbp]
\begin{center}
\includegraphics[width=\textwidth]{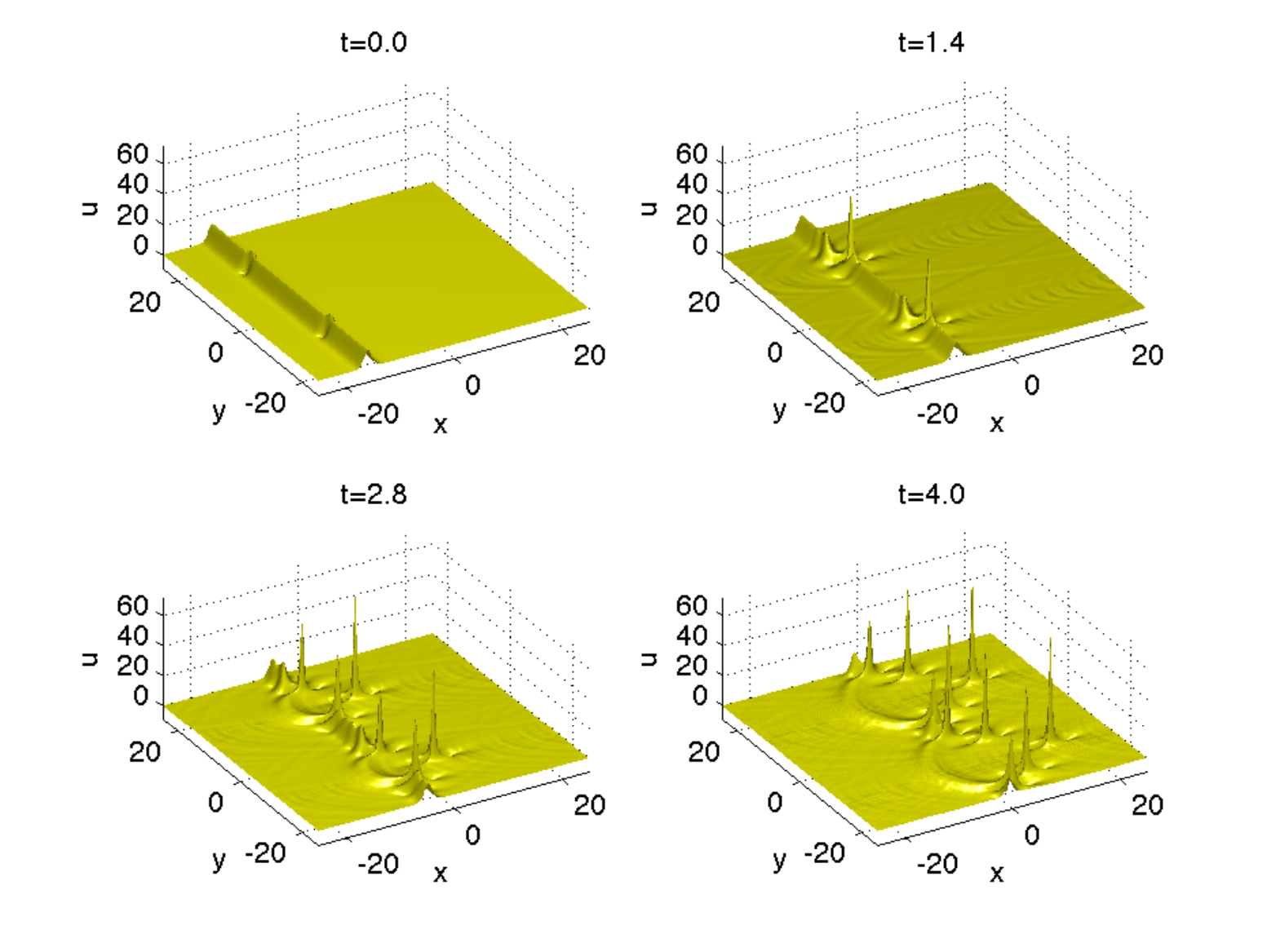} 
\caption{Solution to the  KP I equation  for initial 
data given by the KdV soliton $u_{sol}(x+2L_{x},0)$ 
plus perturbation $u_{p}=6(x-x_{1})\exp(-(x-x_{1})^{2})\left(\exp(-(y+L_{y}\pi/2)^{2})
    +\exp(-(y-L_{y}\pi/2)^{2})\right)$, $x_{1}=-2L_{x}$,  for various values of $t$.}
\label{kpIsolper}
\end{center}
\end{figure}

Exactly this behavior can be seen in Fig.~\ref{kpIsolper} for the perturbed initial data of a 
line soliton with the perturbation (\ref{pert}) and 
$x_{1}=x_{0}=-2L_{x}$, 
i.e., the same setting as studied in Fig.~\ref{kpIIsolper} for KP II. Here the 
initial perturbations develop into 2 lumps which are traveling with 
higher speed than the line soliton. The formation of these lumps 
essentially destroys the line soliton which leads to the formation of 
further lumps. It appears plausible that for sufficiently long times 
one would only be able to observe lumps and small perturbations which 
will be radiated to infinity if studied on $\mathbb{R}^{2}$. 
Notice that the algebraic fall off towards infinity causes the 
above mentioned numerical 
problems. This is reflected by the fact that the Fourier 
coefficients decrease much slower than in the case of KP II such that 
$N_{x}=2^{9}$ and $N_{y}=2^{10}$ Fourier modes are needed to ensure 
that the Fourier coefficients go down by at least 5 orders of 
magnitude in $\xi_{1}$ and $\xi_{2}$. For the relative mass conservation 
one obtains $\Delta=4.3*10^{-5}$ with $N_{t}=6400$.

\begin{figure}
[!htbp]
\begin{center}
\includegraphics[width=0.8\textwidth]{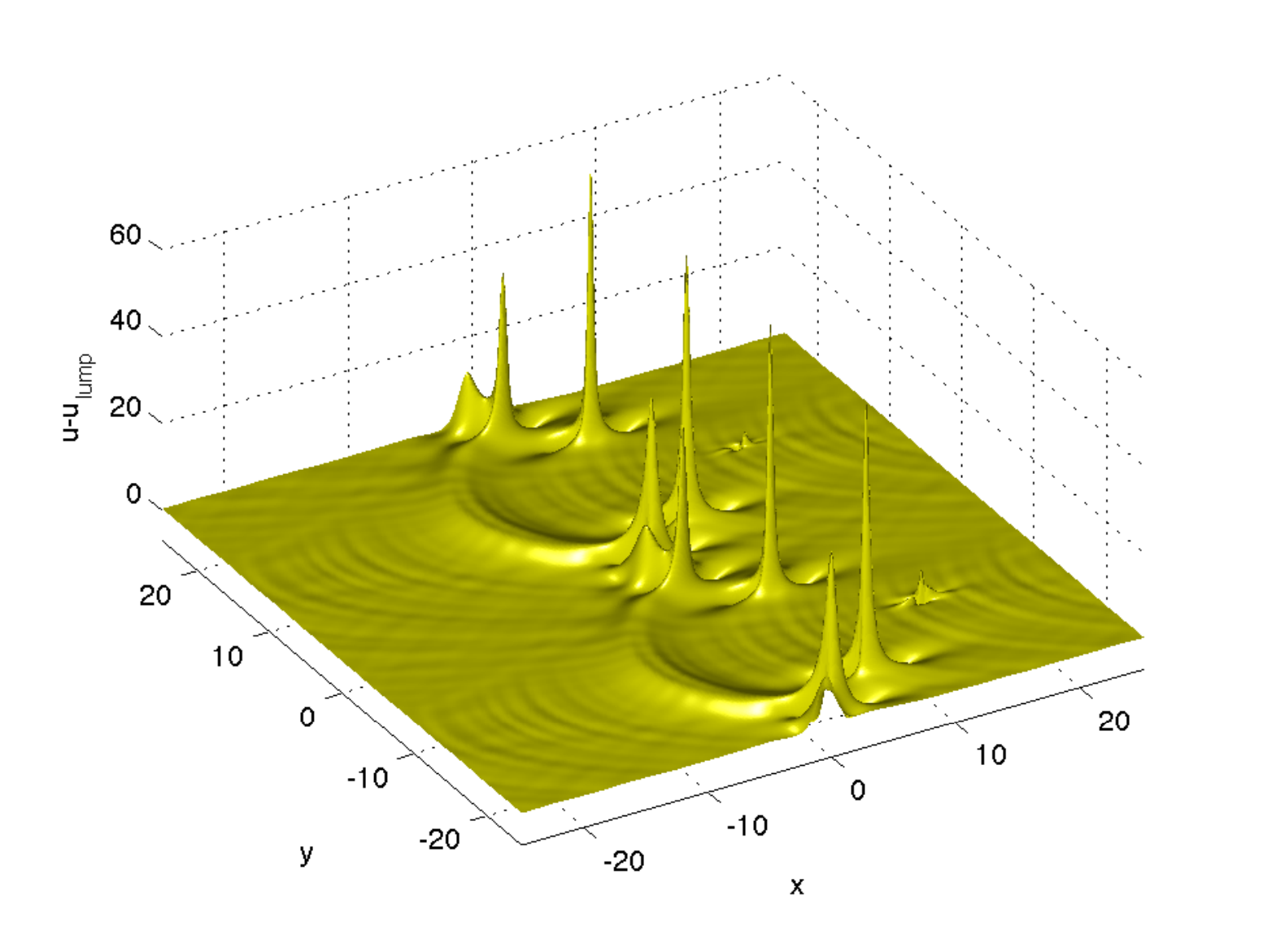} 
\caption{Difference of the solution to the  KP I equation  in 
Fig.~\ref{kpIsolper} for $t=4$ and two lump solutions fitted at the 
peaks farthest to the right. Only very small peaks remain of these 
`lumps' indicating that they develop asymptotically into true lumps.}
\label{kpIsolperl}
\end{center}
\end{figure}

We can give some numerical evidence for the validity of the 
interpretation of the peaks in Fig.~\ref{kpIsolper} as lumps in an 
asymptotic sense. We can identify numerically a certain peak, i.e., 
obtain the value and the location of its maximum. With these 
parameters one can study the difference between the KP solution and 
a lump with these parameters to see how well the lump fits the peak. 
This is illustrated for the two peaks, which formed first and which 
have therefore traveled the largest distance in Fig.~\ref{kpIsolperl}. 
Obviously one cannot expect a true lump in this case since there are 
still some remains of the line soliton to be seen. A multi-lump 
solution might be a better fit, but here we mainly want to illustrate 
the concept which obviously cannot fully apply at  the studied small 
times. Nonetheless Fig.~\ref{kpIsolperl} illustrates convincingly 
that the observed peaks will asymptotically develop into lumps. 

\begin{figure}
[!htbp]
\begin{center}
\includegraphics[width=\textwidth]{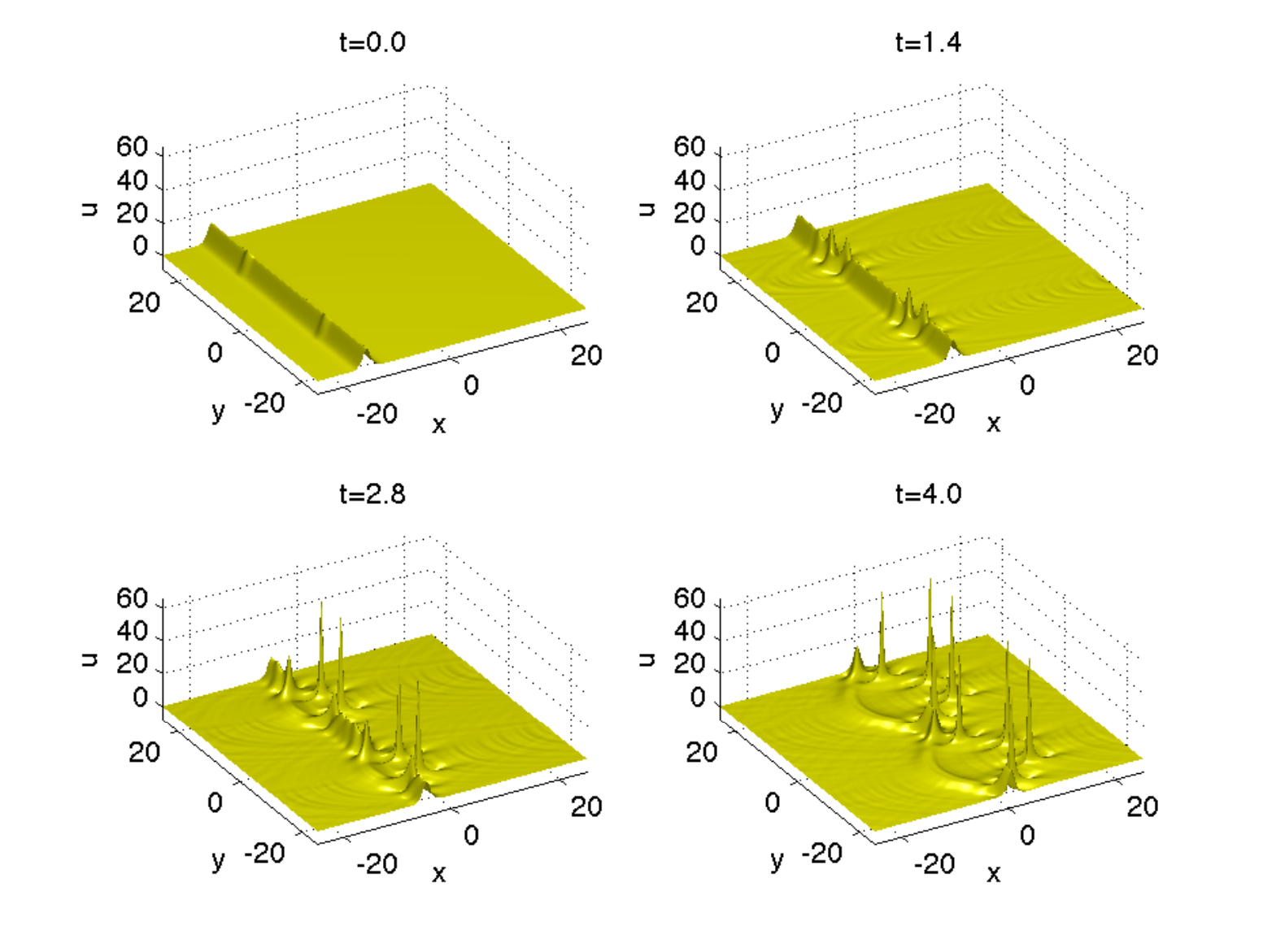} 
\caption{Solution to the  KP I equation  for initial 
data given by the KdV soliton $u_{sol}(x+2L_{x},0)$ 
plus perturbation $u_{p}=-6(x-x_{1})\exp(-(x-x_{1})^{2})\left(\exp(-(y+L_{y}\pi/2)^{2})
    +\exp(-(y-L_{y}\pi/2)^{2})\right)$, $x_{1}=-2L_{x}$,
for various values of $t$.}
\label{kpIsolperp}
\end{center}
\end{figure}

The precise pattern of lumps forming obviously depends on the 
perturbation. In Fig.~\ref{kpIsolperp} one can see the KP I solution 
for the same situation as in Fig.~\ref{kpIsolper} except for the 
different sign of the perturbation $u_{p}$ (\ref{pert}). 
In this case two pairs of lump form first to be followed at later 
times by a collection of lumps. 

\begin{figure}
[!htbp]
\begin{center}
\includegraphics[width=\textwidth]{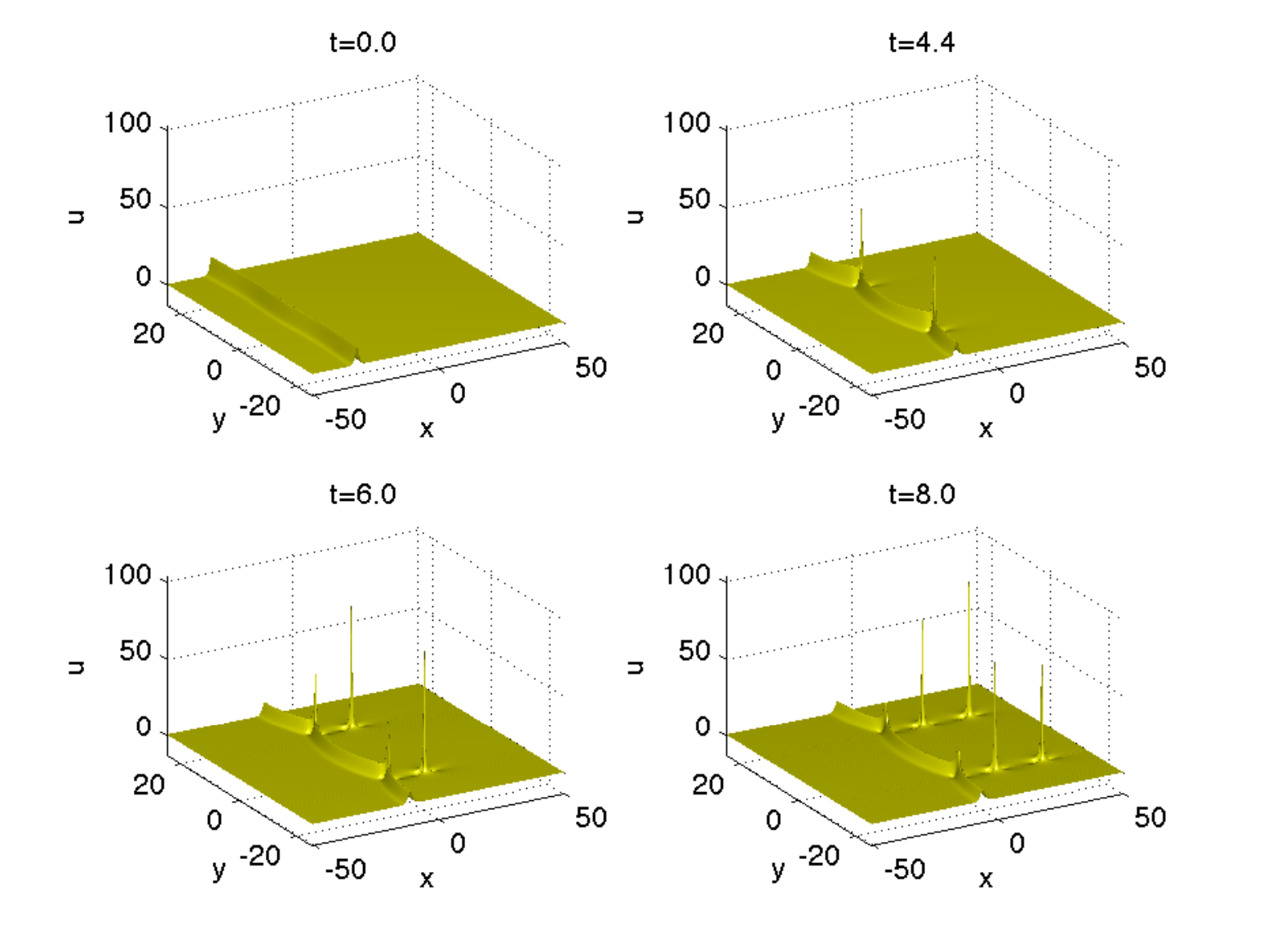} 
\caption{Solution to the  KP I equation  for initial 
data $u_{0}(x) = 12\mbox{sech}^{2}(x+0.4\cos(2y/L_{y}))$ 
for various values of $t$.}
\label{kpIsolperc}
\end{center}
\end{figure}

Perturbing the line soliton in the form $u_{0}(x) = 
12\mbox{sech}^{2}(x+0.4\cos(2y/L_{y}))$, one finds the solution 
shown in Fig.~\ref{kpIsolperc}. Two lumps form at the deformed line 
soliton farthest to the right. Further lumps then form at later times 
where the first detached from the deformed soliton. It appears that 
the perturbation develops in this case into a chain of lumps. The 
computation is carried out with $N_{x}=N_{y}=2^{11}$ modes and 
$10^{4}$ time steps for $L_{x}=16$ and $L_{y}=8$. Due to the 
propagation of lumps of large magnitude for an extended time, the relative mass 
conservation is comparatively low in this case, $\Delta\sim 
3*10^{-3}$.  

\begin{figure}
[!htbp]
\begin{center}
\includegraphics[width=\textwidth]{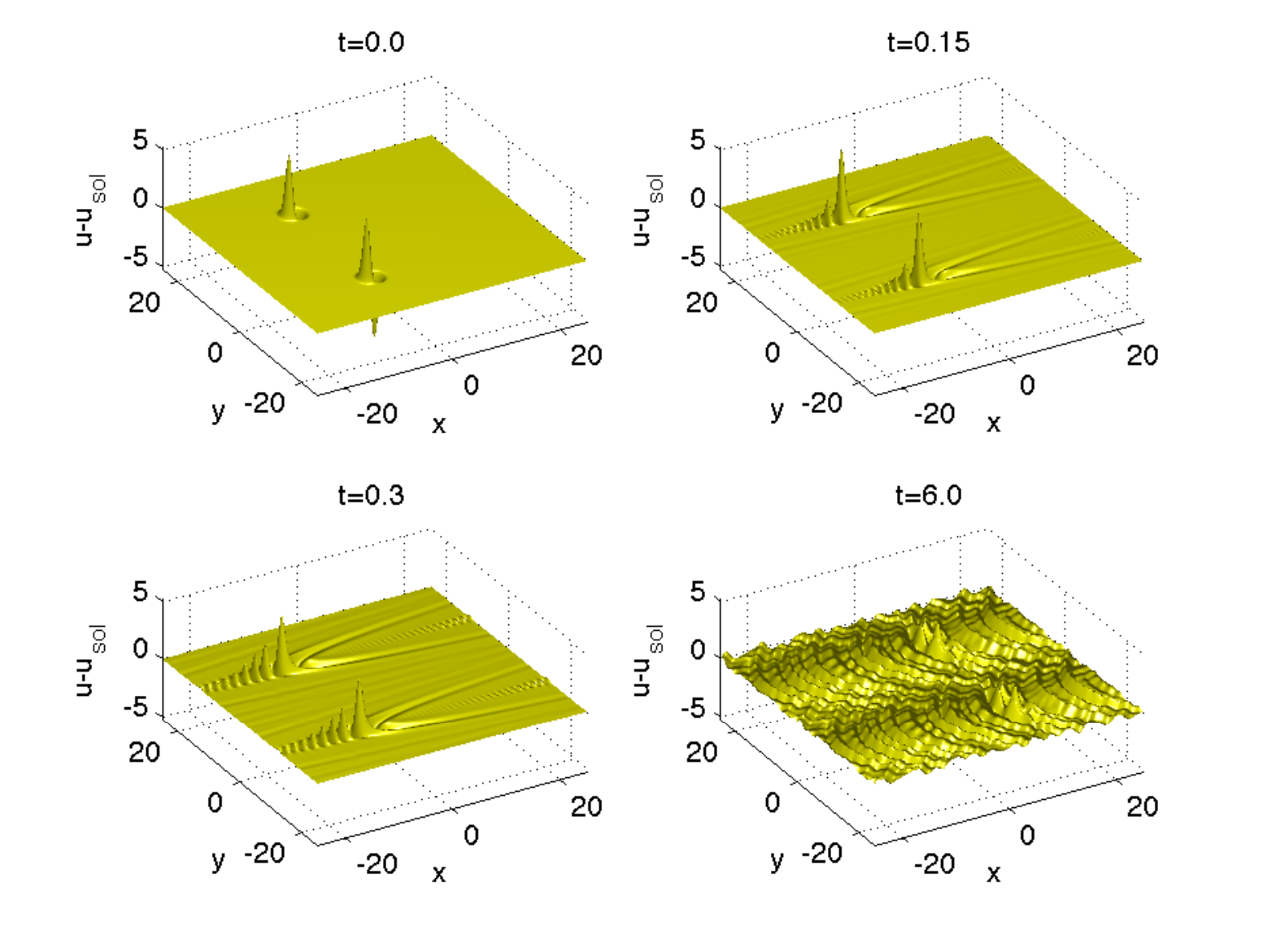} 
\caption{Difference of the solution to the  KP I equation  for initial 
data given by the KdV soliton $u_{sol}(x+2L_{x},0)$ 
plus perturbation $u_{p}=6(x-x_{1})\exp(-(x-x_{1})^{2})\left(\exp(-(y+L_{y}\pi/2)^{2})
    +\exp(-(y-L_{y}\pi/2)^{2})\right)$, $x_{1}=-L_{x}$ and the 
KdV soliton for various values of $t$.}
\label{kpIsolperd}
\end{center}
\end{figure}

The above numerical experiments indicate that the line soliton is 
unstable against general perturbations. Whereas this might be true on long 
time scales, this is not necessarily the case on intermediate time scales. 
To illustrate this `meta-stability' we consider in 
Fig.~\ref{kpIsolperd} a perturbation of the line soliton as before 
$u_{0}(x,y)=12\mbox{sech}^{2}(x+2L_{x})+u_{p}$, but this 
time with $x_{1}=x_{0}/2$, i.e., perturbation and soliton are well 
separated.
The figure shows the difference between KP I solution and line 
soliton. It can be seen that the soliton is essentially stable on the shown time 
scales. The perturbation leads to algebraic tails towards 
positive $x$-values and to dispersive oscillations as studied in 
\cite{klspma}. Due to the imposed periodicity both of these cannot 
escape the computational domain and appear on the respective other 
side. The important point is, however, that though the oscillations 
of comparatively large amplitude hit the line soliton quickly after 
the initial time, its shape is more or less unaffected till $t=6$.  

\begin{figure}
[!htbp]
\begin{center}
\includegraphics[width=\textwidth]{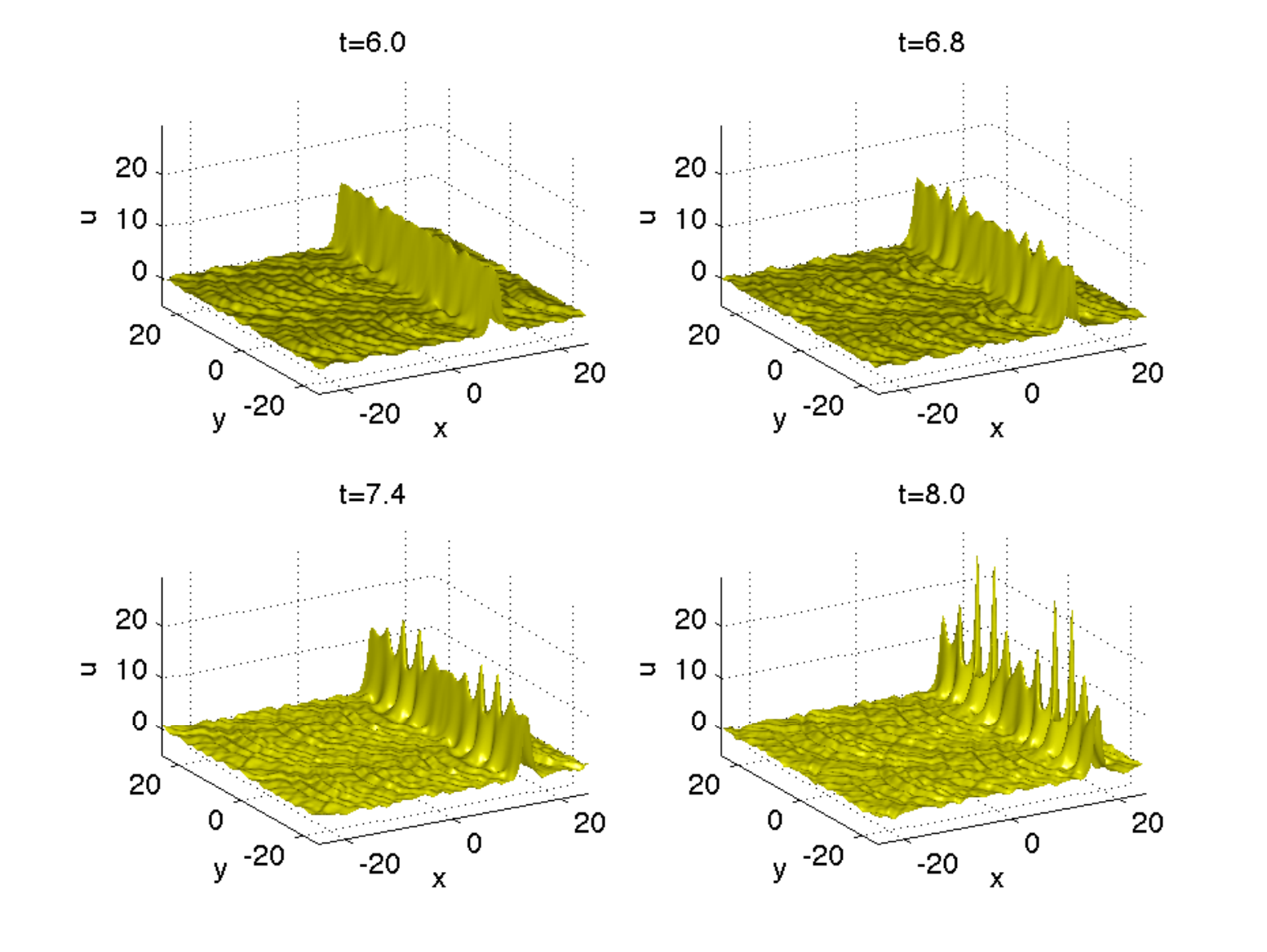} 
\caption{Solution to the  KP I equation  for initial 
data given by the KdV soliton $u_{sol}(x+2L_{x},0)$ 
plus perturbation $u_{p}=6(x-x_{1})\exp(-(x-x_{1})^{2})\left(\exp(-(y+L_{y}\pi/2)^{2})
    +\exp(-(y-L_{y}\pi/2)^{2})\right)$, $x_{1}=-L_{x}$ and the 
KdV soliton for various values of $t$.}
\label{kpIperp3}
\end{center}
\end{figure}

The picture changes, however, if  the code runs for a longer time as 
can be seen in Fig.~\ref{kpIperp3}. The soliton 
develops for times greater than 6 lumps after having stayed close to 
its original shape before despite strong perturbations. The solution 
for times greater than 6 remains the same within the given precision 
if the resolution is changed. 
For instance if we compare the result with $N_{t}=6400$ and 
$N_{t}=12800$, the maximal difference between the numerical solutions 
if of the order $10^{-4}$. The quantity $\Delta$ in the latter case 
is of the order of $10^{-8}$, whereas the Fourier coefficients 
decrease to $10^{-8}$.  Doubling the resolution in $x$ leads to a 
difference of the order of $10^{-9}$ between the solutions.  Thus one 
can conclude that the above solution is in fact correct within 
numerical precision.

\begin{figure}
[!htbp]
\begin{center}
\includegraphics[width=0.8\textwidth]{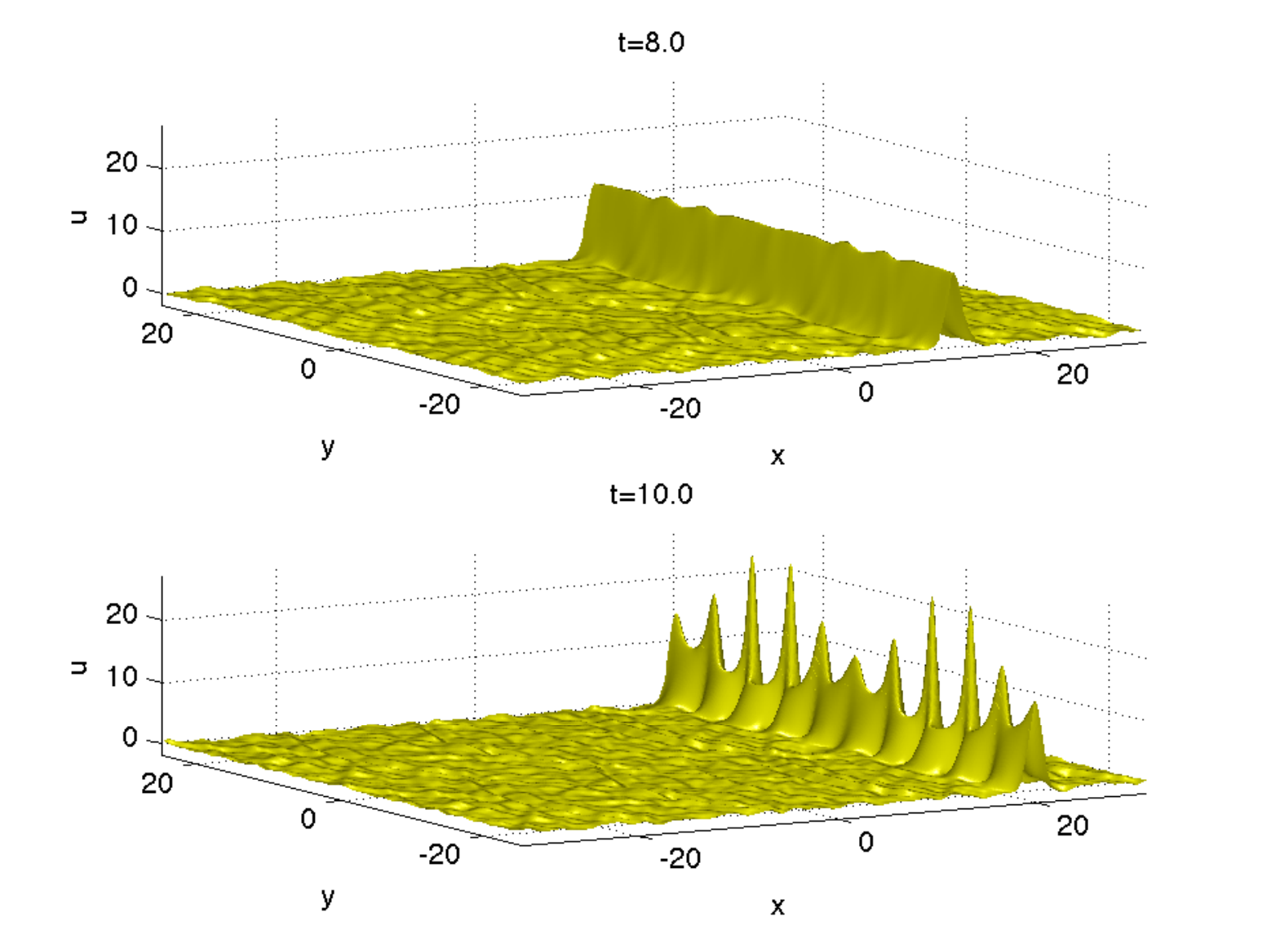} 
\caption{Solution to the  KP I equation  for initial 
data given by the KdV soliton $u_{sol}(x+2L_{x},0)$ 
plus perturbation $u_{p}=6(x-x_{1})\exp(-(x-x_{1})^{2})\left(\exp(-(y+L_{y}\pi/2)^{2})
    +\exp(-(y-L_{y}\pi/2)^{2})\right)$, $x_{1}=-L_{x}$ and the 
KdV soliton as in Fig.~\ref{kpIperp3}  ($L_{x}=8$) above and 
$L_{x}=10$ below.}
\label{kpIperpL10}
\end{center}
\end{figure}

Of course the periodic setting has an important effect here. Changing 
$L_{y}$ from 8 to 10 whilst keeping all other parameters unchanged 
does not affect the solution within the given precision. However the 
picture changes if we pass from $L_{x}=8$ to $L_{x}=10$ as can be 
seen in Fig.~\ref{kpIperpL10}. At time $t=8$ the soliton is still 
mainly unperturbed, and only then the lump formation starts. 
This implies that the instability of the KdV soliton against lumps is 
only triggered when the former is close to the computational 
boundary. Consequently it appears that some form of
`meta-stability' of the KdV soliton is numerically 
established. For large $L_{x}$ a perturbed KdV soliton can exist for 
long times, and the same could be true on 
$\mathbb{R}\times\mathbb{T}$. But we cannot 
make any predictions on the time scales where this would be the case.

\begin{figure}
[!htbp]
\begin{center}
\includegraphics[width=\textwidth]{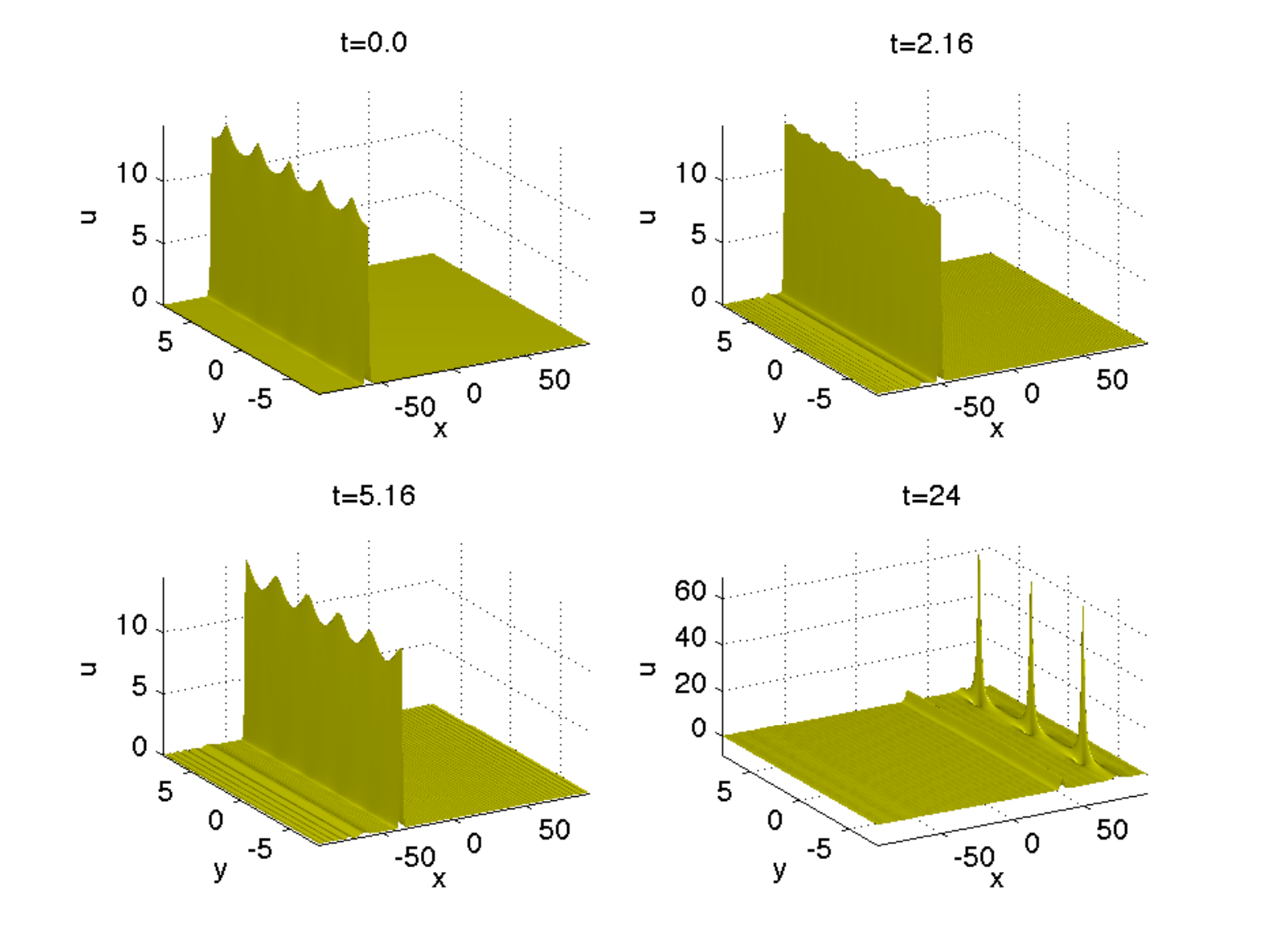} 
\caption{Solution to the  KP I equation  for initial 
data given by the KdV soliton 
plus perturbation given by the Zaitsev solution, 
$u_{0}=u_{sol}(x+L_{x},0)+0.1u_{Zait}(x+L_{x},y,0)$, 
$\alpha=1$, $\beta=0.5$, $L_{x} = 10$, $L_{y} = 5/\delta$ 
for various values of $t$.}
\label{kpIsolzait}
\end{center}
\end{figure}

Haragus and Pego \cite{hape} showed that the Zaitsev solution is an 
attractor for periodic (in $y$) perturbations of the line soliton. 
As we will illustrate below, the Zaitsev solution itself is unstable 
against periodic perturbations. 
To check the results of \cite{hape} one has to consider small 
periodic perturbations which travel with the wave. Above we 
do not have this restriction and find typically the formation of lumps 
that travel at higher speeds than the perturbed solution, thus 
violating the traveling wave condition. To get closer to the 
situation in \cite{hape}, we consider initial data for the exact KdV 
soliton and a small perturbation given by the Zaitsev solution $u_{Zait}$ 
(\ref{Za}) of the form 
$u_{0}(x)=u_{sol}(x+L_{x},0)+0.1u_{Zait}(x+L_{x},y,0)$ with 
$\alpha=1$, $\beta=0.5$, $L_{x} = 10$, $L_{y} = 5/\delta$. The 
relative mass conservation $\Delta$ is of the order $10^{-7}$ and 
smaller until the lumps appear, and decreases then to $10^{-2}$. 
As can be seen 
in Fig.~\ref{kpIsolzait}, the perturbation and the soliton 
essentially travel at the same speed for some time. One finds that 
the solution oscillates in this phase between the soliton and 
presumably the 
Zaitsev solution, until the instability against lump formation sets in.

\begin{figure}
[!htbp]
\begin{center}
\includegraphics[width=0.8\textwidth]{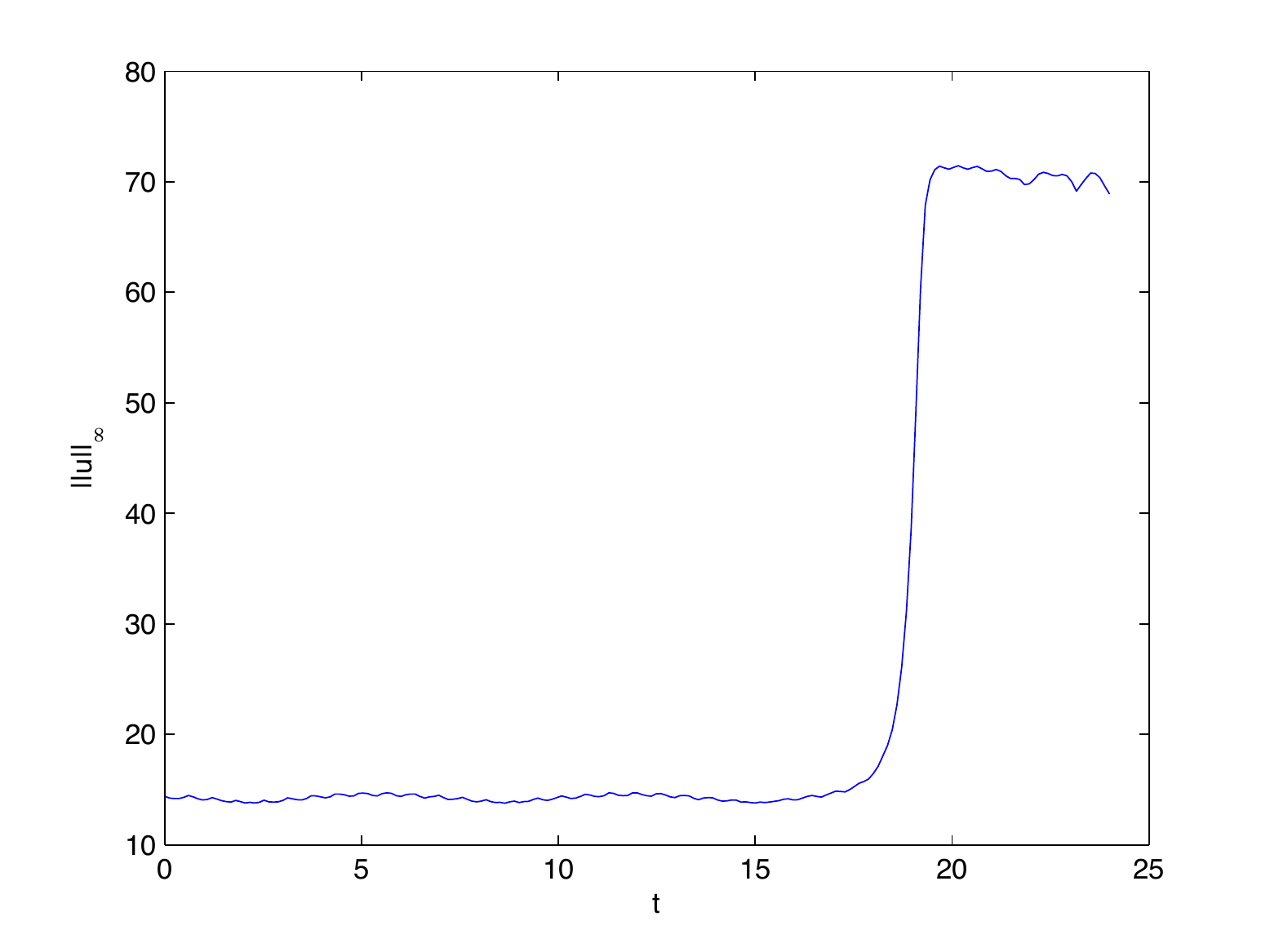} 
\caption{$L^{\infty}$-norm of the solution to the  KP I equation  of 
Fig.~\ref{kpIsolzait}.}
\label{kpIsolzaitmax}
\end{center}
\end{figure}

These oscillations of the traveling wave profile between soliton and 
Zaitsev solution is even more visible if one studies the maximal 
amplitude of the solution in Fig.~\ref{kpIsolzait} in dependence of 
time as shown in Fig.~\ref{kpIsolzaitmax}. Three maxima can be 
recognized before the onset of lump formation. A similar behavior will be 
seen in the next subsection for a perturbation of the Zaitsev 
solution. In contrast to the situation shown in Fig.~\ref{kpIperpL10} 
the solution here does not change if a larger $L_{x}$ is chosen. The 
lump formation starts at roughly the same time as before. 

\subsection{Perturbations of the Zaitsev solution}
The KP I equation admits solutions that are exponentially localized in 
one spatial direction and periodic in the other. The simplest such 
solution (\ref{Za}) was found by Zaitsev \cite{Z}.
The solution was generalized in \cite{muta89}.
We consider perturbations of the form
$u_{p}=6(x+L_{x}/2)\exp(-(x+L_{x}/2)^{2}-y^{2})$ which satisfy the 
constraint (\ref{constraint}). For the numerical experiments we choose 
$\alpha=1$, $\beta=0.5$, $L_{x} = 10$, $L_{y} = 5/\delta$, 
$N_{x}=2^{9}$, $N_{y}=2^{8}$, $N_{t}=10000$. With these settings 
we obtain a numerical mass conservation $\Delta < 10^{-5}$. Notice 
that the unperturbed Zaitsev solution can be numerically propagated 
with the same precision as the line soliton, see \cite{KR}.

\begin{figure}
[!htbp]
\begin{center}
\includegraphics[width=\textwidth]{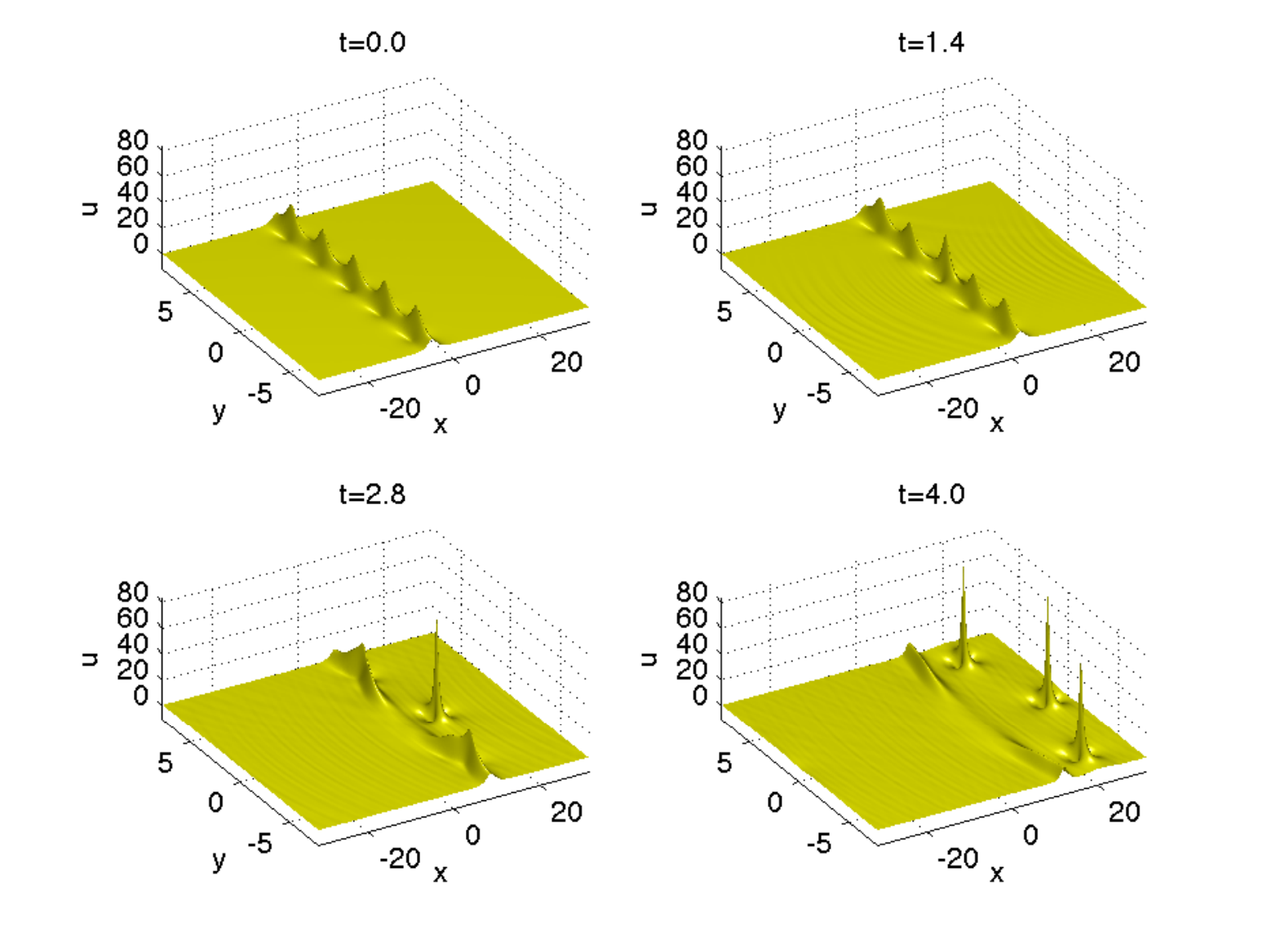} 
\caption{Solution to the  KP I equation  for initial 
data given by the Zaitsev solution $u_{Zait}(x+L_{x}/2,y,0)$ plus a perturbation 
$u_{p}=6(x+L_{x}/2)\exp(-(x+L_{x}/2)^{2}-y^{2})$ for various values of $t$.}
\label{kpzaitp}
\end{center}
\end{figure}

For initial data given by the Zaitsev solution centered at $-L_{x}/2$ 
plus $u_{p}=6(x+L_{x}/2)\exp(-(x+L_{x}/2)^{2}-y^{2})$, 
we find the behavior shown in Fig.~\ref{kpzaitp}. One 
can see that the initial perturbation develops into a lump traveling 
faster than the remaining Zaitsev solution. The latter develops in 
the following further lumps in a similar way as the perturbed KdV 
soliton. 

\begin{figure}
[!htbp]
\begin{center}
\includegraphics[width=\textwidth]{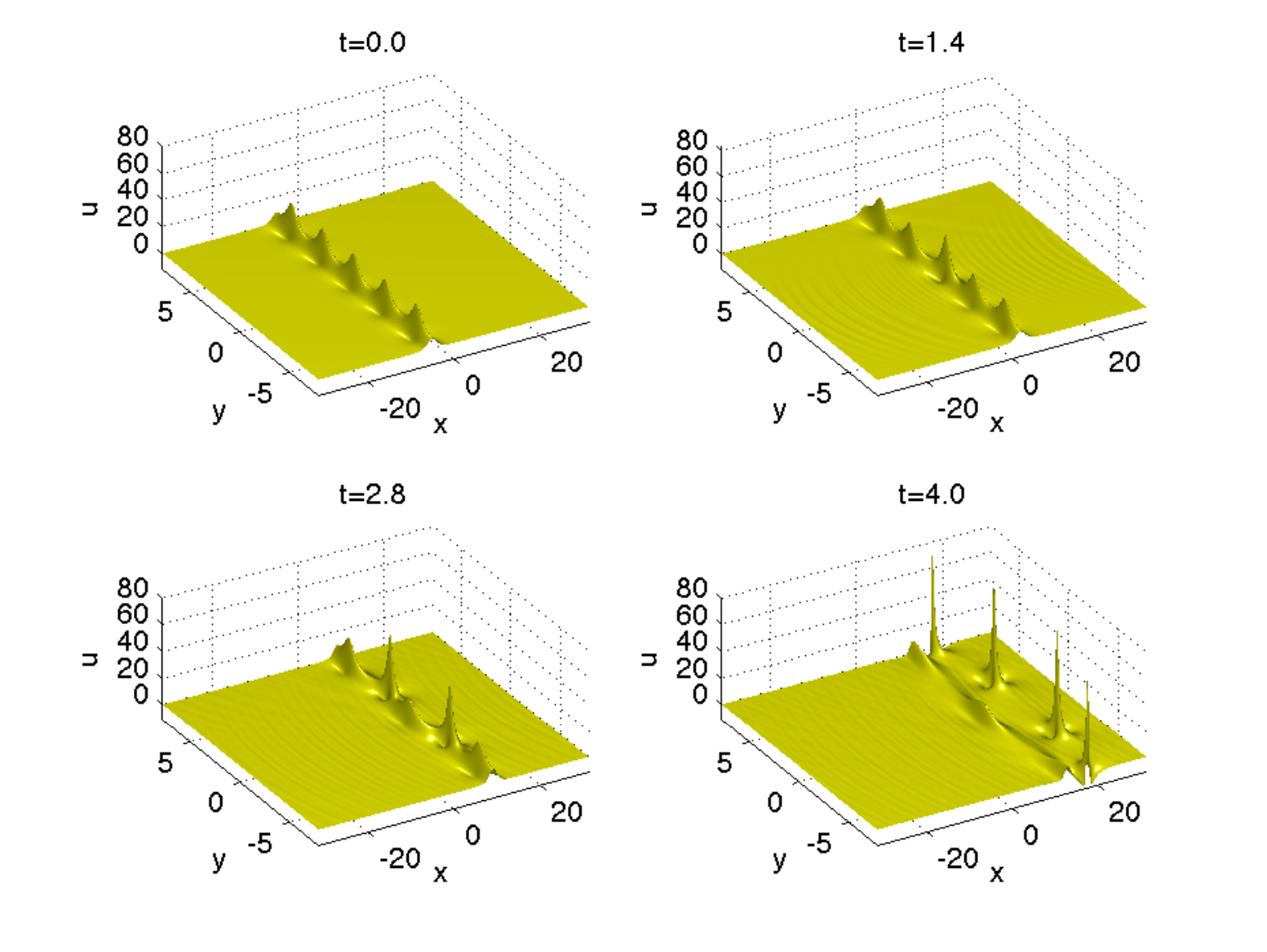} 
\caption{Solution to the  KP I equation  for initial 
data given by the Zaitsev solution $u_{Zait}(x+L_{x}/2,y,0)$ plus a perturbation 
$u_{p}=-6(x+L_{x}/2)\exp(-(x+L_{x}/2)^{2}-y^{2})$ for various values of $t$.}
\label{kpzaitm}
\end{center}
\end{figure}

Thus it appears that the Zaitsev solution is unstable as the line 
soliton against the formation of lumps. The precise pattern of the 
lumps depends on the initial perturbation. In Fig.~\ref{kpzaitm} the 
same initial condition as in Fig.~\ref{kpzaitp} is considered, this 
time however with a change in sign in the perturbation. The 
initial perturbation develops two lumps with further lumps 
forming at later times. 

\begin{figure}
[!htbp]
\begin{center}
\includegraphics[width=\textwidth]{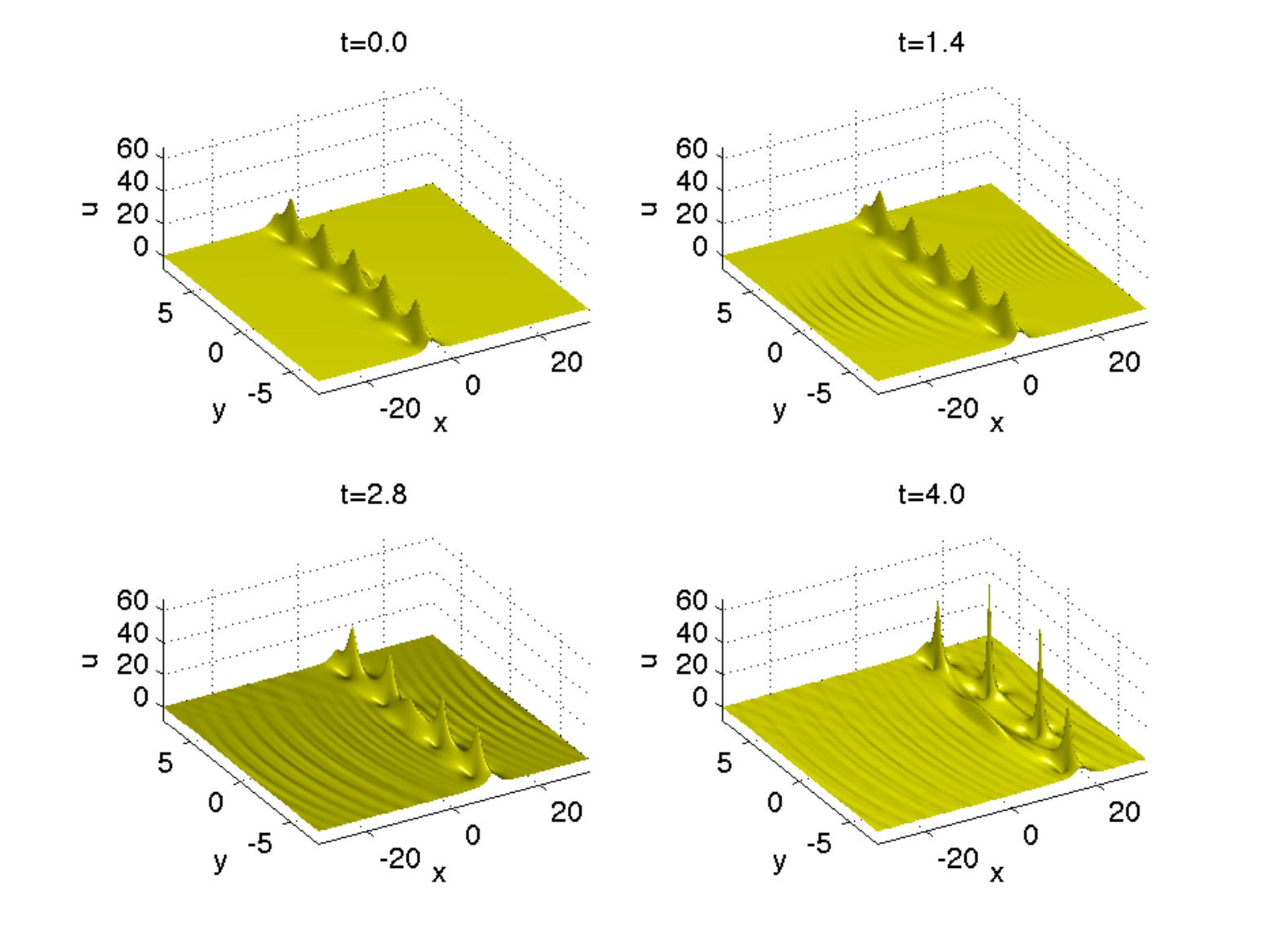} 
\caption{Solution to the  KP I equation  for initial 
data given by the Zaitsev solution plus a perturbation 
$u_{p}=6x\exp(-x^{2}-y^{2})$ for various values of $t$.}
\label{kpzaitmd}
\end{center}
\end{figure}

Notice that the Zaitsev solution is in some sense more unstable than 
the line soliton. First the amplitudes of the considered perturbations are a factor of 
10 smaller than the amplitudes of the solution, and still the decay to 
lumps is almost immediate. More importantly the Zaitsev solution is also 
unstable against a displaced perturbation 
($u_{p}=6x\exp(-x^{2}-y^{2})$) as can be seen in Fig.~\ref{kpzaitmd}. 
As in Fig.~\ref{kpIsolperd} for the perturbed KdV soliton, the initial perturbation develops 
dispersive oscillations and tails, but the former almost immediately destroy 
the Zaitsev solution and leads to the formation of lumps, in contrast 
to the situation for the KdV soliton. 

\begin{figure}
[!htbp]
\begin{center}
\includegraphics[width=\textwidth]{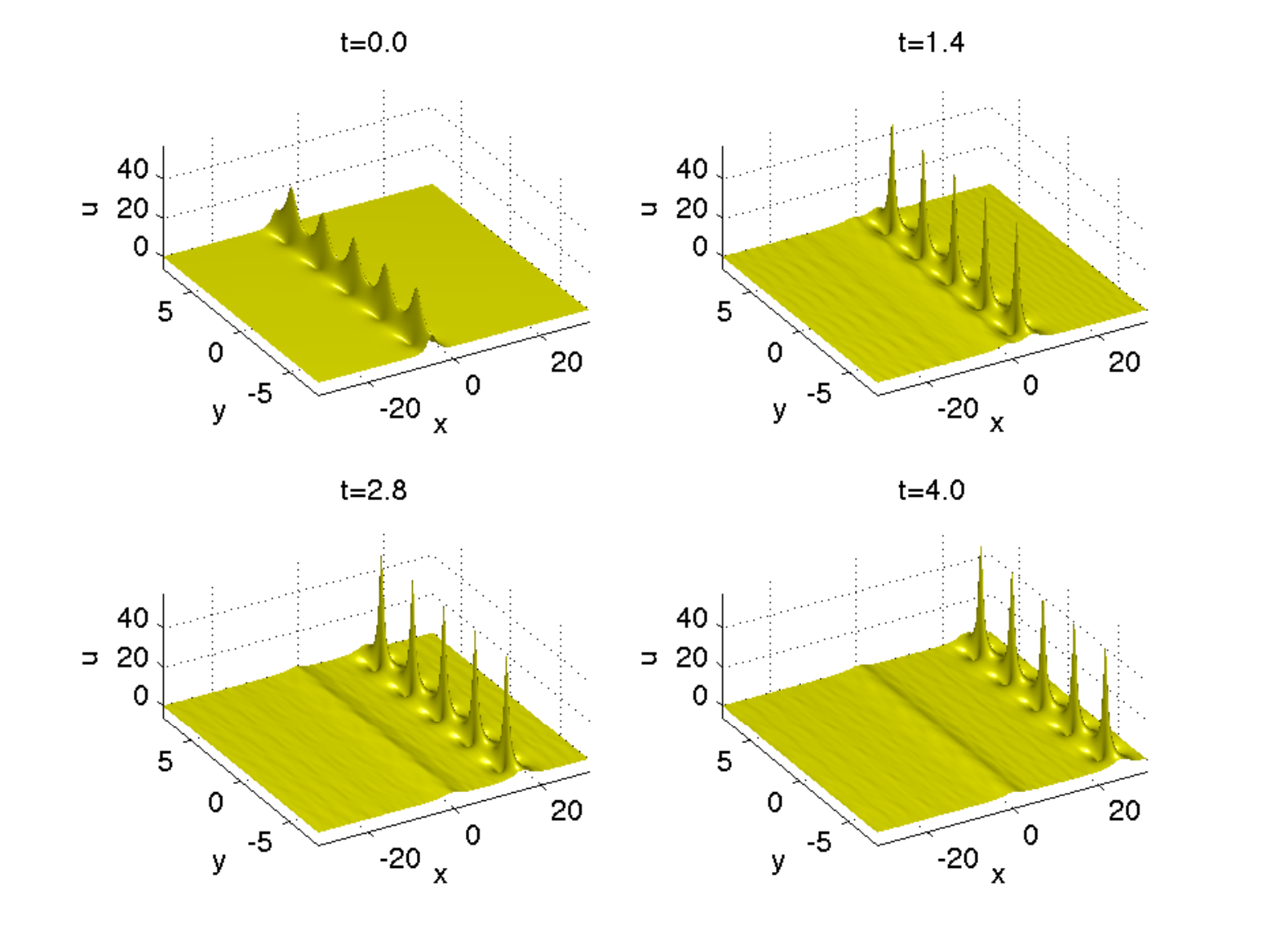} 
\caption{Solution to the  KP I equation  for initial 
data given by the Zaitsev solution plus a perturbation with the same 
period in $y$, $u_{0}(x)=1.1u_{Zait}(x+L_{x}/2,y,0)$ for various values of $t$.}
\label{kpzaitper}
\end{center}
\end{figure}

\begin{figure}
[!htbp]
\begin{center}
\includegraphics[width=0.8\textwidth]{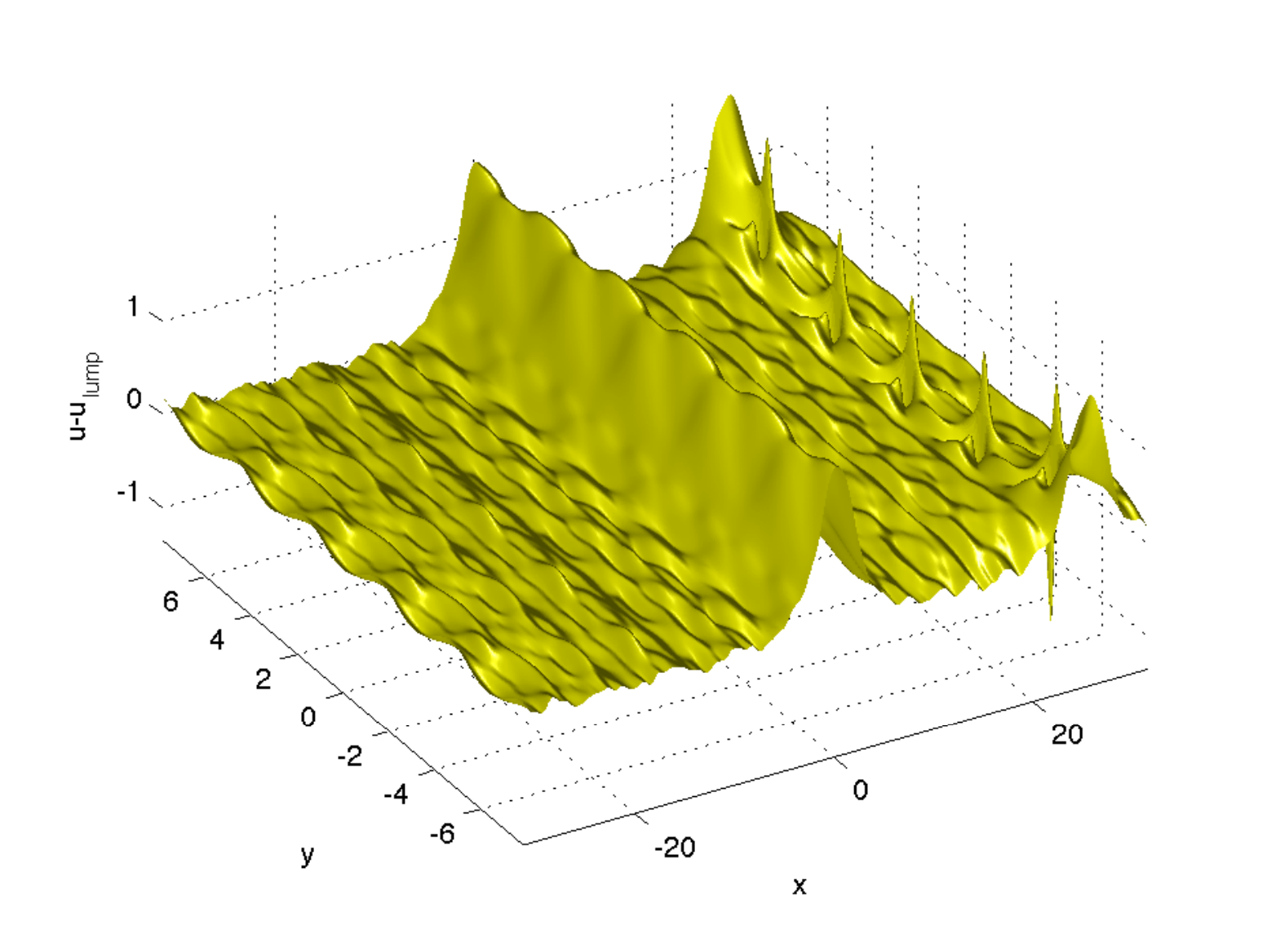} 
\caption{Difference of the solution to the  KP I equation  in 
Fig.~\ref{kpzaitper} for $t=4$ and five lump solutions fitted at the 
peaks.}
\label{kpzaitperml}
\end{center}
\end{figure}

The Zaitsev solution also appears to be unstable against lump 
formation for perturbations with the same period in $y$. In 
Fig.~\ref{kpzaitper} we show the time evolution of initial data from 
the Zaitsev solution as before, but with an amplitude multiplied by a 
factor of $1.1$. This corresponds to a small perturbation of exactly 
the same period as the Zaitsev solution. 
It can be seen that each of 
the maxima of the initial data develops a lump-like structure. 
If we neglect the precise form of a multi-lump solution and just 
subtract a single lump at each of the final peaks in 
Fig.~\ref{kpzaitper}, we obtain Fig.~\ref{kpzaitperml}. This suggests 
that the formed peaks will in fact develop asymptotically into a 
periodic array of lumps.

\begin{figure}
[!htbp]
\begin{center}
\includegraphics[width=\textwidth]{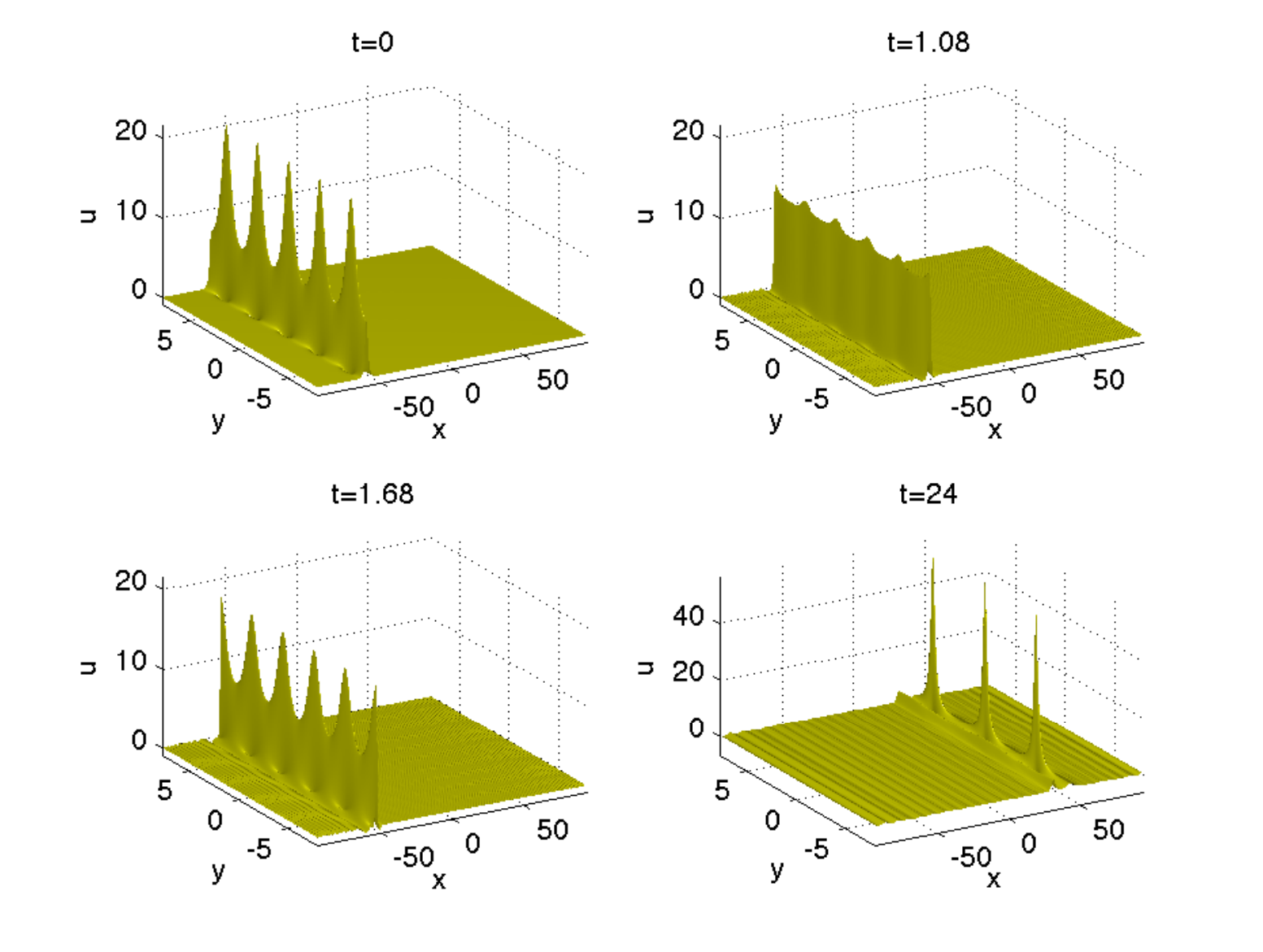} 
\caption{Solution to the  KP I equation  for initial 
data given by the Zaitsev solution plus a perturbation with the same 
period in $y$, $u_{0}(x)=0.9u_{Zait}(x+L_{x}/2,y,0)$ for various values of $t$.}
\label{kpzaitpersm}
\end{center}
\end{figure}

This is however only the case if the Zaitsev solution is multiplied 
with a factor greater than 1, i.e., if the humps are enlarged as 
above. In this case they seem to evolve into  lumps. If instead the 
Zaitsev solution is multiplied by a factor smaller than 1, 
e.g.~$0.9$, one gets for the time evolution of such data 
Fig.~\ref{kpzaitpersm}. The computation is carried out for 
$L_{x}=30$, $L_{y}=2.5$ with 
$N_{x}=2^{11}$, $N_{y}=2^{9}$ and $N_{t}=20000$ time steps to yield a 
$\Delta\sim 3*10^{-4}$. In this case one again observes oscillations 
between the Zaitsev solution  and presumably the KdV soliton, until finally lump formation sets in. 

\begin{figure}
[!htbp]
\begin{center}
\includegraphics[width=0.8\textwidth]{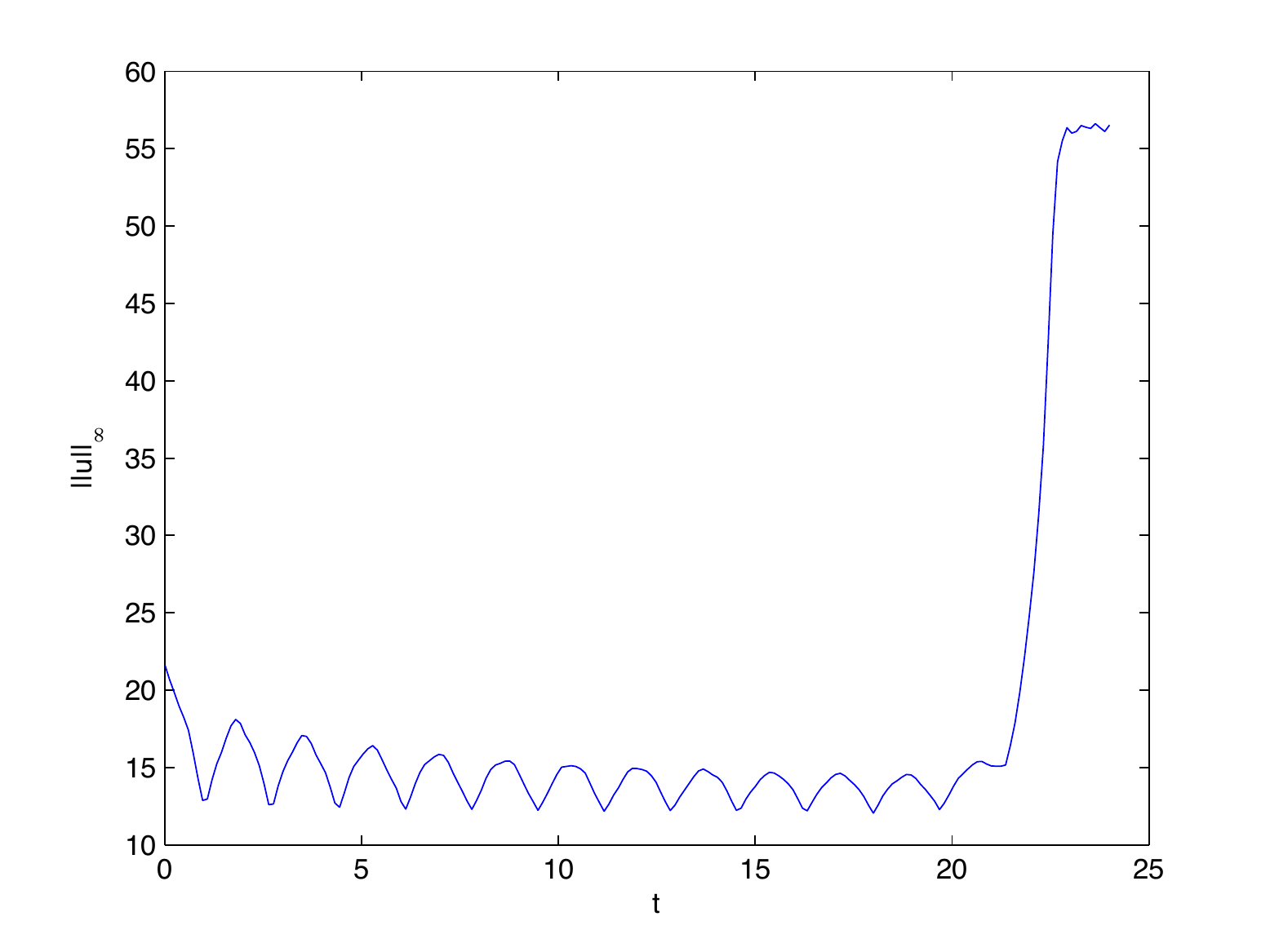} 
\caption{$L^{\infty}$-norm of the solution $u$ of 
Fig.~\ref{kpzaitpersm} in dependence of time.}
\label{kpzaitpersmmax}
\end{center}
\end{figure}

\begin{figure}
[!htbp]
\begin{center}
\includegraphics[width=0.8\textwidth]{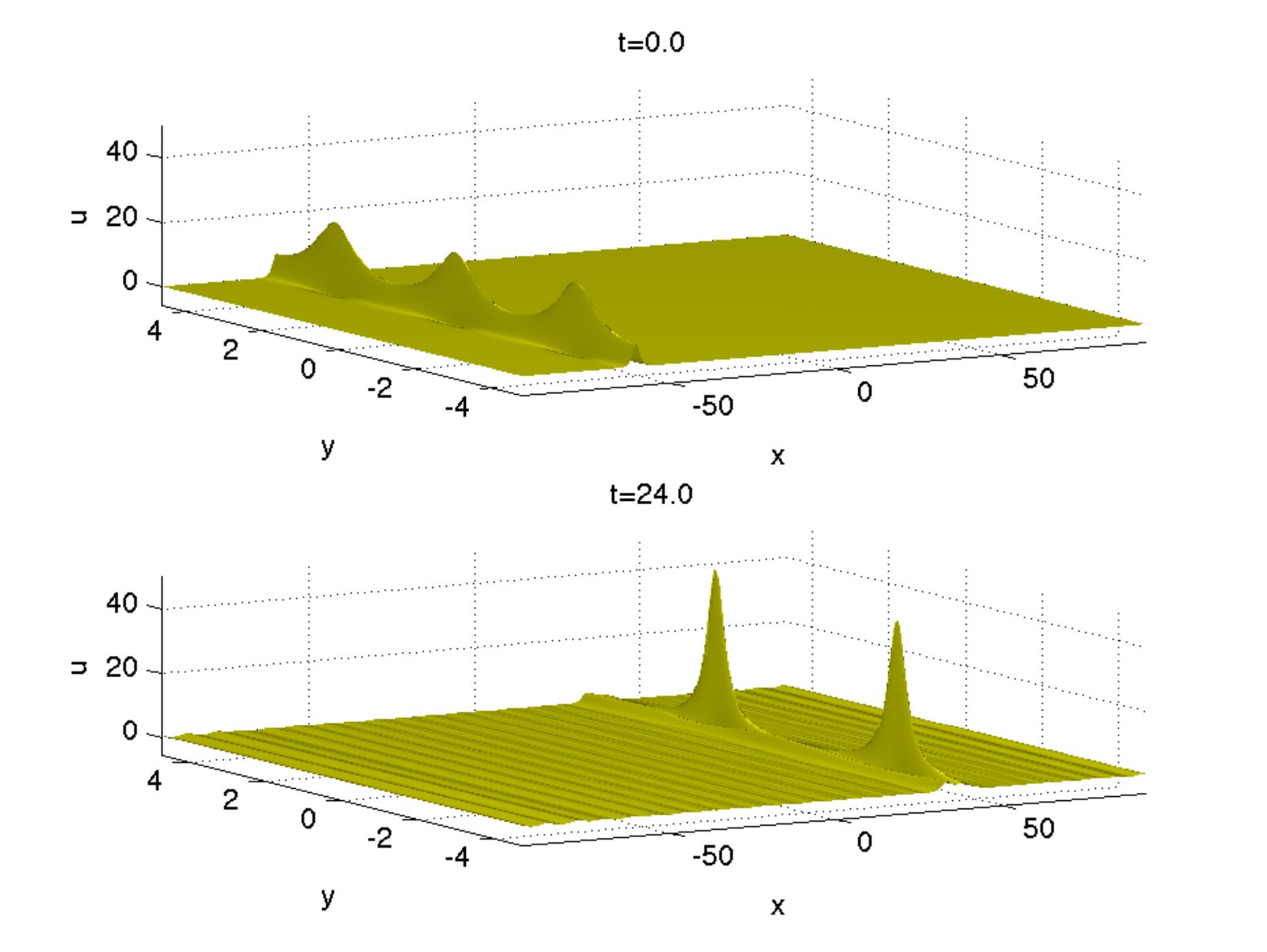} 
\caption{Solution to the KP I equation for the situation of 
Fig.~\ref{kpzaitpersm} with a smaller $L_{y}=1.5$.}
\label{kpzaitpersmally3}
\end{center}
\end{figure}

\begin{figure}
[!htbp]
\begin{center}
\includegraphics[width=0.8\textwidth]{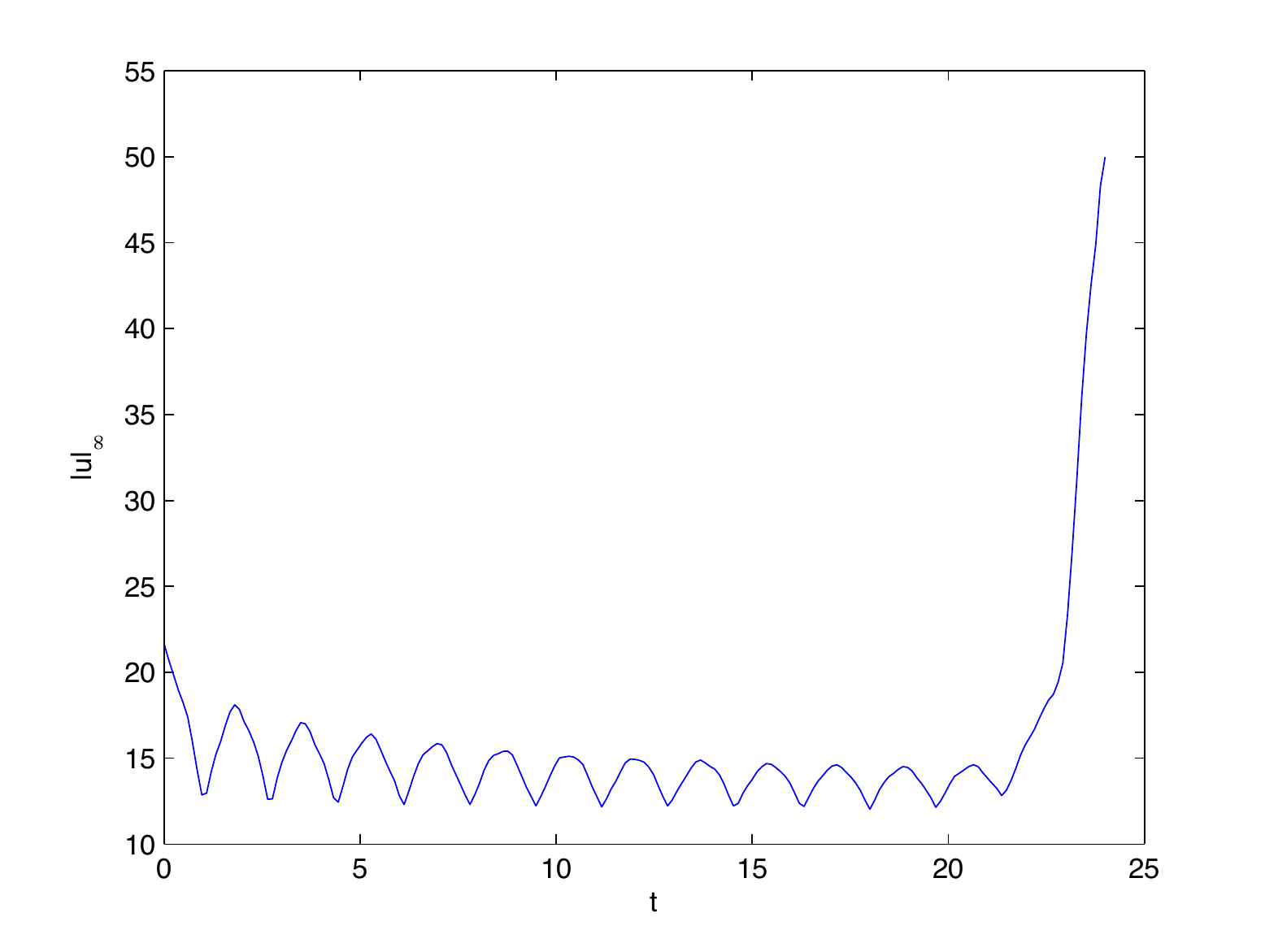} 
\caption{$L^{\infty}$-norm of the solution $u$ of 
Fig.~\ref{kpzaitpersmally3} in dependence of time.}
\label{kpzaitpersmally3max}
\end{center}
\end{figure}

In Fig.~\ref{kpzaitpersmmax} it can be seen that this oscillatory 
phase can be quite long until lump formation starts. It appears that 
this oscillatory state is thus meta-stable in some sense, and the 
pattern does not change within the given precision 
if the code is run with $N_{t}=30000$.
The 
mechanism of the final decay to lumps is not clear. If one runs the 
code with the same settings on a larger domain ($L_{x}=35$), lump 
formation sets in at roughly the same time. The same is true on a 
slightly reduced domain ($L_{x}=25$). Thus in contrast to the 
perturbation of the KdV soliton studied in Fig.~\ref{kpIperp3}, the 
proximity of the domain end in $x$-direction does not appear to be 
decisive. Instead it seems that numerical inaccuracies in the 
periodicity in $y$ trigger the lump formation. In fact if the same 
situation as in Fig.~\ref{kpzaitpersm} is studied with a smaller 
$L_{y}=1.5$, the initial data show 3 maxima, which, however, develop 
into just two lumps as can be seen in Fig.~\ref{kpzaitpersmally3}. 
This indicates some symmetry breaking, presumably due to effects from 
the boundary of the computational domain in $y$-direction. The computation 
is carried out with $N_{x}=2^{11}$, $N_{y}=2^{8}$, 
$N_{t}=20000$ with a final $\Delta\sim 2*10^{-5}$. It does not change 
within the expected limit if a smaller resolution in space or time is 
chosen.
The $L^{\infty}$-norm of the solution in Fig.~\ref{kpzaitpersmally3} 
can be seen in Fig.~\ref{kpzaitpersmally3max}. In comparison to 
Fig.~\ref{kpzaitpersmmax} it can be recognized that the lump 
formation starts at a slightly later time. Thus we cannot conclude 
what would happen if the same initial data would be treated on 
$\mathbb{R}\times \mathbb{T}$ instead of $\mathbb{T}^{2}$, but it 
seems that the instabilities here could be due to a sensitivity to 
inaccuracies in the $y$-periodicity.

\subsection{Numerical study of blowup in generalized KP equations}
In this subsection we will study the time evolution of initial data 
$u_{0}(x)$ for the generalized KP equation (\ref{genpKP}).
It is known \cite{MST2} that the solution will blow up for 
initial data with negative energy for $p\geq 4/3$. 

In the following we will consider the initial data
\begin{equation}
    u_{0}(x)=6\partial_{xx}\exp(-\alpha(x^{2}+y^{2})),\quad \alpha\in 
    \mathbb{R}^{+}
    \label{u0}.
\end{equation}
These data obviously satisfy the constraint (\ref{constraint}).

\begin{figure}
[!htbp]
\begin{center}
\includegraphics[width=0.8\textwidth]{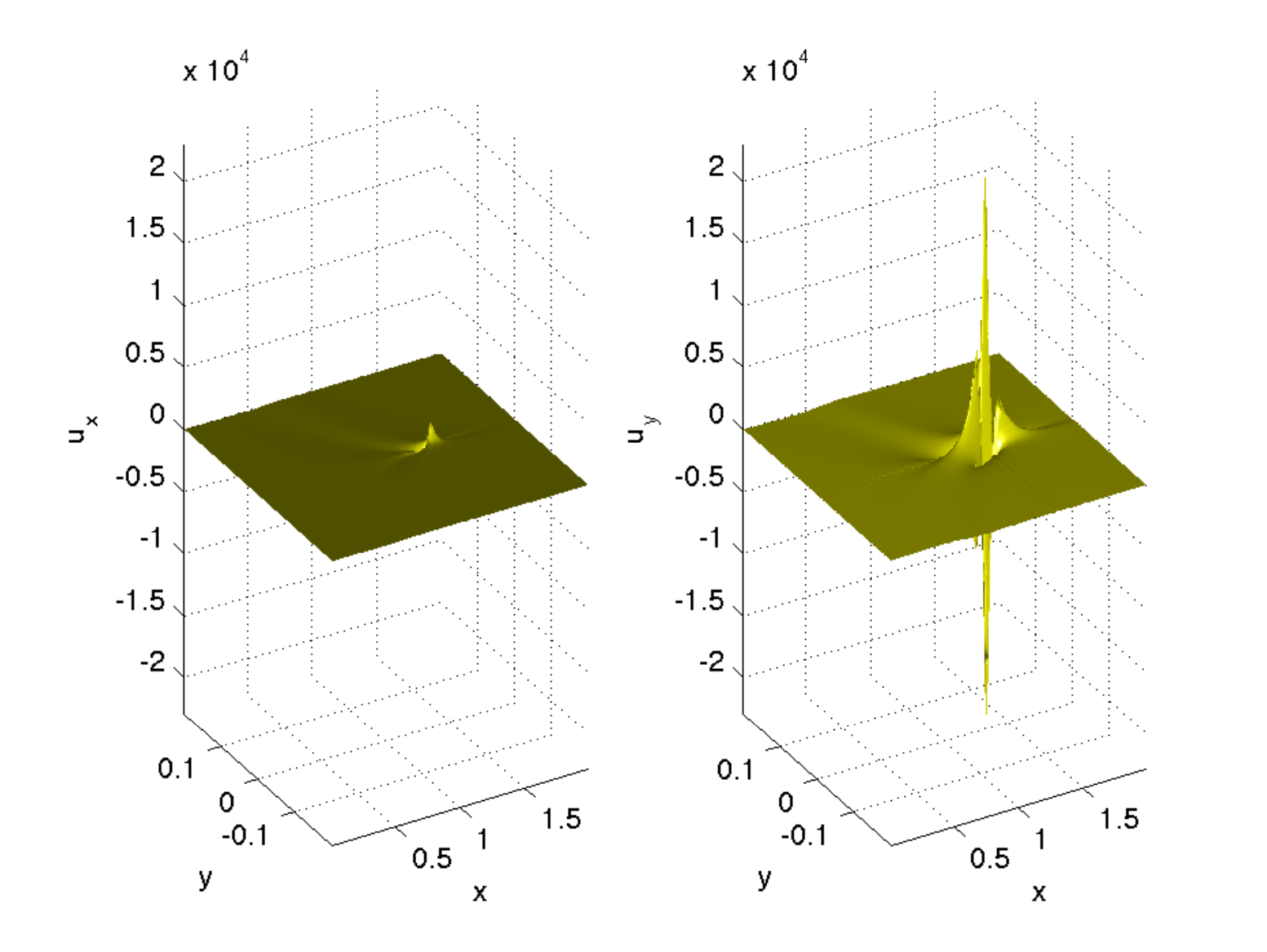} 
\caption{Gradient of the solution to the generalized KP I equation 
($p=2$) for the 
initial data $u_{0}(x,y)=6\partial_{xx}\exp(-x^{2}-y^{2})$ 
for $t=0.0473$.}
\label{kp3duc}
\end{center}
\end{figure}

For the supercritical $p=2$ we find that these data  in fact imply a 
negative energy for $\alpha=1$. The computation is carried out with $L_{x}=5$, 
$L_{y}=2$, $N_{x}=2^{11}$, $N_{y}=2^{13}$ modes and $N_{t}=5000$ time 
steps. The code breaks at a finite time $t_{c}$ since the solution appears to 
develop a singularity. It is stopped once the quantity $\Delta$ 
measuring numerical mass 
conservation becomes larger than $10^{-4}$. It can be checked that the 
Fourier coefficients at this time ($t_{c}\sim0.0473$ for the considered 
initial data) still decrease by 5 orders of magnitude. This implies 
that the solution is still accurate to better than plotting accuracy 
at this time. 
\begin{figure}
[!htbp]
\begin{center}
\includegraphics[width=\textwidth]{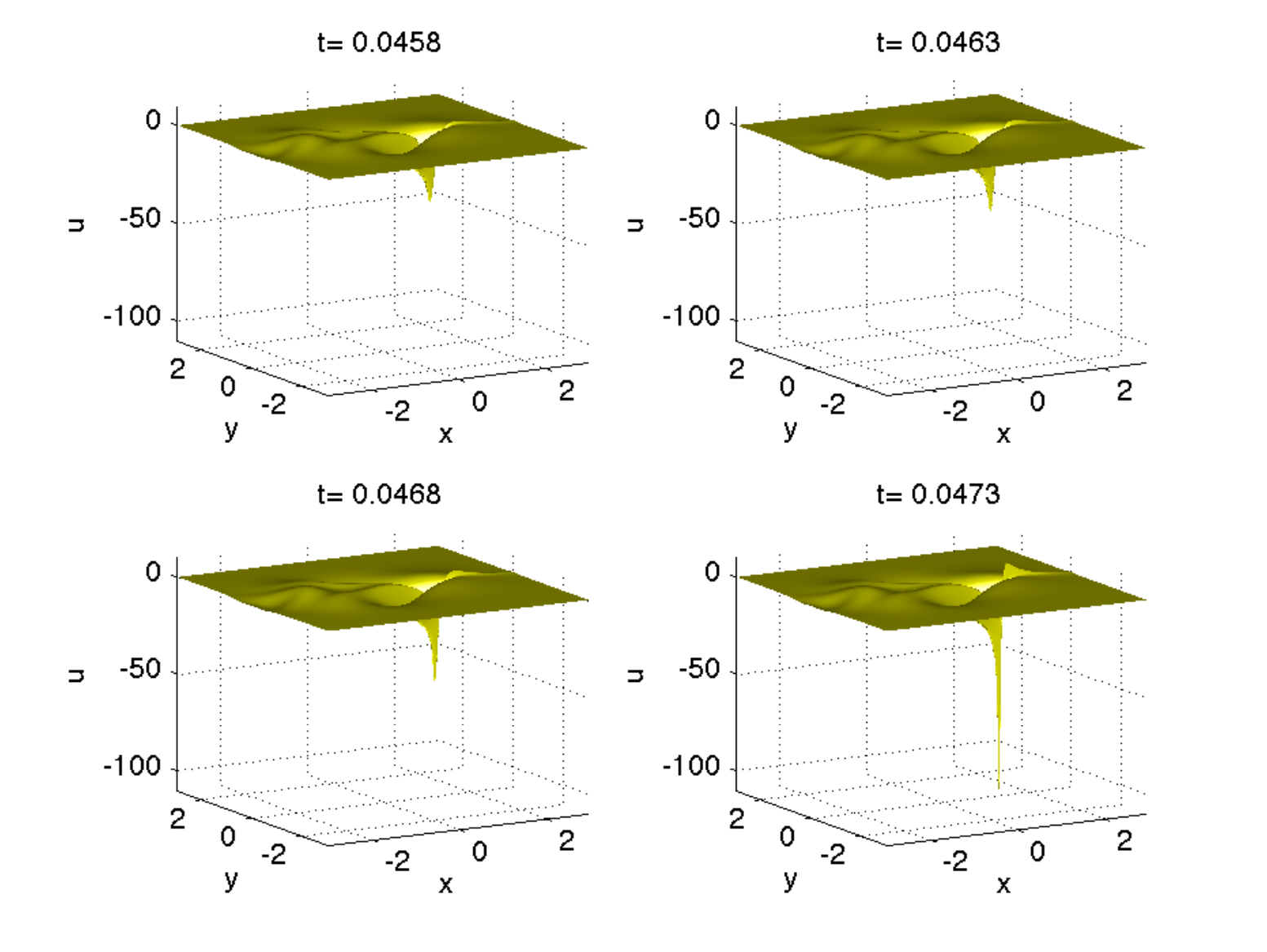} 
\caption{Solution to the generalized KP I equation ($p=2$) for the 
initial data $u_{0}(x,y)=6\partial_{xx}\exp(-(x^{2}+y^{2}))$ 
for various values of $t$.}
\label{kp3u}
\end{center}
\end{figure}
It can be seen that the initial minimum develops apparently into a 
singularity. 
The precise formation of the latter is not isotropic in $x$ and 
$y$ as shown by the gradient of $u$ close to the critical 
time in Fig.~\ref{kp3duc}. It can be seen that the $y$-derivative 
near the catastrophe is more than an order of magnitude bigger than 
the $x$-derivative. Thus it is foremost the $y$-derivative of the 
solution that diverges which then leads to a divergence of the 
solution itself. The much stronger gradients in $y$-direction are the 
reason why a higher resolution has to be chosen in $y$ than in 
$x$. With the used parameters we ensure essentially the same 
resolution in both spatial directions.

\begin{figure}
[!htbp]
\begin{center}
\includegraphics[width=0.8\textwidth]{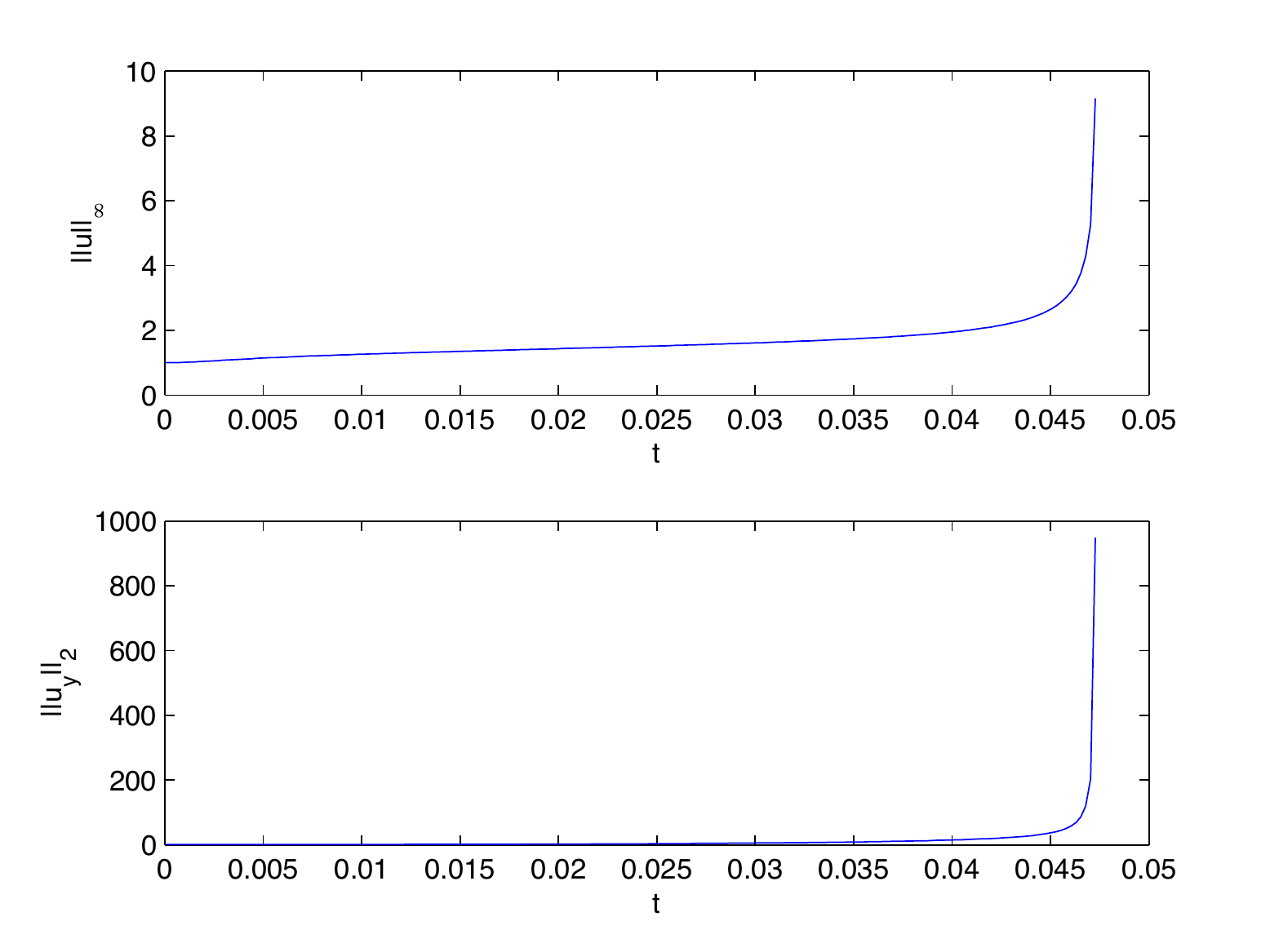} 
\caption{The $L^{\infty}$-norm of the solution to the KP I 
equation ($p=2$) for the 
initial data $u_{0}(x,y)=6\partial_{xx}\exp(-x^{2}-y^{2})$ and the 
$L^{2}$-norm of $u_{y}$ normalized to 1 for $t=0$ in dependence of time.
}
\label{kp3norm}
\end{center}
\end{figure}

This behavior indicates that it is sensible to trace the  time dependence of 
the $L^{\infty}$-norm of $u$ and the $L^{2}$-norm of $u_{y}$ which 
are shown in Fig.~\ref{kp3norm}. Both norms diverge at the critical 
time. It can be seen that both norms are monotonically increasing, 
but that there is a steep divergence at the critical point. It is 
also obvious that the norm of $u_{y}$ diverges much more rapidly than $u$. 
From the obtained data one cannot decide the exact character of the 
divergence near the critical time (exponential, algebraic).

\begin{figure}
[!htbp]
\begin{center}
\includegraphics[width=\textwidth]{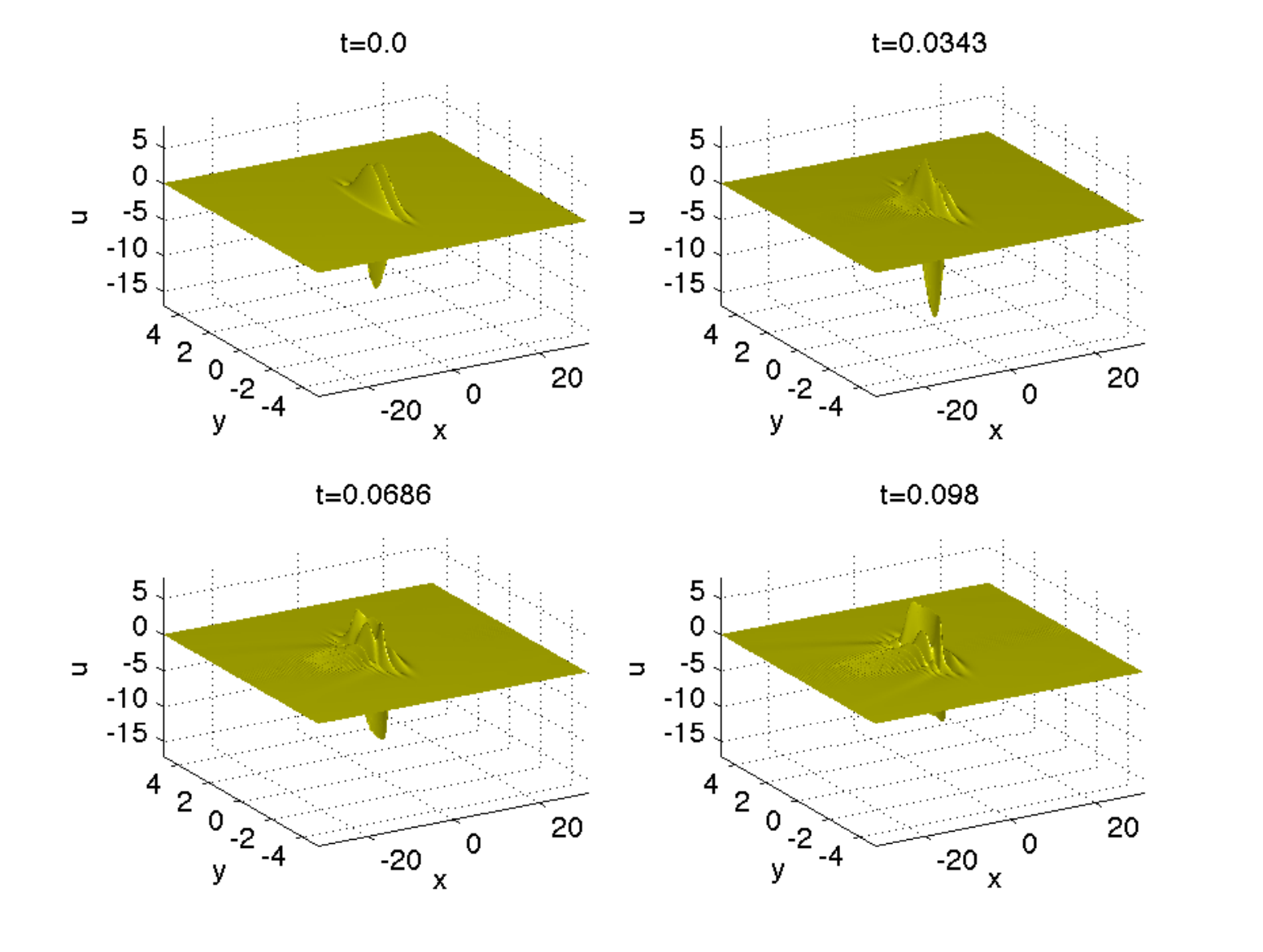} 
\caption{Solution to the generalized KP II equation ($p=2$) for the 
initial data $u_{0}(x,y)=6\partial_{xx}\exp(-x^{2}-y^{2})$ for 
various times.}
\label{kpII}
\end{center}
\end{figure}

For the generalized KP II equation  we find for the same 
$p$ and the same initial data that there is no indication for blowup. 
The time evolution of the solution can be seen in Fig.~\ref{kpII}. It 
can be seen that the typical oscillations in $x$-direction due to the 
Airy term in the KP equation appear, see \cite{klspma}. The 
algebraic fall off to infinity typical for KP 
solutions even for Schwartzian initial data is clearly visible in the 
form of tails to the left. 

\begin{figure}
[!htbp]
\begin{center}
\includegraphics[width=0.8\textwidth]{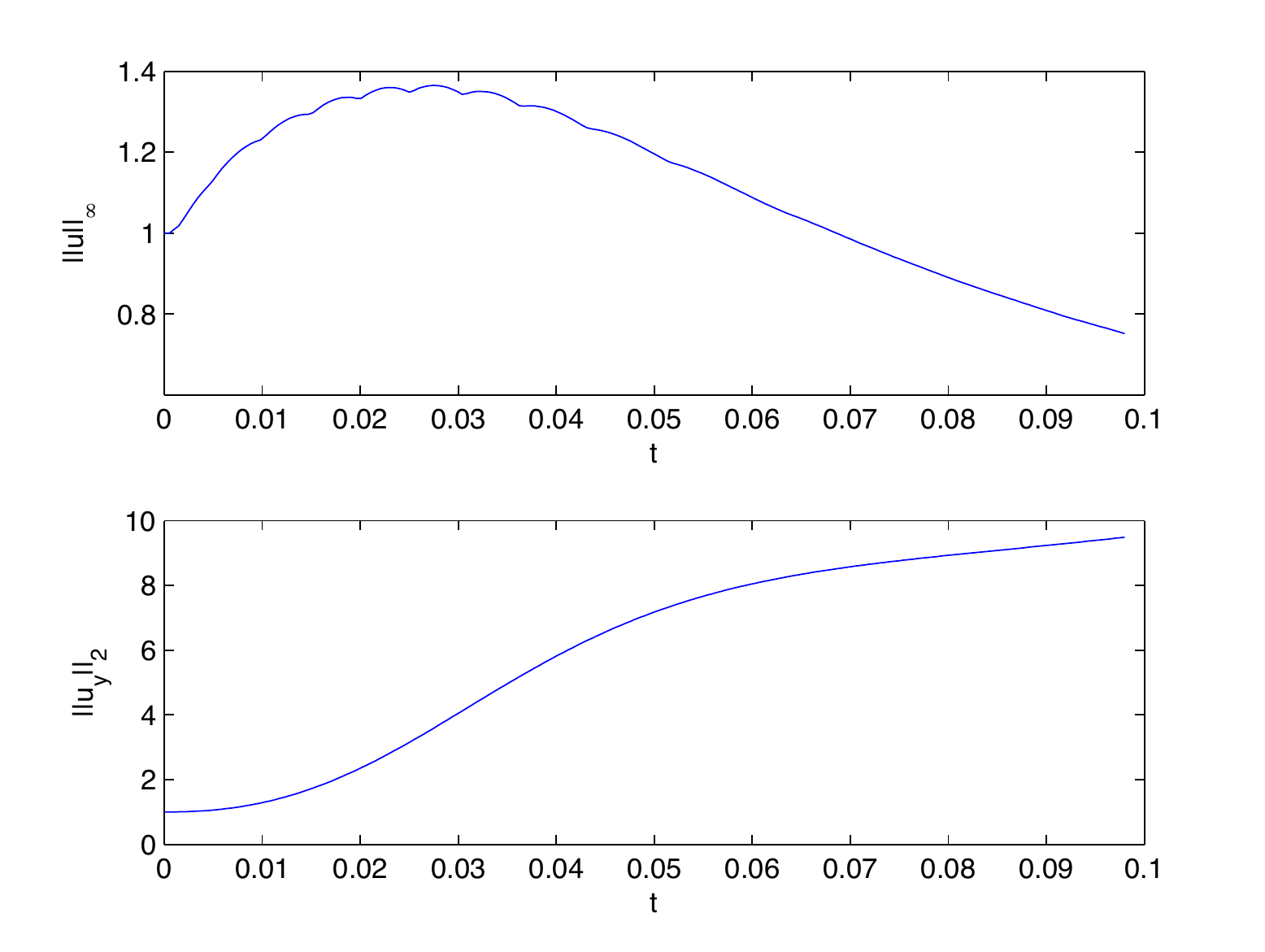} 
\caption{The $L^{\infty}$-norm of the solution to the generalized KP 
II equation ($p=2$) for the 
initial data $u_{0}(x,y)=6\partial_{xx}\exp(-x^{2}-y^{2})$ and the 
$L^{2}$-norm of $u_{y}$  normalized to 1 for $t=0$ in dependence of time.
}
\label{kpIInorm}
\end{center}
\end{figure}

In contrast to the KP I case, there is no indication for a blow up. In 
Fig.~\ref{kpIInorm} the $L^{\infty}$-norm of $u$ and of the 
$L^{2}$-norm of $u_{y}$ are shown. The former seems to decrease since 
the negative hump is essentially radiated away to infinity 
via the tails or to lead to the dispersive oscillations. The 
$L^{2}$-norm of $u_{y}$ appears to approach slowly a finite 
asymptotic value.

\begin{figure}
[!htbp]
\begin{center}
\includegraphics[width=0.8\textwidth]{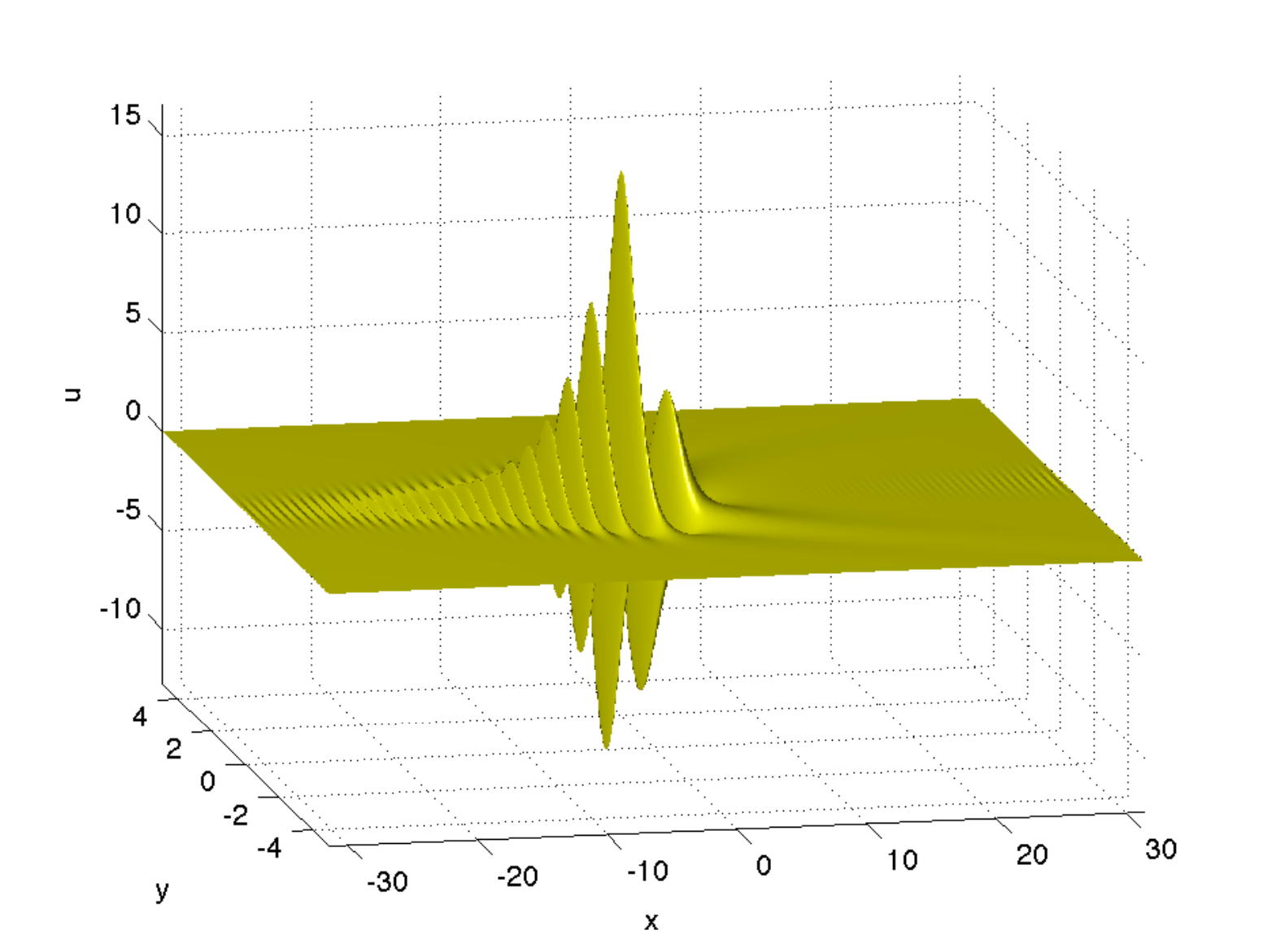} 
\caption{Solution to the  KP I equation ($p=1$) for the 
initial data $u_{0}(x,y)=12\partial_{xx}\exp(-x^{2}-y^{2})$ for 
$t=0.15$.}
\label{kpI}
\end{center}
\end{figure}

For the standard KP I equation ($p=1$) similar initial data lead as 
expected to a regular solution as can be seen in Fig.~\ref{kpI}. The 
energy for these initial data is positive. The 
solution just shows the characteristic oscillations and the tails. 
This computation and the ones below are carried out with 
$N_{x}=2^{10}$ and $N_{y}=2^{8}$ and $N_{t}=1000$ time steps. The Fourier 
coefficients decrease to $10^{-4}$ which ensures spatial 
resolution, and the relative mass conservation is numerically  
satisfied to better than $10^{-6}$. Thus with numerical precision we can exclude 
the formation of a singularity up the considered time $t=0.15$.
The $L^{\infty}$-norm of $u$ and the $L^{2}$-norm of $u_{y}$ both 
decrease for large times. Thus there is no indication of the appearance of a 
singularity as expected.
\begin{figure}
[!htbp]
\begin{center}
\includegraphics[width=0.8\textwidth]{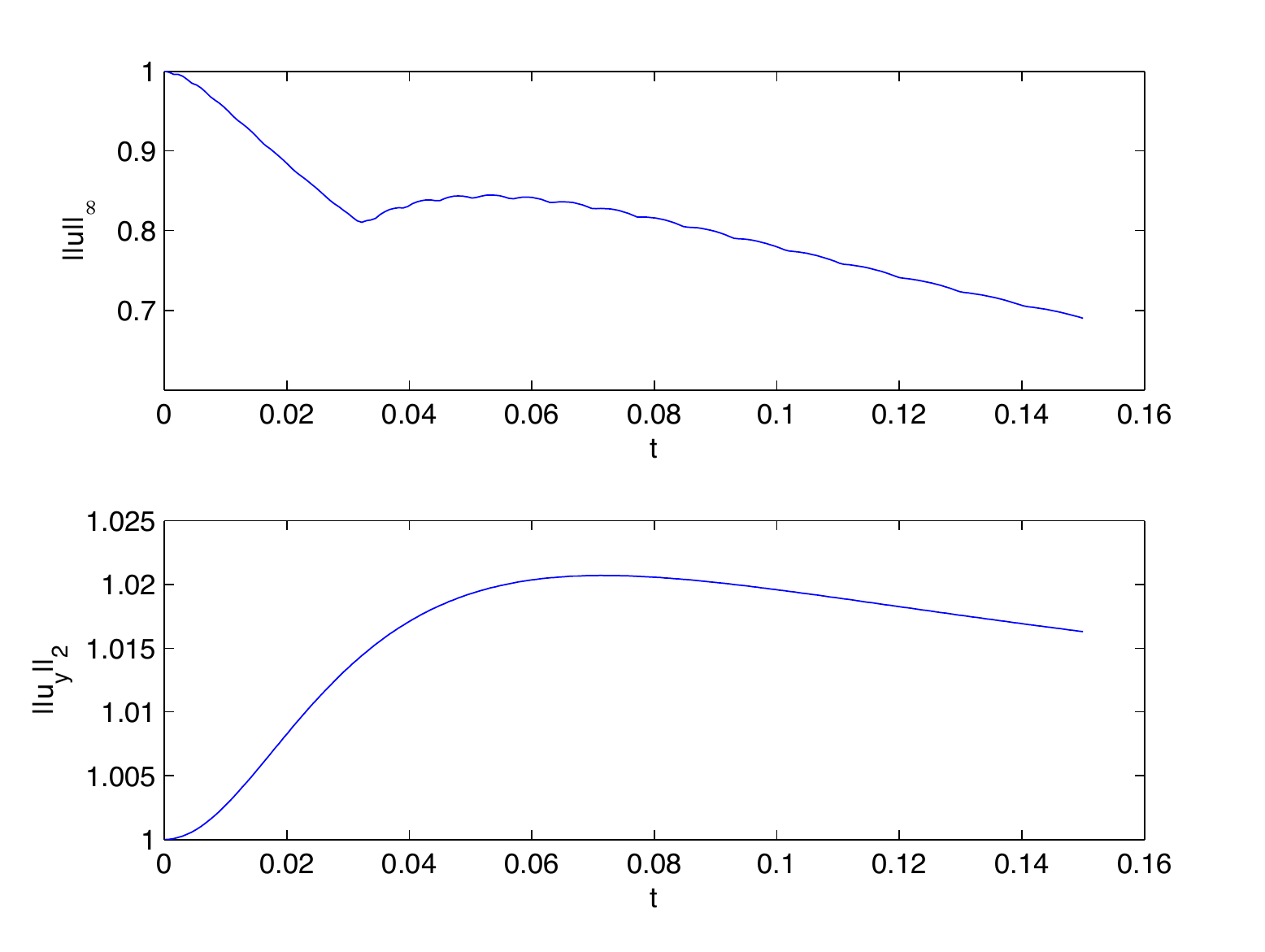} 
\caption{The $L^{\infty}$-norm of the solution to the KP I 
equation ($p=1$) for the 
initial data $u_{0}(x,y)=12\partial_{xx}\exp(-x^{2}-y^{2})$ and the 
$L^{2}$-norm of $u_{y}$ normalized to 1 for $t=0$ in dependence of time.
}
\label{kpInorm}
\end{center}
\end{figure}

\begin{figure}
[!htbp]
\begin{center}
\includegraphics[width=0.8\textwidth]{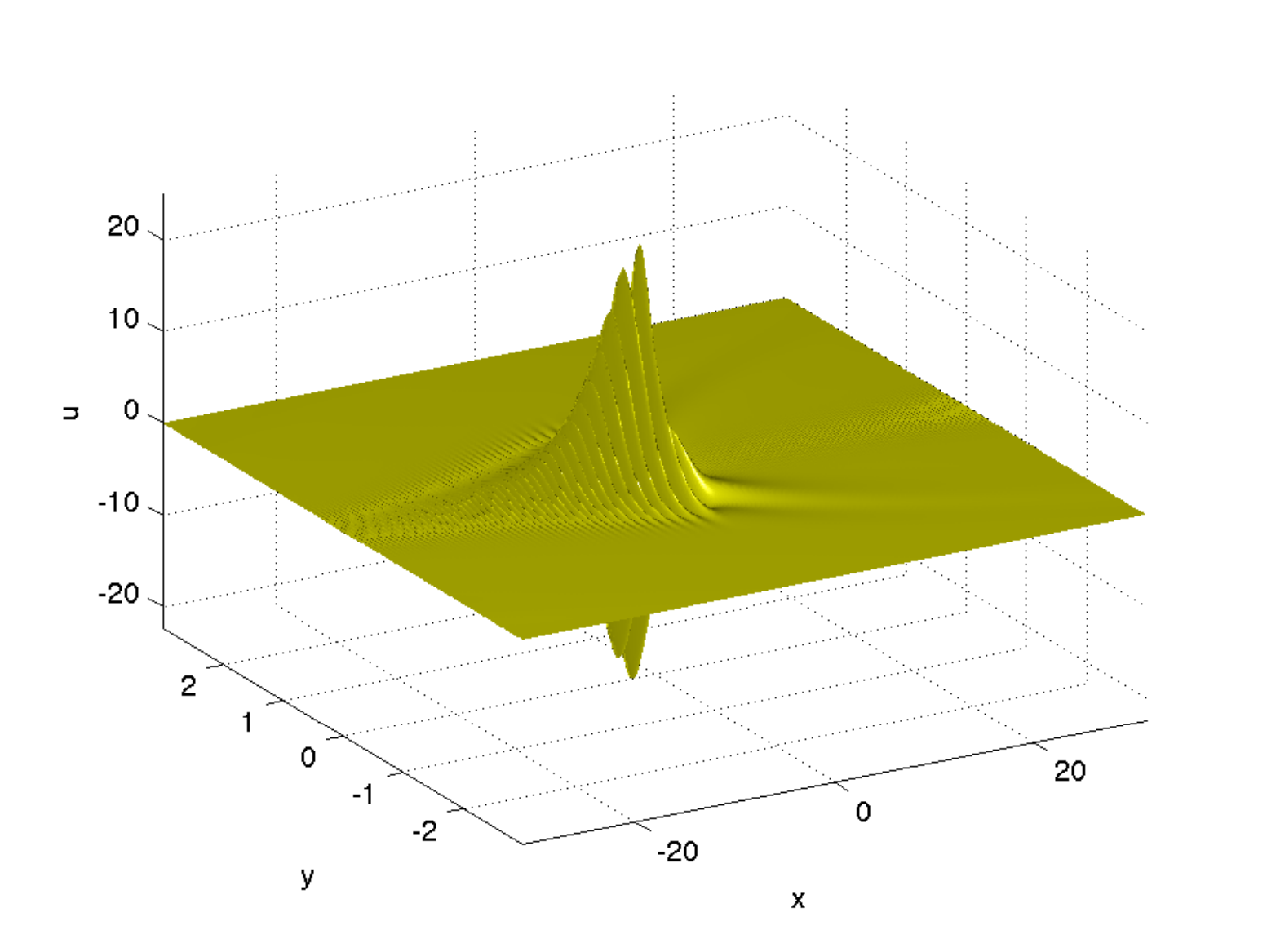} 
\caption{Solution to the generalized KP I equation with critical 
exponent $p=4/3$ for the 
initial data $u_{0}(x,y)=6\partial_{xx}\exp(-4(x^{2}+y^{2}))$ for 
$t=0.05$.}
\label{kpIc}
\end{center}
\end{figure}

For the critical exponent $p=4/3$ we obtain for the initial data 
(\ref{u0}) with $\alpha=4$, which leads to a negative energy, the
solution shown in Fig.~\ref{kpIc}. 

\begin{figure}
[!htbp]
\begin{center}
\includegraphics[width=0.8\textwidth]{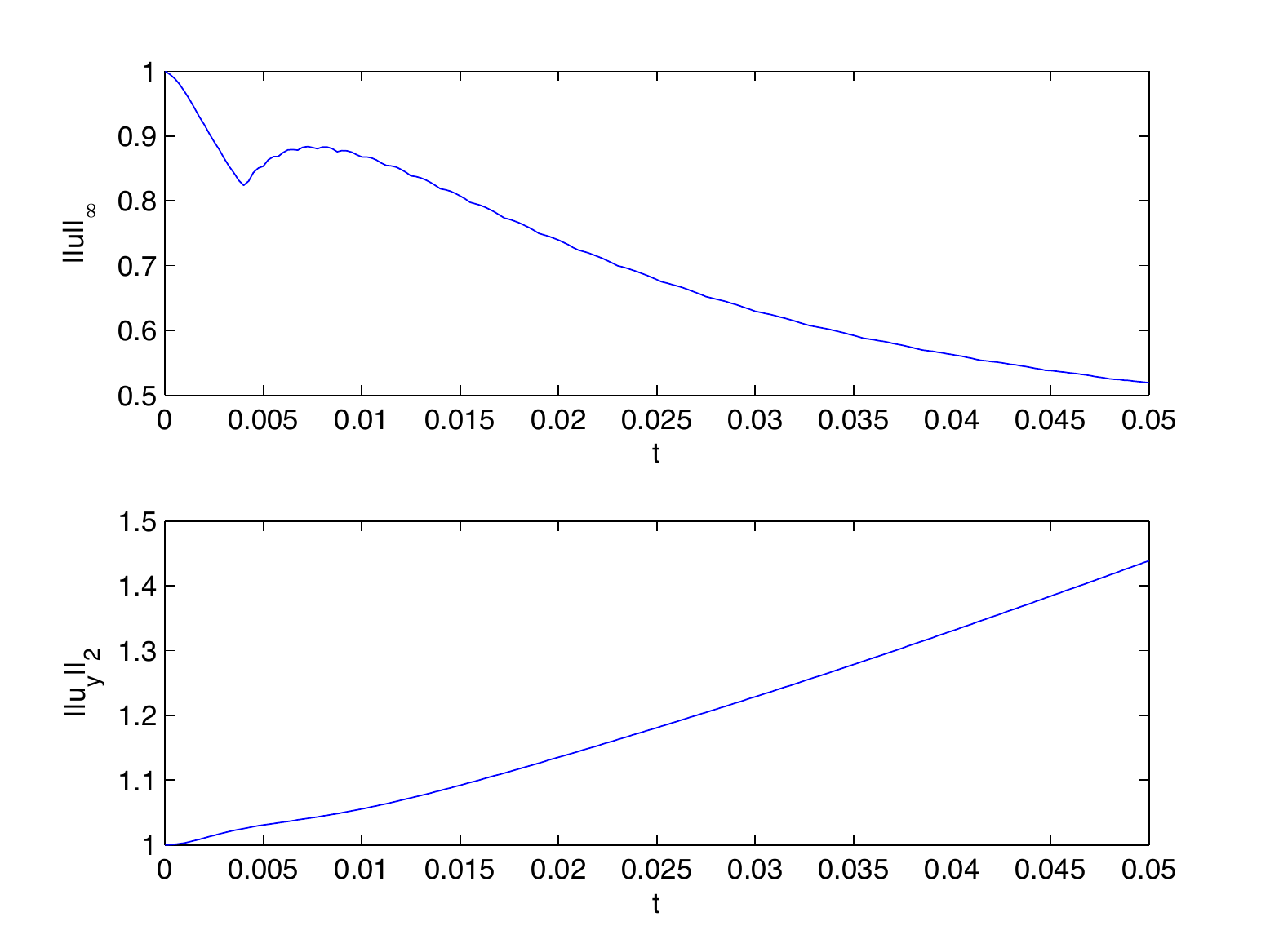} 
\caption{The $L^{\infty}$-norm of the solution to the generalized KP I 
equation with critical exponent $p=4/3$ for the 
initial data $u_{0}(x,y)=6\partial_{xx}\exp(-4(x^{2}+y^{2}))$ and the 
$L^{2}$-norm of $u_{y}$ normalized to 1 for $t=0$ in dependence of time.
}
\label{kpIcnorm}
\end{center}
\end{figure}

There is no indication of a 
divergence of the solution, and the $L^{\infty}$-norm of $u$ actually 
decreases as can be seen in Fig.~\ref{kpIcnorm}. However the 
$L^{2}$-norm of $u_{y}$ seems to increase without bound in accordance 
with the theory. It appears that no blowup happens for finite times, but 
this can of course not be decided with numerical methods alone. We 
just did not see an indication of a blowup for even longer times. 

%\begin{Acknowledgements}
%vvvv
%\end {Acknowledgements}

\begin{acknow}
We thank M.~Haragus for useful discussions and hints. The second Author thanks L. Molinet and N. Tzvetkov for useful suggestions and for a lengthy and fruitful joint research on KP equations. 
This work has been supported in part by the project FroM-PDE funded by the European
Research Council through the Advanced Investigator Grant Scheme, the Conseil R\'egional de Bourgogne
via a FABER grant, the ANR via the program ANR-09-BLAN-0117-01 and the Wolfgang Pauli Institute in Vienna.
\end{acknow}

\bibliographystyle{amsplain}

\end{document}